\documentclass[11pt,reqno]{amsart}
\usepackage{mathtools}
\usepackage{enumerate}
\usepackage[backend=bibtex,style=numeric]{biblatex}
\usepackage[colorlinks=true,linkcolor=blue,urlcolor=blue]{hyperref}
\usepackage{dsfont}
\usepackage{geometry}
\usepackage[dvips]{graphics}
\usepackage{graphicx}
\usepackage{epsfig}
\usepackage{todonotes}
\usepackage{epstopdf}
\usepackage{amsmath}
\usepackage{amsfonts}
\usepackage{amssymb}
\usepackage{amsthm}
\usepackage{esint}
\usepackage{eurosym}
\usepackage{multicol}
\usepackage{enumitem}
\newcommand{\E}{\mathbb{E}}
\newcommand{\R}{\mathbb{R}}

\newcommand{\N}{\mathbb{N}}
\usepackage{verbatim}
\usepackage{hyperref}
\usepackage{dsfont}
\usepackage{bbm}
\usepackage{dsfont}
\usepackage{bm}
\usepackage{mathrsfs}
\usepackage{relsize}
\usepackage{textcomp,enumitem}
\usepackage[leqno]{amsmath}
\usepackage{xcolor}
\usepackage[displaymath,mathlines]{lineno}
\allowdisplaybreaks
\usepackage{ifthen}
\makeatletter
\makeindex
\makeatother
\newcommand{\BB}{\mathcal B}

\newcommand{\dd}{\,\mathrm d}

\newcommand{\ip}[2]{\left\langle #1,#2\right\rangle}

\def\leq{\leqslant}

\def \E{\mathbb{E}}
\def \N{\mathbb{N}}
\def \R{\mathbb{R}}

\newcommand{\leqnomode}{\tagsleft@true}
\newcommand{\reqnomode}{\tagsleft@false}
\newcommand{\miezz}{\frac{1}{2}}

\newcommand{\eps}{\varepsilon}

\newcommand{\norm}[1]{\ensuremath{\left\| #1 \right\|}}

\newcommand{\norminf}[1]{\ensuremath{\left\Arrowvert #1 \right\Arrowvert_\infty}}

\newcommand{\dm}[1]{\ensuremath{\frac{\delta #1}{\delta m}}}

\newcommand{\daw}[1]{\mathbf{d}_{#1}}

\newcommand{\bo}[1]{\boldsymbol{#1}}

\newcommand{\ds}{\displaystyle}

\newcommand{\red}[1]{\textcolor{red}{#1}}

\newcommand{\norlh}[1]{\norm{#1}_{\mathcal L(H)}}

\newcommand{\nlh}[1]{\|#1\|_{\mathcal L(H)}}

\newcommand{\nch}[1]{\|#1\|_{C(H)}}
\newcommand{\nh}[1]{\|#1\|_{H}}

\newcommand{\EQ}[1]{\E\Big[\sup\limits_{t\in[0,T]}\|#1(t)\|_H^2\Big]}
\newcommand{\cex}[1]{\E\big[#1\mid\mathscr F_t^0\big]}
\newcommand{\xan}{X^{\alpha,n}}

\newcommand{\stkout}[1]{\ifmmode  \text{\sout{\ensuremath{#1}}}\else\sout{#1}\fi}

\newtheorem{theorem}{Theorem}[section]
\newtheorem{assumption}[theorem]{Assumption}
\newtheorem{lemma}[theorem]{Lemma}
\newtheorem{proposition}[theorem]{Proposition}

\newtheorem{definition}[theorem]{Definition}

\makeatletter
\def\th@remark{%
\thm@headfont{\bfseries}
\normalfont 
}
\makeatother

\theoremstyle{remark}
\newtheorem{remark}{Remark}[section]

\def \R{\mathbb{R}}

\definecolor{red}{rgb}{1.0,0.0,0.0}
\def\red#1{{\textcolor{red}{#1}}}
\definecolor{blu}{rgb}{0.0,0.0,1.0}

\usepackage{amsmath,mathtools}
\usepackage{amsmath, amsthm, amssymb, amsfonts, enumerate}

\usepackage[colorlinks=true,linkcolor=blue,urlcolor=blue]{hyperref}
\usepackage{dsfont}
\usepackage{color}
\usepackage{geometry}
\usepackage{todonotes}
\usepackage{epstopdf}
\usepackage[commandnameprefix=always]{changes}
\usepackage{amsmath,mathtools}
\usepackage{changes}
\usepackage{relsize}
\def \R{\mathbb{R}}

\definecolor{red}{rgb}{1.0,0.0,0.0}

\definecolor{blu}{rgb}{0.0,0.0,1.0}

\definecolor{gre}{rgb}{0.03,0.50,0.03}

\definecolor{darkviolet}{rgb}{0.58, 0.0, 0.83}

\def \eps{\varepsilon}

\title[LQ Mean Field Games and Master Equations in Hilbert Spaces]{Linear quadratic Mean Field Games and Master Equations in Hilbert spaces with common noise}
\author[Ghilli]{Daria Ghilli}
\author[Ricciardi]{Michele Ricciardi}

\address{D.~Ghilli: Dipartimento di Scienze Economiche e Aziendali, University of Pavia, Via San Felice al Monastero 5, Pavia, Italy}
\email{\href{mailto: daria.ghilli@unipv.it}{daria.ghilli@unipv.it}}
\address{M.~Ricciardi: Dipartimento di Economia e Finanza, LUISS University of Rome, Viale Romania 32, Roma, Italy.}
\email{\href{mailto: ricciardim@luiss.it}{ricciardim@luiss.it}}

\date{\today}

\numberwithin{equation}{section}

\begin{document}

\begin{abstract}
We study linear-quadratic Mean Field Games with common noise in an infinite-dimensional Hilbert space. The state dynamics are driven by a possibly unbounded linear operator, and the interaction with the population enters through the mean of its conditional distribution, both in the dynamics and in the cost functional. We allow general quadratic costs including linear and non-separable terms. Exploiting the linear-quadratic structure, we characterize mild solutions of the Mean Field Game system and of the associated Master Equation through systems of operator-valued Riccati equations, linear ordinary differential equations, and stochastic differential equations. A central difficulty is the analysis of a possibly non-self-adjoint Riccati equation arising from the coupling with the population mean. Under suitable structural assumptions, we establish global-in-time existence, uniqueness, and uniform a priori estimates for the corresponding mild solutions. The analysis relies on a fixed-point argument, stability estimates, and Yosida approximations of the unbounded generator. We also prove that the solution of the Master Equation, evaluated along the equilibrium flow of conditional distributions, recovers the Mean Field Game value function, which is not obvious due to the mild formulation used for both solutions. Finally, a verification theorem shows that the mild Mean Field Game solution coincides with the value function of the representative agent and identifies the optimal feedback control.
\end{abstract}

\maketitle

\smallskip


\smallskip

{\textbf{MSC2010 subject classification}}:  49L20, 70H20, 93E20,  47D03

\tableofcontents
\section{Introduction} 
Mean Field Games (MFGs) are a mathematical framework introduced by Lasry and Lions to model strategic interactions among a large number of agents (players) who are individually insignificant but collectively influence the system. Each agent aims to optimize their own objective, which depends on their state and the distribution of all other agents' states. As the number of agents grows, the system becomes computationally intractable, but MFGs provide a limiting description where the population's behavior is represented by a distribution (or ``mean field").
The foundations of this theory trace back to 2006, with seminal contributions by Lasry and Lions on the one hand, and Huang, Caines, and Malhamé on the other. Since then, substantial progress has been made in both theoretical and applied aspects of MFGs. Key references for this study include classical works such as \cite{BeFrYa,CarmonaDelarueBook, CarmonaDelarue2}.
A fundamental aspect of MFG theory is rigorously proving the convergence of the Nash equilibria to the solution of the MFG system as the number of players tends to infinity. Under certain assumptions and in specific cases (see \cite{CarmonaDelarueBook}, Part II, Chapter 6, Section 6.1), it has been shown that any solution to a Mean Field Game corresponds to an $\eps$-Nash equilibrium of the associated $N$-player game. Moreover, in some instances  (see \cite{CarmonaDelarue2}, Part II, Chapter 6, Sections 6.2 and 6.3, and  \cite{CarmonaDelarue2}, Chapter 8, Section 8.2), the convergence of Nash equilibria from the finite 
$N$-player game to the solution of the corresponding Mean Field Game has been established.

Later, Lions proved in \cite{cursolions} that the solutions of the MFG system correspond to trajectories of a new infinite-dimensional partial differential equation (PDE) in the space of measures, termed the Master Equation (ME). 
The ME is a central object in MFGs since it  describes the value function of a representative agent as a function of her/his state and the population's distribution. Thus, unlike the coupled Hamilton-Jacobi-Bellman (HJB) and Fokker-Planck equations in classical MFGs, the ME directly encodes the equilibrium for all possible initial conditions, making it a powerful tool for analyzing the limiting behavior of large systems. Its usefulness for convergence lies in its ability to approximate the Nash equilibria of finite-player games as the number of players tends to infinity, bypassing the need to solve increasingly complex systems.

This paper builds upon the work  \cite{federico2024linearquadratic}, where the authors made the first steps towards the development of a theory for linear-quadratic (LQ) MFGs in infinite-dimensional Hilbert spaces. The main result of \cite{federico2024linearquadratic} is the well-posedness of LQ MFGs in a particular setting: the functional is just quadratic, with interaction through the mean appearing just in the objective functional and without common noise. The present paper extends this framework in three directions, by proving i) the well-posedness of LQ MFGs in a  substantially more general setting, allowing linear terms, non-separable state-control terms, interaction through the mean in both the dynamics and the cost, and common noise; ii)  existence and uniqueness for the Master Equation;  iii)  two verification theorems, the first of them showing that the solution of the Master Equation, evaluated along the equilibrium flow of conditional distributions, recovers the MFG value function; the second proving that the mild MFG solution coincides with the control-theoretic value function and identifies the optimal feedback control.
As already highlighted, we focus on the case where the state dynamics and costs are linear and quadratic, respectively, and moreover, where the interaction among agents depends solely on the mean of the population distribution. This kind of structure has been widely studied in the literature (see \cite{BeFrYa, BeSuYa, CH, federico2024linearquadratic}), since it is suitable to study many problems arising in the applications. In the present paper we do not study these applications, which  will be topics of a forthcoming paper (see \cite{DariaMichele2}).

The primary motivation stems from the challenges posed by infinite-dimensional state spaces, which naturally arise in applications where the standard SDE \eqref{SDE} in a finite dimensional space is not sufficient to describe the complexity of the generic agent's dynamics. This can happen, for example, if the dynamics are delayed or depend on additional variables such as time, age, or spatial position. In the former case, the dynamics are described by a non-Markovian SDE, whereas in the latter case they satisfy a stochastic partial differential equation (SPDE).

We just mention that, unlike finite-dimensional settings, the operators appearing in the dynamics (in particular the operator $A$ acting on the state, see \eqref{SDE}) are very often unbounded in various applications. This unboundedness complicates the analysis, necessitating weaker notions of solutions and sophisticated approximation techniques. In the present paper, we adopt this approach by studying a suitable notion of the so-called $\textit{mild solutions}$ to the MFG and ME. This concept is widely used in infinite-dimensional control.

The literature on infinite-dimensional MFGs, Master Equations, and Nash systems has developed rapidly in recent years, although several fundamental questions remain open. The present paper addresses some of these open issues in the linear-quadratic setting, which have been described above. We also mention that the results of the present paper provide the analytical foundation for the study of the corresponding Nash system, its convergence, the convergence of optimal trajectories and a verification theorem for the Nash system, which will be addressed in \cite{DariaMichele2}. 
We mention \cite{federico2024linearquadratic}  and \cite{FGSMFG} as the first two papers establishing a new theory for infinite-dimensional MFGs: the first in the LQ setting (in a particular case for the objective and dynamics, as explained before in detail) and the second extending the analysis to more
general nonlinear cases under global Lipschitz regularity assumptions on the Hamiltonian. Note also that in the previous papers no convergence was addressed. In \cite{Liu}  the MFG and the NS (but no ME) are studied  in the LQ case in a different setting than ours but only for a small time
horizon (whereas  we prove global-in-time existence). Regarding the setting, in the present paper we deal with more general functionals and with the common noise, which are not treated in \cite{Liu}. On the other hand in \cite{Liu} more general noises than ours are considered. Moreover, an approximation result, showing that the solution to the $N$-player game can be approximated by the solution to the MFG, is proved. \cite{FoZh} addresses the LQ case, in a specific example and in a different setting than ours; there, both the Nash system and the ME are studied.  Finally, in  \cite{CaFoMoSu, FoZh} a specific application to
systemic risk with delay in the control variable is studied, and a verification theorem and a convergence result are also proved.

We now describe more in detail our methods and results. Due to the quadratic structure and the dependence on the population distribution only occurring through its mean, it is possible to search for an explicit solution to the MFG and the ME and rewrite the MFG and the above ME as two systems for the coefficients of the respective ansatzes. For the MFG we get a system of one uncoupled backward Riccati equation, a coupled forward-backward system of SDEs  and a scalar backward SDE. The stochasticity is due to the presence of the common noise. For the ME we get a system of one uncoupled backward Riccati equation, a second backward Riccati equation and  backward linear ODEs. Note that for the Master Equation the mean is treated as an independent state variable rather than as a stochastic process evolving along a prescribed equilibrium flow. Consequently, the equations for the coefficients are deterministic ODEs. The effect of the common noise is instead encoded in the second-order terms with respect to the measure variable in the Master Equation; stochasticity reappears when the deterministic solution of the Master Equation is evaluated along the random conditional distribution generated by the common noise. All the above-mentioned equations are written on the Hilbert space. Exploiting the LQ structure, we introduce a notion of mild solution to the MFG and to the ME through the mild formulations of the corresponding equations for the coefficients.
More in detail, for the MFG system, we search for a solution of the form
\begin{equation}\label{eq:value-fixed-law2}
u(t,x)=\frac12\ip{P(t)x}{x}+\ip{r(t)}{x}+s(t).
\end{equation}
where $P(\cdot), r(\cdot), s(\cdot)$ are time-dependent coefficients which can be found by solving the above-mentioned system of one uncoupled Riccati equation for $ P$ and a forward-backward system of SDEs for $(Y,r)$, where $Y$ denotes the mean of the population's distribution. Note that the function defined in \eqref{eq:value-fixed-law2} is random in the presence of common noise. The above-mentioned equations are completed by the backward SDE for $s$, which can be solved once the other coefficients are known.

The solution of the Riccati equation is standard in the literature, whereas the main  difficulty lies in solving the forward-backward system for $(Y,r)$. This is carried out by a decoupling argument, assuming
$$
r=(K-P)Y+\ell$$
for some deterministic $(K,\ell)$. The main point is that the equation for $K$ is uncoupled and the equation for $\ell$ can be solved once $K$ is found. In order to rigorously justify this decoupling argument, we use the Yosida approximations. This is a standard trick in infinite dimensional control, which allows to replace the unbounded operator $A$ by bounded approximations $A_n$. In this way, the solutions of the approximating equations are strict and it is possible to differentiate. Most of the work then focuses on proving the convergence in $n$ of the solutions of such approximating equations.

We underline that the Riccati equation we get for $K$ encodes a major difficulty, since in the most general case it is not symmetric. We are able to treat this case with an additional technical assumption, that is, Assumption \ref{ipotesi3} case (2-b). We refer to subsection \ref{subsec:lqc} for further comments on such an assumption. We just highlight here that it is a technical generalization in order to treat the non-symmetric case, and the applications we will analyze in \cite{DariaMichele2} will be studied under different and more specific assumptions thanks to the specificity of the setting in each application treated. To the best of our knowledge, global well-posedness results applicable to the non-self-adjoint Riccati equation arising here are not available in the existing infinite-dimensional MFG literature. Anyway, some ideas developed in our results appear in a similar form in the literature of the finite dimensional case, see \cite{Freiling2002} for other details.

Since all our results strongly rely on the study of the described ODEs and SDEs for the coefficients, an important part of this paper is dedicated to a separate study of our equations. First, we study the deterministic case, see subsection \ref{subsec:dc}. In order to establish the global-in-time well-posedness of our ODEs we proceed by first proving it in short time  by the Banach-Caccioppoli fixed point theorem and then by iterating the argument and proving the global-in-time well-posedness result (see Proposition \ref{prop:wellposedness_general}). The iteration is possible just thanks to some a priori estimate on the solution given in Lemma \ref{LaStimaAPrioriPerEccellenza}, which is one of the main results of the paper. Moreover, a useful stability result is contained in Proposition \ref{prop:stabilita}, which is crucial to exploit the convergence of the Yosida approximations mentioned above. Then in subsection \ref{subsec:sc} we study the stochastic case. The proof of the well-posedness of our SDEs is strongly inspired by results already available in the literature, but proved with hypotheses not satisfied in our case. The main difficulty lies in proving estimates on the solutions and in re-adapting the stability results, similar to the one proved in the deterministic case.

Analogously, we proceed for the ME. Note that the well-posedness of the ME can be derived from the above-mentioned results. We just mentioned that its  explicit solution has the following form
$$
U(t,x,y)=\frac{1}{2}\langle P(t)x,x\rangle +\langle \Upsilon(t)y, x\rangle +\frac{1}{2}\langle \Gamma(t)y,y \rangle +\langle\psi(t),x\rangle+\langle\phi(t),y\rangle+\mu(t),
$$
where $x \in H$ is the state variable living in the Hilbert space $H$ endowed with scalar product $\langle \cdot, \cdot \rangle$, $y$ denotes the mean and $P(\cdot), \Upsilon(\cdot), \Gamma(\cdot), \psi(\cdot), \phi(\cdot), \mu(\cdot)$ are time dependent coefficients to be determined by solving the above mentioned system through the above described results.

Moreover, the Yosida approximations are the main technical tool, together with careful estimates, in order to prove the two verification theorems,  that is, Theorem \ref{thm:v1} and Theorem \ref{thm:v2}, in Section \ref{sec:8}.

Finally, we remark that there are other types of solutions to Master equations. We mention e.g. \cite{MoZh}, dealing with \textit{weak solutions} of a quite general ME in finite dimensional spaces with non smooth data. The two notions, mild solution and weak solution, are developed for different purposes and under different technical frameworks. They are complementary approaches to tackling the same fundamental problem, that is, the lack of classical solutions. More in detail, the mild solutions approach is the one typically used for infinite dimensional state equations, which is not the setting of \cite{MoZh}, since $x$ lives in a finite dimensional space, even though the resulting Master equation is written on the Wassertein space. Moreover, it allows us to use the properties of the related semigroup operator in the setting of the present paper. The weak solution approach instead does not focus on the infinite dimensionality of the state dynamics nor on the possibility of writing the state equation with a linear operator. The weak solution approach is inspired by the integration by parts formula and is based on a Lipschitz regularity estimate of the solution with respect to the distribution variable to get global solutions.

The paper is organised as follows. 
In Section \ref{sec:notation} we give the notation which will hold throughout the paper. In Section \ref{sec:formulation} we describe the setting and give the formulation of the MFG and the ME together with the main assumptions. In Section \ref{sec:quadratic} we write the systems which arise from the MFG and the ME by our ansatzes. In Section \ref{sec:preliminari} we give some preliminary results on well-posedness of Riccati equations, ODEs, SDEs, first in the deterministic case and then on in the stochastic case. We also prove stability results which are crucial to use the Yosida approximation. In Section \ref{sec:mfg} we prove the well-posedness of the MFG and in Section \ref{sec:ME} that of the ME.  In Section \ref{sec:8} we prove the two verification theorems.

\section{Notation}\label{sec:notation}
Henceforth, $(H,\langle\cdot,\cdot\rangle_H)$ will always denote a separable real Hilbert space, not necessarily finite-dimensional. Different Hilbert spaces will be denoted by the letters $K,K_0,V$.

We introduce the notation for the following spaces and sets, which will often be used throughout the paper:
\begin{itemize}
\item $\mathcal L(H;K)$ denotes the set of bounded linear operators $\Lambda:H\to K$. The adjoint operator of $\Lambda\in\mathcal L(H;K)$ is denoted by $\Lambda^*\in\mathcal L(K;H)$. Moreover, if $K=H$ we simply write $\mathcal L(H)$ instead of $\mathcal L(H;H)$;
\item $\Sigma(H)$ denotes the set of operators $\Lambda\in\mathcal L(H)$ such that $\Lambda^*=\Lambda$. An operator with such a property is called a \emph{self-adjoint operator};
\item $\Sigma^+(H)$ denotes the subset of $\Sigma(H)$ formed by the \emph{non-negative operators}, i.e. the operators $\Lambda$ such that $\langle\Lambda h,h\rangle\ge0$ for all $h\in H$. We recall that, if $\Lambda\in\Sigma^+(H)$, then
\begin{equation}\label{eq:normsim}
\nlh\Lambda=\sup\limits_{||h||\leq 1}\langle\Lambda h,h\rangle.
\end{equation}
If $\Lambda\in\Sigma^+(H)$, we say that $\Lambda\ge0$. In the same way, we say that $\Lambda_1\ge\Lambda_2$ if $\Lambda_1-\Lambda_2\in\Sigma^+(H)$.

In the same way, we denote by $\mathcal L^+(H)$ the set of operators $\Lambda\in\mathcal L(H)$, not necessarily self-adjoint, such that $\langle\Lambda h,h\rangle\ge0$ for all $h\in H$. Observe that
$$
\Lambda\in\mathcal L^+(H)\iff \Lambda+\Lambda^*\in\Sigma^+(H),
$$
and observe also that the formula \eqref{eq:normsim} is not valid for $\Lambda\in\mathcal L^+(H)$.

With an abuse of notation, we say that $\Lambda\ge0$ also if $\Lambda\in\mathcal L^+(H)$. Observe that the correct notation would be $\Lambda+\Lambda^*\ge0$, but they both imply $\langle\Lambda h,h\rangle\ge0$ for all $h\in H$. Moreover, as before, we say that $\Lambda_1\ge\Lambda_2$ if $\Lambda_1-\Lambda_2\in\mathcal L^+(H)$.
\item $\mathcal L_2(K;H)$ denotes the space of Hilbert-Schmidt operators from a Hilbert space $K$ to another Hilbert space $H$, endowed with the norm
$$
\norm{T}_{\mathcal L_2(K;H)}=\left(\sum\limits_{j=1}^{+\infty}\norm{Te_j}_H^2\right)^\miezz,\qquad T\in\mathcal L_2(K;H)\,.
$$
where $\{e_j\}_{j\ge1}\subset K$ is an orthonormal basis of $K$. Recall that, if $T$ is a Hilbert-Schmidt operator, then $TT^*\in\mathcal L_1(H)$, where $\mathcal L_1(H)$ is the space of trace-class operators.

\item We denote by $L^\infty_s([0,T];\mathcal L(H;K))$ the space of strongly measurable functions
$\Lambda:[0,T]\to\mathcal L(H;K)$, endowed with the norm
\[
\|\Lambda\|_{L^\infty(\mathcal L(H;K))}
:=
\operatorname*{ess\,sup}_{t\in[0,T]}
\|\Lambda(t)\|_{\mathcal L(H;K)}
<+\infty.
\]
Here, strong measurability means that $t\mapsto\Lambda(t)x$ is strongly measurable from $[0,T]$ to $K$ for every $x\in H$. As usual, functions which coincide almost everywhere are identified.
Observe that, if $\Lambda\in L^\infty_s([0,T];\mathcal L(H;K))$, then
$\Lambda^*\in L^\infty_s([0,T];\mathcal L(K;H))$ and
\[
\|\Lambda^*\|_{L^\infty(\mathcal L(K;H))}
=
\|\Lambda\|_{L^\infty(\mathcal L(H;K))}.
\]

In the same way, we define the sets $L^\infty(\Sigma (H))$ and $L^\infty(\Sigma^+(H))$.

\item $C([a,b];E)$, with $a,b\in\R$, $a<b$ and where $E$ is a Banach space, consists of continuous functions $f:[a,b]\to E$. It is endowed with the usual norm $\norm{f}_{\mathcal C([a,b];E)}:=\max\limits_{t\in[a,b]}\norm{f(t)}_E\,.$

When there is no ambiguity for the time interval $[a,b]$, we simply write $\norm{f}_{C(E)}$ instead of $\norm{f}_{C([a,b];E)}$. Moreover, if $f\in C([a,b];\R)$, we adopt the standard notation $\norminf{f}$ instead of $\norm{f}_{C(\R)}$;

\item We denote by $C_s([0,T];\mathcal L(H;K))$ the set of strongly continuous functions from $[0,T]$ into $\mathcal L(H;K)$, i.e. the set of functions $\Lambda:[0,T]\to\mathcal L(H;K)$ such that $t\mapsto\Lambda(t)x$ is a continuous function from $[0,T]$ to $K$ for all $x\in H$. In the same way we define $C_s([0,T];\Sigma(H))$, $C_s([0,T];\Sigma^+(H))$, $C_s([0,T];\mathcal L^+(H))$.
Moreover, we say that $\Lambda_n\to\Lambda$ in $C_s([0,T];\mathcal L(H;K))$ if $\Lambda_n(\cdot)x\to\Lambda(\cdot)x$ in $C([0,T];K)$, for all $x\in H$.
Observe that
\begin{equation}\label{inclusionemanoncredosiaembedding}
C_s([0,T];\mathcal L(H;K))
\subset L^\infty_s([0,T];\mathcal L(H;K)),
\end{equation}
since, for every $x$, $t\to\Lambda(t)x$ is continuous on a compact interval, and hence bounded; the conclusion then follows from the Banach-Steinhaus Theorem. We stress that convergence in $C_s([0,T];\mathcal L(H;K))$ does not mean convergence with respect to the uniform operator norm $\|\cdot\|_{L^\infty(\mathcal L(H;K))}$: indeed, $\Lambda_n\to\Lambda$ in $C_s([0,T];\mathcal L(H;K))$ means that
\[
\sup_{t\in[0,T]}\|(\Lambda_n(t)-\Lambda(t))x\|_K\longrightarrow0
\qquad\text{for every }x\in H.
\]
The norm $\|\cdot\|_{L^\infty(\mathcal L(H;K))}$ will only be used to obtain uniform bounds.

\item $\mathcal P(H)$ denotes the space of Borel probability measures on $H$. Since $H$ can be unbounded (e.g. $H=\R$, in the finite dimensional case), we work with the space $\mathcal P_2(H)$, i.e. the space of Borel probability measures on $H$ with finite second-order moment. This space is endowed with the following distance:
$$
\daw2(m_1,m_2):=\inf\limits_{\pi\in\Pi(m_1,m_2)}\left(\int_{H^2}\norm{x-y}_H^2\,\pi(dx,dy)\right)^\miezz\,,\qquad m_1,m_2\in\mathcal P_2(H)\,,
$$
where the set $\Pi(m_1,m_2)$ consists of Borel probability measures $\pi\in\mathcal P_2(H^2)$ with marginals $m_1$ and $m_2$, i.e. with $\pi(B\times H)=m_1(B)$, $\pi(H\times B)=m_2(B)$ for any Borel set $B\subseteq H$.
\end{itemize}

To avoid ambiguity, we explicitly indicate the Banach space associated with each norm. The only exception is for $\R $, where we adopt the notation $\norm{v}$ to indicate the standard Euclidean norm.

\medskip

Finally, we use the bold notation $\bo w$ to indicate a vector of $H^N$ defined by $\bo w =
(w_1, . . . , w_N)$, where $w_i$ is an already defined vector of $H$. The same notation will be used for vectors of $V^N$. Moreover, we adopt the notation $\bar w_{-i}$ to indicate the average of all the variables but the $i$-th,

$$
\bar w_{-i}=\frac 1{N-1}\sum\limits_{j\neq i}w_j.
$$

\section{Formulation of the problem}\label{sec:formulation}

\subsection{The Mean Field Game}
We describe the asymptotic behavior of an $N$-player differential game, whose trajectories are affected by a $K$-valued idiosyncratic noise $(W(s))_s$, different and independent for any player, and a $K_0$-valued common noise $(W_0(s))_s$, equal for all the players.

Let $(\Omega,\mathscr F,\mathbb P)$ be a probability space, and let $W$ and $W_0$ be two independent cylindrical Wiener processes. We denote by $(\mathscr F_t)_t$ the filtration generated by both $W$ and $W_0$, namely
$$
\mathscr F_t:=\sigma\big(W(s),W_0(s):0\le s\le t\big),
$$
completed in the standard way. Moreover, we denote by $(\mathscr F_t^0)_t$ the filtration generated only by the common noise $W_0$, i.e.
$$
\mathscr F_t^0:=\sigma\big(W_0(s):0\le s\le t\big),
$$
again completed in the standard way.

Consider a final time $T>0$ for the game. The control for the representative player is called $\alpha(\cdot)$, and we require
$$
\alpha\in L^2_{\mathscr P}(\Omega\times[0,T];V),
$$
i.e. $\alpha$ is a square integrable progressively measurable process with values in a Hilbert space $V$, where $\mathscr P$ denotes the progressive $\sigma$-algebra associated with the filtration $(\mathscr F_t)_t$. In the same way, $\mathscr P_0$ will denote the progressive $\sigma$-algebra associated with the filtration $(\mathscr F_t^0)_t$. The norm in these sets is the one induced by $L^2([0,T]\times\Omega;V)$, which will be shortly indicated as $\|\cdot\|_{L^2(V)}$.

The dynamic describing the evolution of the population is called $X(\cdot)$ and solves the $H$-valued SDE in $[0,T]$:
\begin{equation}\label{SDE}
\begin{cases}
dX(s)=\big[AX(s)+B(s)\alpha(s)+B_0(s)\bar m(s)\big]\,ds+\sigma(s) dW_s+\sigma_0(s) dW_s^0\,,\\
\mathscr L(X(0))=m_0\,.
\end{cases}
\end{equation}
Similarly, the dynamic $X^{t,x}(\cdot)$ of the representative player, starting at time $t$ in the position $x$, is given by \eqref{SDE}, with initial condition $X^{t,x}(t)=x$.

Here $m_0\in\mathcal P_2(H)$ is the initial distribution of the population, and $m\in L^2_{\mathscr P_0}([0,T]\times\Omega;\mathcal P_2(H))$ is a random probability measure adapted to the filtration, which at equilibrium represents the conditional law of the population given the common noise $W_0$.

The notation $\bar m(s)$ denotes the (random) average of $m(s)$, i.e. the first order moment:
$$
\bar m(s)=\int_H \xi\, m(s,d\xi)\in L^2_{\mathscr F_s}(\Omega;H).
$$

Finally, $A$, $B$ and $B_0$ are three $H$-valued linear operators. We will be more specific about the assumptions on $A$, $B$, $B_0$ and on the diffusions $\sigma$ and $\sigma_0$ in the following subsections.

We stress the fact that, in many applications, the operator $A$ turns out to be unbounded (see e.g. \cite{federico2024linearquadratic,FoZh,MRMR}), and defined just in a subset $\mathcal D(A)\subset H$. Hence, the equation \eqref{SDE} cannot be understood in a classical sense. To this purpose, we adopt the so-called \emph{mild} formulation, whose definition will be given later.

For a generic player, who starts her/his game in the position $x\in H$ at time $t\in[0,T]$, the (random) value function is defined in the following way. If the cost functional $J$ is given by
\begin{equation}\label{eq:J}
J(t,x,\alpha)=\mathbb E\left[\int_{t}^T\left(L(s,\alpha(s),X^{t,x}(s),\bar m(s))\right)ds+G(X^{t,x}(T),\bar m(T))\mid \mathscr F_t^0\right],
\end{equation}
then
\begin{equation}\label{def:uvaluefunction}
u(t,x)=\inf\limits_{\alpha_\cdot\in L^2_{\mathscr P}}J(t,x,\alpha)\,,
\end{equation}
where $L:[0,T]\times V\times H^2\to\R$ is the Lagrangian cost for the control and $G:H^2\to\R$ is the final cost.

We underline that  $L$ and $G$ depend on the distribution just through $\bar m(s)$. We decide to restrict ourselves to this case inspired by \cite{BeSuYa,federico2024linearquadratic} and by many economic applications (such as vintage capital mean field games models and systemic risk with delay), where this dependence occurs. In the present paper, we do not deal with such applications, which instead will be studied in a forthcoming paper \cite{DariaMichele2}.

In the meantime, observe that in most of the literature the Lagrangian $L$ usually has the following separable structure
$$
L(t,\alpha,x,m)=L(t,\alpha,x)-F(t,x,m),
$$
i.e., there is no coupling between the control of the representative agent and the interaction term with the population. Since in various economic applications that we aim at studying this coupling appears, we decide to keep a non separable (in this regard) Lagrangian.

A Nash equilibrium occurs if the prediction of the players is correct, i.e. if the conditional law of $X^{opt}(s)$ given $W_0$ is $m(s)$, where $X^{opt}$ represents the optimal process, obtained with the optimal control $\alpha^{opt}$. 

We consider the couple $(u,m)$. A standard formal computation, obtained applying Ito's formula and the dynamic programming principle, gives the following stochastic MFG system;

\begin{equation}\label{mfg}
\begin{cases}
\ds du(t)=\Big[-\mathrm{tr}\big((\mathcal A+\mathcal A_0)D^2u\big)+\mathcal H(t,x,Du,\bar m(t))-\mathrm{tr}(\sigma_0\mathrm{Jac}\,v)\Big]dt+\langle v,dW_t^0\rangle\,,\\
\ds dm(t)=\Big[\mathrm{tr}\big((\mathcal A+\mathcal A_0)D^2m\big)+\mathrm{div}(mD_p\mathcal H(t,x,Du,\bar m(t)))\Big]dt-\mathrm{div}(m\sigma_0dW_t^0)\,,\\
\ds m(0)=m_0\,,\qquad u(T,x)=G(x,\bar m(T))\,,
\end{cases}
\end{equation}
where we have defined $\mathcal A:=\frac12\sigma\sigma^*$, $\mathcal A_0:=\frac12\sigma_0\sigma_0^*$, and the \emph{Hamiltonian} of the system $\mathcal H:[0,T]\times H^3\to\R$ as
\begin{equation}\label{def:hamilton}
\mathcal H(t,x,p,y)=\sup_{\alpha\in V} \Big\{-\langle Ax+B(t) \alpha+B_0(t) y, p\rangle - L(t,\alpha,x,y)\Big\},
\end{equation}
Observe that the previous MFG system is stochastic since we have the additional unknown $v(\cdot,\cdot)$, which is a random variable of $(t,x)\in[0,T]\times H$ with values in $K_0$. This term is crucial, since the Hamilton-Jacobi-Bellman SDE for $u$ is backward, and the additional unknown $v$ has to be chosen in order to ensure that the process $u(t)$ remains progressively measurable with respect to $(\mathscr F_t^0)_t$.

In the case without common noise, i.e. when $\sigma_0=0$, the function $m(s)$ is a deterministic measure representing the law of population at time $s$. In this case $v=0$ and we recover the deterministic MFG system on $H$:
\begin{equation}\label{mfgdeterministic}
\begin{cases}
\ds -u_t-\mathrm{tr}(\mathcal AD^2_{xx}u)+\mathcal H(t,x,D_{x}u,\bar m)=0\,,\\
\ds m_t-\mathrm{tr}(\mathcal AD^2_{xx}m)-\mathrm{div}(m\mathcal H_p(t,x,D_xu,\bar m(t)))=0\,,\\
\ds m(0)=m_0\,,\qquad u(T,x)=G(x,\bar m(T))\,.
\end{cases}
\end{equation}
As in the Nash system, the presence of the unbounded operator $A$ in the Hamiltonian $\mathcal H$ forces us to consider just mild solutions for the systems \eqref{mfg} and \eqref{mfgdeterministic}. A particular case of \eqref{mfg} (i.e. for a particular choice of the cost functions) was studied in \cite{federico2024linearquadratic}, where existence and uniqueness in a suitable sense are given. See also \cite{FGSMFG} for the case where the Hamiltonian has a sub-linear growth. We refer to the introduction for further references.

We state the assumptions needed to study the above MFG system \eqref{mfg}.  In particular, the assumptions on $A$ allow us  to consider unbounded operators $A$ and to give a suitable notion of solutions for \eqref{SDE} and \eqref{mfg}.
 Note that, as already said in the introduction, in most of the applications the operator $A$ is unbounded.
\begin{assumption}\label{ipotesi}
Let $K$, $K_0$, $V$ and $H$ be separable real Hilbert spaces. The following hypotheses on the data are required:
\begin{itemize}
\item There exists a closed dense subset $\mathcal D(A)\subseteq H$ s.t. $A:\mathcal D(A)\to H$ is a linear, possibly unbounded, operator. Moreover, $A$ generates a strongly continuous semigroup $e^{tA}:H\to H$, for $t\ge0$;
\item $B\in C([0,T];\mathcal L(V;H))$, $B_0\in C([0,T];\mathcal L(H))$, whereas $\sigma\in C([0,T];\mathcal{L}_2(K;H))$, $\sigma_0\in C([0,T];\mathcal{L}_2(K_0;H))$, $L\in C([0,T]\times V\times H^2;\R)$, $m_0\in\mathcal P_2(H)$.
\end{itemize}
\end{assumption}

 More in detail, because $A$ is unbounded and thanks to the above mentioned assumption, the SDE \eqref{SDE} can have in general just \emph{mild} solutions, whose definition is stated below for a more general case.

\begin{definition}\label{defmildX}
Let $A:D(A)\to H$ satisfy Assumption \ref{ipotesi}. Consider the following SDE:
\begin{equation}\label{SDEfor}
\begin{cases}
dZ(t)=\Big[AZ(t)+B(t)Z(t)+j(t)\Big]dt+\sigma_t^z\,dW_t^0,\\
Z(0)=Z_0,
\end{cases}
\end{equation}
where $B\in L^\infty_{\mathscr P_0}([0,T]\times\Omega;\mathcal L(H))$, $j\in L^2_{\mathscr P_0}([0,T]\times\Omega;H)$, $\sigma^z\in L^2_{\mathscr P_0}([0,T]\times\Omega;\mathcal L_2(K_0;H))$, $Z_0\in L^2_{\mathscr F_0^0}(\Omega;H)$. Then $Z\in L^2_{\mathscr P_0}(\Omega;C([0,T];H))$ is a mild solution of \eqref{SDEfor} if, for all $t\in[0,T]$,
\begin{align*}
Z(t)=e^{tA}Z_0+\int_0^t e^{(t-s)A}\bigl[B(s)Z(s)+j(s)\bigr]\,ds+\int_0^t e^{(t-s)A}\sigma_s^z\,dW_s^0 .
\end{align*}
\end{definition}

In this way, we consider in the SDE \eqref{SDE} a mild solution $X(\cdot)$, in the sense of Definition \ref{defmildX}.\\

As regards the MFG system, a definition of a mild solution is not at all obvious and will be given later, for a particular choice of the cost functions $\mathcal H$ and $G$. 

\subsection{The Master Equation} 
In \cite{CarmonaDelarueBook} (see Part II,
Chapter 6, Section 6.1) it was proved, at least in the finite dimensional case, that the optimal strategies in the MFG system provide approximate $\eps$-Nash equilibria in the $N$-player games.

However, due to the lack of compactness properties of the problem, it is very difficult to prove that the optimal strategies of the Nash system converge to the ones of the MFG system. For this purpose, Lions introduced in \cite{cursolions} a new infinite dimensional equation in the space of measures, whose trajectories are the solutions $(u,m)$ of the MFG system \eqref{mfg}.

This PDE, called \emph{Master Equation}, involves derivatives in the space of measures, which need to be defined in order to give meaning to the equation. Hence, we give the following definition.
\begin{definition}
Let $H$ be a Hilbert space and $U:\mathcal P_2(H)\to\R$. We say that $U$ is $\mathcal C^1$ with respect to the measure variable if $\exists\, K:\mathcal P_2(H)\times H\to\R$ such that, for all $m_1,m_2\in\mathcal P_2(H)$, one has
\begin{equation}\label{def:dmV}
\lim\limits_{s\to0^+}\frac{U(m_1+s(m_2-m_1))-U(m_1)}s=\int_H K(m_1,\xi)\,(m_2-m_1)(d\xi)\,.
\end{equation}
Since $K$ is defined up to an additive constant, we use the notation $\dm{U}$ to indicate the unique $K$ satisfying \eqref{def:dmV} and such that
$$
\int_H K(m,\xi)\,dm(\xi)=0\qquad\forall m\in\mathcal P_2(H)\,.
$$
We also define the function $D_mU:\mathcal P_2(H)\times H\to H$: if $\ds\dm{U}$ is $C^1$ in the space variable, then
$$
D_mU(m,\xi):=D_\xi\dm{U}(m,\xi)\qquad\forall(m,\xi)\in\mathcal P_2(H)\times H\,.
$$
The function $D_mU$ is called the intrinsic derivative of $U$ in the measure variable.
\end{definition}
This definition was already given, in the finite dimensional case, in \cite{MA}.

In the same way, if $D_mU$ is $C^1$ with respect to $m$, we can define $\frac{\delta^2 U}{\delta m^2}:\mathcal P_2(H)\times H^2\to\R$, and, if $\frac{\delta^2 U}{\delta m^2}$ is $C^2$ in the variables $(\xi,\xi')\in H^2$, we define $D^2_{mm}U:\mathcal P_2(H)\times H^2\to\mathcal L(H)$ as
$$ D^2_{mm}U(m,\xi,\xi'):=D^2_{\xi,\xi'}\frac{\delta^2U}{\delta m^2}(m,\xi,\xi')\qquad \forall(m,\xi,\xi')\in\mathcal P_2(H)\times H^2.
$$

The function which will solve the Master Equation can be defined starting from the MFG system. Consider an initial condition $(t_0,m_0)\in [0,T)\times\mathcal P_2(H)$, and define a function $U:[0,T]\times H\times\mathcal P_2(H)\to\R$. For $t=0$, we simply let
$$
U(0,x,m_0)=u(0,x),
$$
where $(u,m)$ is the solution of \eqref{mfg} with $m(0)=m_0$. Observe that, since $u(t)$ is adapted with respect to $\mathscr F_t^0$, then $u(0,\cdot)$ is deterministic, and the definition above is well-posed as a real-valued function.

For a generic time $t_0\in(0,T]$, the ideas are similar. We consider the solution $(u,m)$ of the system \eqref{mfg} in $[t_0,T]$, with initial condition $m(t_0)=m_0$, and with $W_0(t)$ replaced by $W_0(t)-W_0(t_0)$. Then we set, as in \cite[Chapter 5]{MA},
\begin{equation}\label{eq:ME_mfg}
U(t_0,x,m_0)=u(t_0,x)\,.
\end{equation}

Of course, this definition is well-posed if \eqref{mfg} admits a unique mild solution, which will be proved subsequently.

If we formally derive the equation satisfied by $U$, we obtain an infinite dimensional equation in the space of measures, called Master Equation. Its formulation is the following one:

\begin{equation}\label{Master}
\begin{cases}
\displaystyle-U_t-\mathrm{tr}\big((\mathcal A+\mathcal A_0)D^2_{xx}U\big)+\mathcal H\big(t,x,D_xU,\bar m\big)
-\int_H \mbox{tr}\big((\mathcal A+\mathcal A_0)D_\xi D_mU\big)\,dm(\xi)\vspace{0.2cm}\\\vspace{0.2cm}
\ds\quad+\int_H\big\langle D_p\mathcal H\big(t,\xi,D_xU(t,\xi,m),\bar m\big),D_mU\big\rangle\,dm(\xi)-2\int_H\mathrm{tr}(\mathcal A_0D_xD_mU)dm(\xi)\\
\quad\ds-\iint_{H^2}\mathrm{tr}\big(\mathcal A_0D^2_{mm}U(t,x,m,\xi,\xi')\big)dm(\xi)dm(\xi')=0,\\\\
U(T,x,m)=G(x,\bar m)\,.
\end{cases}
\end{equation}
As already said for the MFG system \eqref{mfg}, this equation has to be intended in a mild sense, due to the presence of the unbounded operator $A$ in the Hamiltonian $\mathcal H$.

\subsection{The linear quadratic case}\label{subsec:lqc}
We restrict ourselves to a specific choice of the cost functionals  and we require the following hypotheses on the cost functions $L$ and $G$:
\begin{assumption}\label{ipotesi2}
We take $L:[0,T]\times V\times  H^2\to\R$ and $G:H^2\to\R$ with the following form:
\begin{align}
L(t,\alpha,x,y):&=\frac12\langle R_0(t)\alpha,\alpha\rangle+\frac{1}{2}\langle Q_0(t)x,x\rangle +\miezz\langle Z_0(t)y,y\rangle +\langle S_0(t)y,x\rangle\label{eqn:Fqua}\\
&\quad+\langle L_0(t)\alpha,x\rangle+\langle M_0 (t)\alpha,y\rangle+\langle\theta_0(t),\alpha\rangle
+\langle \eta _0(t),x\rangle +\langle \zeta _0(t),y\rangle + \lambda _0(t)\,,\notag\\
G (x,y):&=\frac{1}{2}\langle Q_T x,x\rangle +\langle S_T y,x\rangle+\miezz\langle Z_T y,y\rangle+\langle \eta_T ,x\rangle+\langle \zeta_T ,y\rangle+\lambda_T \,.\label{eqn:Gqua}
\end{align}
The hypotheses on the coefficients are the following:
\begin{enumerate}
\item $Q_0, Z_0\in C([0,T];\Sigma(H))$, $R_0\in C([0,T];\Sigma(V))$, $Q_T,Z_T\in\Sigma(H)$;
\item $S_0 \in C([0,T];\mathcal L(H))$, $L_0 ,M_0 \in C([0,T];\mathcal L(V;H))$, $S_T\in\mathcal L(H)$;
\item $\eta_0 ,\zeta_0 \in C([0,T];H)$, $\lambda_0 \in C([0,T];\R)$ and $\theta_0 \in C([0,T];V)$;
\item $\eta_T ,\zeta_T \in H$, $\lambda_T \in\R$;
\item~\label{cinque} $\exists\delta_R>0$ s.t. $R_0\ge\delta_R I_V$, where $I_V$ is the identity operator on $V$
\end{enumerate}
\end{assumption}

Observe that the self-adjointness hypotheses on $Q_0,Z_0,R_0,Q_T,Z_T$ are not restrictive, since we can replace $Q_0$ with $\frac12(Q_0+Q_0^*)$, and in the same way with $R_0$, $Z_0$, $Q_T$ and $Z_T$, without changing the Lagrangian cost $L$ and the final cost $G$.

From the quadratic structure of $L$ we can obtain a precise formula for the Hamiltonian $\mathcal H$.

To avoid too heavy notation, we omit the dependence on $t$ for the coefficients. We include it just in the proofs where such a dependence is relevant and should be emphasized.

A straightforward computation gives
\begin{equation}\label{eq:LewisHamilton}
\mathcal H (t,x,p,y)=\frac12\langle\mathcal Bp,p\rangle -\langle (A-D_1)x,p\rangle-\langle D_2y,p\rangle+\langle c_3,p\rangle- F(t,x,y),
\end{equation}
where we adopt the notation
\begin{equation}\label{EQUAZIONI}
\begin{split}
&\mathcal B:=BR_0^{-1}B^*,\qquad\qquad\ \ \,  D_1=BR_0^{-1}L_0^*,\\
& D_2=B_0-BR_0^{-1}M_0^*,\qquad c_3=BR_0^{-1}\theta_0,
\end{split}
\end{equation} and where
\begin{equation*}
\begin{split}
F(t,x,y)=\frac{1}{2}\langle Qx,x\rangle +\miezz\langle Zy,y\rangle +\langle Sy,x\rangle
+\langle \eta,x\rangle +\langle \zeta,y\rangle + \lambda\,,
\end{split}
\end{equation*}
with
\[
\begin{array}{cc}
\begin{aligned}
Q      &= Q_0 - L_0R_0^{-1}L_0^*,\\
Z      &= Z_0 - M_0R_0^{-1}M_0^*,\\
S      &= S_0 - L_0R_0^{-1}M_0^*,
\end{aligned}
&\qquad\qquad
\begin{aligned}
\eta   &= \eta_0 - L_0R_0^{-1}\theta_0,\\
\zeta  &= \zeta_0 - M_0R_0^{-1}\theta_0,\\
\lambda&= \lambda_0 - \frac12\langle R_0^{-1}\theta_0,\theta_0\rangle.
\end{aligned}
\end{array}
\]

As already said, $\mathcal B$ is a short notation for $\mathcal B (t)$, $D_1$ for $D_1 (t)$, $Q$ for $Q (t)$ and so on. Note that the existence of $R_0^{-1}$ is guaranteed by hypothesis (\ref{cinque}) of Assumption \ref{ipotesi2}.

We  note that the optimal $\alpha\in V$ that makes \eqref{eq:LewisHamilton} true is given by
\begin{equation}\label{eq:alfaottimo}
\alpha^{opt}(t,x,p,y)=-R_0^{-1}\big(B^*p+L_0^*x+M_0^*y+\theta_0\big).
\end{equation}

Observe that the hypotheses in Assumption \ref{ipotesi2} strongly generalize the ones in \cite{federico2024linearquadratic}, including a more general mixed term and linear terms, also in the control and hence provide a more general linear-quadratic setting.

Our main result is well-posedness  (as regards mild solutions) of the Mean Field Games and the Master Equation (see the next section for more details). In order to obtain global existence and uniform a priori estimates on the norms of the coefficients, which are crucial for global existence, we need a strengthening of the hypotheses. 

\begin{assumption}\label{ipotesi3}
Assume the following hypotheses:
\begin{enumerate}
\item $Q\in C([0,T];\Sigma^+(H))$ and $Q_T\in\Sigma^+(H)$;
\item Assume that at least one of these two assumptions is satisfied:
\begin{itemize}[label={}, align=left, leftmargin=*]
\item[(2-a, \emph{symmetric case}):] $\exists d_2\in\R$ s.t. $D_2=d_2I_H$, where $I_H$ is the identity operator on $H$. Moreover, we have $S\in C([0,T];\Sigma(H))$ and $S_T\in\Sigma(H)$ and $Q+S, \mathcal B\in C([0,T];\Sigma^+(H))$, $Q_T+S_T\in \Sigma^+(H)$;
\item[(2-b, \emph{non-symmetric case}):] $A$ generates a symmetric group of operators $(e^{tA})_{t\in\R}$ and $Q_T+S_T$ is coercive, i.e. there exists $\delta>0$ such that
$$
\langle (Q_T+S_T)x,x\rangle\ge\delta\|x\|_H^2\,.
$$
Moreover, for a certain $\beta\in\R$ and for all $t\in[0,T]$, it holds $M(t)\ge \beta R$, where the operators $M(t)\in\mathcal L(H^2;H^2)$ and $R\in\mathcal L(H^2;H^2)$ are defined as
$$
M\begin{pmatrix}u\\v\end{pmatrix}=\begin{pmatrix}
Q+S & 0\\
-D_2 & \mathcal B
\end{pmatrix}\begin{pmatrix} u\\v\end{pmatrix},\qquad R\begin{pmatrix}u\\v\end{pmatrix}=\begin{pmatrix} 0 & I_H\\
I_H & 0\end{pmatrix}\begin{pmatrix} u \\ v\end{pmatrix}=(v,u);
$$
\item[(2-c, \emph{trivial case}):] $S=-Q$, $S_T=-Q_T$.
\end{itemize}
\end{enumerate}
\end{assumption}

We give some comments about the previous assumption. First of all, in \cite{federico2024linearquadratic} the operators $S$ and $S_T$ needed to be positive, a condition equivalent to the standard Lasry-Lions monotonicity assumptions, i.e.
\begin{align*}
&\int_H \big(F(t,x,\bar m_1)-F(t,x,\bar m_2)\big)(m_1-m_2)(dx)\ge0\,,\qquad \forall\,(t,m_1,m_2)\in (0,T)\times \mathcal P_2(H)^2\\
&\int_H \big(G(x,\bar m_1)-G(x,\bar m_2)\big)(m_1-m_2)(dx)\ge0\,,\qquad \forall\,(m_1,m_2)\in\mathcal P_2(H)^2\,.
\end{align*}
More recently, several well-posedness results for Mean Field Games have been obtained under the so-called displacement monotonicity conditions; see, e.g., \cite{MeszarosMou2024}. In our notation, the displacement monotonicity condition for the running cost requires that
\[
\begin{aligned}
	&\mathbb E\Big[
	\big\langle D_xL(t,Z^1,X^1,\bar m_1)-D_xL(t,Z^2,X^2,\bar m_2),X^1-X^2\big\rangle\\
	&\qquad\qquad+
	\big\langle D_\alpha L(t,Z^1,X^1,\bar m_1)-D_\alpha L(t,Z^2,X^2,\bar m_2),Z^1-Z^2\big\rangle
	\Big]\ge0,
\end{aligned}
\]
for suitable square-integrable random variables $X^1,X^2$ and $Z^1,Z^2$, with $\mathscr L(X^i)=m_i$, $i=1,2$. Analogously, the displacement monotonicity condition for the terminal cost is
\[
\mathbb E\Big[
\big\langle D_xG(X^1,\bar m_1)-D_xG(X^2,\bar m_2),X^1-X^2\big\rangle
\Big]\ge0.
\]

In the symmetric case, these displacement monotonicity conditions are stronger than the positivity assumptions in Assumption \ref{ipotesi3}-(2-a). Indeed, for the quadratic running cost in Assumption \ref{ipotesi2}, the displacement monotonicity condition is equivalent to
\[
Q\ge0,\qquad
Q+S-\frac14M_0R_0^{-1}M_0^*\ge0,
\]
whereas for the terminal cost it is equivalent to
\[
Q_T\ge0,\qquad Q_T+S_T\ge0.
\]
Since $M_0R_0^{-1}M_0^*\ge0$, these conditions imply in particular
\[
Q\ge0,\qquad Q+S\ge0,\qquad Q_T\ge0,\qquad Q_T+S_T\ge0,
\]
which are precisely the corresponding positivity requirements in Assumption \ref{ipotesi3}-(2-a). Hence, in the symmetric case, our assumptions allow a class of problems which is not necessarily displacement monotone. Notice moreover that, if $M_0=0$, the displacement monotonicity conditions for the running cost reduce exactly to $Q\ge0$ and $Q+S\ge0$.

Observe that we also include in our hypotheses a very special case, namely the condition $S=-Q$ and $S_T=-Q_T$. In this case the well-posedness of the MFG system and the Master Equation can be directly proved, as we will see in Section \ref{sec:mfg} and \ref{sec:ME}, even without any positivity assumptions on the coefficients and the operator $M(\cdot)$. We stress the fact that this particular case is strongly motivated by applications, and this will be shown in our upcoming work \cite{DariaMichele2}.

The non-negativity of $Q$ and $Q_T$ is a standard hypothesis, assumed in  \cite{federico2024linearquadratic} and in most of the papers related to linear quadratic cases to ensure the well-posedness of the Mean Field Games system.\\

To conclude, we give some comments about the second part of the assumption. As regards the hypothesis $M\ge \beta R$, observe that in the symmetric case this hypothesis is always satisfied for $\beta=-\frac{d_2}2$, since it is equivalent to require $\mathcal B(t),Q(t)+S(t)\in\Sigma^+(H)$ for all $t\in[0,T]$. This hypothesis aims to provide a well-posedness result also for non-symmetric Riccati equations of the form \eqref{eq:genericRiccati}, a case almost unexplored in the literature. In particular, we are not aware of any result on non-symmetric Riccati equations in infinite dimension.

The coercivity hypothesis on $Q_T+S_T$, required for the non-symmetric case, can be replaced by the following condition: $Q_T+S_T\ge0$ (i.e. the inequality above with $\delta=0$), and there exists $\delta_0>0$ such that
$$
M(t)\ge\lambda R+\delta_0\begin{pmatrix} I_H & 0\\ 0 & 0\end{pmatrix}.
$$
This condition, in the symmetric case, would be equivalent to requiring $\langle (Q+S)x,x\rangle\ge\delta_0\|x\|_H^2$. We do not delve into the details of this generalization, since it does not add any significant improvement for the applications analyzed in the literature.

Another comment concerns the hypothesis on $A$ required in case $(2-b)$. Although it is always satisfied in the finite dimensional case, in many applications where the infinite dimensional framework is required this hypothesis is not satisfied. It certainly includes, e.g., the case of an operator in $L^2(\R)$ defined by
$$
Af(x)=f'(x),\qquad \mbox{where}\qquad e^{tA}f(\cdot)=f(\cdot+t),\quad t\in\R,
$$
but in many examples, as the delayed systemic risk studied in \cite{DariaMichele2}, the semigroup $e^{tA}$ is not defined for $t<0$. We decide anyway to keep this generic formulation, which allows one to study at least the finite dimensional cases and a class of infinite dimensional ones.\\

\medskip

The previous hypotheses, i.e. Assumptions \ref{ipotesi}, \ref{ipotesi2}, \ref{ipotesi3}, allow us to prove the main result of this paper, i.e. the well-posedness of the Mean Field Games system and the Master Equation.

\section{The linear quadratic case}\label{sec:quadratic}

We now focus on the MFG system \eqref{mfg} and the Master Equation \eqref{Master} when Assumption \ref{ipotesi2} holds. The idea is to look for solutions that are quadratic polynomials, in order to rewrite the systems \eqref{mfg} and \eqref{Master} as smaller systems of equations for the coefficients of the corresponding quadratic ansatz.

\subsection{The Mean Field Games system}

For the MFG system, we search for a solution of \eqref{mfg} of the form
\begin{equation}\label{eq:value-fixed-law}
u(t,x)=\frac12\ip{P(t)x}{x}+\ip{r(t)}{x}+s(t),
\end{equation}
where $P\in C_s([0,T];\Sigma(H))$, $r\in L^2_{\mathscr P_0}([0,T]\times\Omega;H)$, $s\in L^2_{\mathscr P_0}([0,T]\times\Omega;\R).$



In the presence of common noise, the value function is random. 

We now derive formally the equations for \(P\), \(r\), and \(s\).
We assume that $r$ and $s$ admit the semimartingale decompositions
\[
d r(t)=a_t^r\,dt+Z_t^r\,dW_t^0,
\qquad
d s(t)=a_t^s\,dt+\langle Z_t^s,dW_t^0\rangle,
\]
where \(Z^r\in L^2_{\mathscr P_0}(\Omega;\mathcal L_2(K_0;H))\) and \(Z^s\in L^2_{\mathscr P_0}(\Omega; K_0)\). Hence
\[
du(t,x)
=
\frac12\langle P'(t)x,x\rangle\,dt
+\langle a_t^r,x\rangle\,dt
+a_t^s\,dt
+
\left\langle (Z_t^r)^\ast x+Z_t^s,dW_t^0\right\rangle.
\]
Thus the martingale field in the stochastic HJB is
\begin{equation}\label{eq:v}
v(t,x)=(Z_t^r)^\ast x+Z_t^s.
\end{equation}

Moreover,
\[
\nabla_xu(t,x)=P(t)x+ r(t),
\qquad
D_{xx}^2u(t,x)=P(t),\qquad \mathrm{Jac}\,v(t,x)=(Z_t^r)^*.
\]
We adopt the notation $Y_t:=\bar m(t)$, to emphasize the nature of the process $\bar m(t)$. Using the structure of the Hamiltonian, the stochastic HJB equation in \eqref{mfg} reads
\[
\begin{aligned}
&
du(t,x)
+
\Bigg\{
\mathrm{tr}\bigl((\mathcal A(t)+\mathcal A_0(t))P(t)\bigr)
+
\mathrm{tr}\bigl((Z_t^r)^\ast\sigma_0\bigr)-
\frac12
\left\langle
\mathcal B(t)(P(t)x+ r(t)),P(t)x+ r(t)
\right\rangle
\\
&\quad
+
\left\langle P(t)x+ r(t),(A-D_1(t))x+D_2(t)Y_t-c_3(t)\right\rangle
+
\frac12\langle Q(t)x,x\rangle
+
\langle S(t)Y_t,x\rangle
+
\frac12\langle Z(t)Y_t,Y_t\rangle
\\
&\quad
+
\langle \eta(t),x\rangle
+
\langle \zeta(t),Y_t\rangle
+
\lambda(t)
\Bigg\}dt-
\left\langle v(t,x),dW_t^0\right\rangle=0.
\end{aligned}
\]

We now match the terms of equal degree in \(x\). The quadratic terms in \(x\) give
\[
\frac12\langle P'(t)x,x\rangle
+
\left\langle P(t)x,(A-D_1(t))x\right\rangle
-
\frac12\langle \mathcal B(t)P(t)x,P(t)x\rangle
+
\frac12\langle Q(t)x,x\rangle
=0.
\]
Equivalently (as always we do not write the dependence on $t$ for brevity),
\[
P'
+
P(A-D_1)
+
(A^\ast-D_1^\ast)P
-
P\mathcal BP
+
Q
=0.
\]

The linear terms in \(x\) in the drift give
\[
\langle a_t^r,x\rangle
+
\langle  r(t),(A-D_1(t))x\rangle
+
\langle P(t)x,D_2(t)Y_t-c_3(t)\rangle
-
\langle \mathcal B(t)P(t)x, r(t)\rangle
+
\langle S(t)Y_t,x\rangle
+
\langle \eta(t),x\rangle
=0.
\]
Therefore
\[
a^r
=
-(A^\ast-D_1^\ast-P\mathcal B) r-(PD_2+S)Y+Pc_3-\eta.
\]
Hence \(r\) satisfies
\[
dr
=
-
\left[
(A^\ast-D_1^\ast-P\mathcal B) r
+
(PD_2+S)Y
-
Pc_3
+
\eta
\right]dt
+
Z^r\,dW_t^0.
\]

Finally, the zero-order terms give
\[
a^s
+
\mathrm{tr}\bigl((\mathcal A+\mathcal A_0)P\bigr)
+\mathrm{tr}\bigl((Z^r)^\ast\sigma_0\bigr)+
\langle  r,D_2Y-c_3\rangle
-\frac12\langle \mathcal B r, r\rangle
+
\frac12\langle ZY,Y\rangle
+
\langle\zeta,Y\rangle
+
\lambda
=0.
\]
Thus
\[
\begin{aligned}
d s
&=
-
\Big[
\mathrm{tr}\bigl((\mathcal A+\mathcal A_0)P\bigr)
+
\mathrm{tr}\bigl((Z^r)^\ast\sigma_0\bigr)
+
\langle  r,D_2Y-c_3\rangle
\\
&\qquad
-
\frac12\langle \mathcal B r, r\rangle
+
\frac12\langle ZY,Y\rangle
+
\langle\zeta,Y\rangle
+
\lambda
\Big]dt
+
\langle Z^s,dW_t^0\rangle.
\end{aligned}
\]

The terminal conditions for $r$ and $s$ follow by matching the linear terms in \(x\) and the scalar terms in the
terminal condition \(u(T,x)=G(x,Y_T)\). Since
\[
G(x,Y_T)
=
\frac12\langle Q_Tx,x\rangle
+
\langle  S_TY_T,x\rangle
+
\frac12\langle Z_TY_T,Y_T\rangle
+
\langle\eta_T,x\rangle
+
\langle\zeta_T,Y_T\rangle
+
\lambda_T,
\]
we obtain
\[
P(T)=Q_T,\qquad r(T)= S_TY_T+\eta_T,\qquad  s(T)
=
\frac12\langle Z_TY_T,Y_T\rangle
+
\langle\zeta_T,Y_T\rangle
+
\lambda_T.
\]

Given $P$, the optimally controlled state satisfies
\begin{equation}\label{eq:opt-state}
\begin{aligned}
dX_t^*
&=\Bigl[(A-D_1(t))X_t^*+D_2(t)Y_t-
\BB\bigl(P(t)X_t^*+ r(t)\bigr)-c_3(t)\Bigr]\dd t \\
&\quad +\sigma(t)\dd W_t+\sigma_0(t)\dd W_t^0.
\end{aligned}
\end{equation}
Taking conditional expectation with respect to the common noise gives the forward equation
\begin{equation*}
dY_t=\Bigl[(A-D_1(t)+D_2(t))Y_t-\BB(t) P(t)Y_t-\BB(t) r(t)-c_3(t)\Bigr]dt
+\sigma_0(t)\, dW_t^0.
\end{equation*}

Hence, from the MFG system we have obtained the following system for the coefficients $(P,r,s)$, for $Y$ and for the diffusions $(Z_t^r,Z_t^s)$:
\begin{align}
&\begin{cases}\label{eq:P-Riccati}
P'+P(A-D_1)+(A^*-D_1^*)P-P\BB P+Q=0,\\
P(T)=Q_T,
\end{cases}\\
&\begin{cases}\label{eq:r-backward}
dr=-\Bigl[(A^*-D_1^*-P\BB) r+(PD_2+S)Y-Pc_3+\eta\Bigr] dt
+Z^r\, dW_t^0,\\
r(T)= S_TY_T+\eta_T,
\end{cases}\\
&\begin{cases}\label{eq:Y-r-forward}
dY=\Bigl[(A-D_1+D_2)Y-\BB PY-\BB  r-c_3\Bigr] dt
+\sigma_0\, dW_t^0,\\
Y_0=\bar m_0,
\end{cases}\\
&\begin{cases}\label{eq:s}
d s=
-\Big[\mathrm{tr}\bigl((\mathcal A+\mathcal A_0)P\bigr)
+\mathrm{tr}\bigl((Z^r)^\ast\sigma_0\bigr) \\
\qquad+\langle  r,D_2Y-c_3\rangle-\frac12\langle\mathcal B  r, r\rangle+\frac12\langle ZY,Y\rangle+\langle\zeta,Y\rangle+\lambda\Big]dt+\langle Z^s,dW_t^0\rangle,\\
s(T)=\frac12\langle Z_TY_T,Y_T\rangle+\langle\zeta_T,Y_T\rangle+\lambda_T.
\end{cases}
\end{align}

Observe that the MFG system \eqref{mfg} and the system \eqref{eq:P-Riccati}-\eqref{eq:s} are equivalent if $A\in\mathcal L(H)$, since the previous formal computations become rigorous. However, the equations \eqref{eq:P-Riccati}-\eqref{eq:s} can be understood in a mild sense: for the equation \eqref{eq:Y-r-forward} we already have Definition \ref{defmildX}, for \eqref{eq:P-Riccati} and \eqref{eq:r-backward} we state below a general definition.

\begin{definition}\label{defquadric}
Let $A_1:D(A_1)\to H$ and $A_2:D(A_2)\to H$ be (possibly unbounded) operators satisfying the hypotheses of $A$ in Assumption \ref{ipotesi}. Consider the following ODE:
\begin{equation}\label{ODEgen}
\begin{cases}
\Lambda'+\Lambda A_1+A_2\Lambda-\Lambda\mathcal B\Lambda+B_1\Lambda+\Lambda B_2+J=0,\\
\Lambda(T)=J_T,
\end{cases}
\end{equation}
where $J,\mathcal B,B_1,B_2\in L^\infty_s([0,T];\mathcal L(H))$, $J_T\in\mathcal L(H)$. Then $\Lambda\in L^\infty_s([0,T];\mathcal L(H))$ is a mild solution of the ODE \eqref{ODEgen} if, for a.e. $t\in[0,T]$,
\begin{align*}
\Lambda(t)=&\,e^{(T-t)A_2}J_Te^{(T-t)A_1}\\
\quad&+\int_t^T e^{(s-t)A_2}\Big(J(s)-\Lambda(s)\mathcal B(s)\Lambda(s)+B_1(s)\Lambda(s) +\Lambda(s)B_2(s)\Big)e^{(s-t)A_1}\,ds .
\end{align*}
\end{definition}

\begin{remark}\label{rem:finedellafine}
Observe that every mild solution $\Lambda\in L^\infty_s([0,T];\mathcal L(H))$
of \eqref{ODEgen} admits a representative belonging to
$C_s([0,T];\mathcal L(H))$. Indeed, the right-hand side of the mild formulation
in Definition \ref{defquadric}, defined for every $t\in[0,T]$, is strongly
continuous with respect to $t$, by the strong continuity of the semigroups
generated by $A_1$ and $A_2$, and we can easily conclude with the dominated convergence theorem.
\end{remark}

\begin{definition}\label{deflinback}
Let $A:D(A)\to H$ satisfy Assumption \ref{ipotesi}. Consider the following SDEs:
\begin{equation}\label{SDEgen}
\begin{cases}
dz(t)=-\Big[Az(t)+Bz(t)+j\Big]dt-\sigma_t^z\,dW_t^0,\\
z(T)=z_T,
\end{cases}
\end{equation}
where $B\in L^\infty_{\mathscr P_0}([0,T]\times\Omega;\mathcal L(H))$, $j\in L^2_{\mathscr P_0}([0,T]\times\Omega;H)$, $z_T\in L^2_{\mathscr F_T^0}(\Omega;H)$. Then $(z,\sigma^z)\in L^2_{\mathscr P_0}([0,T]\times\Omega;H)\times L^2_{\mathscr P_0}([0,T]\times\Omega;\mathcal L_2(K_0;H))$ is a mild solution of \eqref{SDEgen} if, for all $t\in[0,T]$,
\begin{align*}
z(t)=e^{(T-t)A}z_T+\int_t^T e^{(s-t)A}\bigl[B(s)z(s)+j(s)\bigr]\,ds+\int_t^T e^{(s-t)A}\sigma_s^z\,dW_s^0 .
\end{align*}
\end{definition}

Observe that \eqref{SDEgen} is a backward stochastic differential equation. Therefore, unlike the forward equation \eqref{SDEfor}, the process $\sigma^z$ is part of the unknown. Its role is to guarantee that the solution $z$ is adapted to the filtration $(\mathscr F_t^0)_t$, despite the terminal condition being prescribed at time $T$.

Observe that the previous definition also covers the deterministic case, namely when $\sigma^z=0$ and all the coefficients are deterministic. In this situation, \eqref{SDEgen} reduces to a linear ODE and we simply denote its solution by $z\in C([0,T];H),$ instead of using the stochastic notation $(z,0)$. \\

We are now able to give a suitable definition of mild solution for the linear quadratic MFG system.
\begin{definition}\label{def:LQsolMFG}
Let Assumptions \ref{ipotesi} and \ref{ipotesi2} hold.
A \emph{linear-quadratic-mean mild solution}, briefly an \emph{LQM mild
solution}, of the MFG system \eqref{mfg} is a collection \(\bigl(P,Y,r,Z^r,s,Z^s\bigr)\) solving \eqref{eq:P-Riccati}-\eqref{eq:s} in the sense of Definitions \ref{defmildX}, \ref{defquadric} and \ref{deflinback}.
\end{definition}

Observe that the SDE \eqref{eq:Y-r-forward} for $(Y_t)_t$ is forward, hence it is adapted to the filtration $(\mathscr F_t^0)_t$. Conversely, the SDEs \eqref{eq:r-backward} and \eqref{eq:s} are backward. This means that the unknowns $Z_t^r$ and $Z_t^s$ have to be chosen in order to make the processes $(r(t))_t$ and $(s(t))_t$ adapted to $(\mathscr F_t^0)_t$.

\begin{remark}
	We will also say, with a slight abuse of terminology, that $(u,m)$ is an LQM mild solution of the MFG system if $u$ admits the representation \eqref{eq:value-fixed-law} and $\bar m=Y$, where $Y$ solves \eqref{eq:Y-r-forward}. Notice that this terminology does not determine $m$ uniquely, since the system
	\eqref{eq:P-Riccati}-\eqref{eq:s} determines only its mean $Y$. In particular,
	different measure-valued processes having the same mean $Y$ correspond to the
	same solution at the level of the coefficients.
	
	Anyway, we decide to keep this formulation since it will be convenient in the verification result, where we will compare the Master Equation solution evaluated along a measure-valued process $m(t)$ with the MFG value function.
\end{remark}

\subsection{The Master Equation} In the same way, we look for a solution of \eqref{Master} of the form
\begin{equation}\label{U_explicit}
\begin{split}
& U(t,x,m)=\tilde U\left(t,x,\int_H\xi\,m(d\xi)\right)\,,\qquad\hbox{with $\tilde U:[0,T]\times H^2\to\R$ defined as}\\
& \tilde U(t,x,y)=\frac{1}{2}\langle P(t)x,x\rangle +\langle \Upsilon(t)y, x\rangle +\frac{1}{2}\langle \Gamma(t)y,y \rangle +\langle\psi(t),x\rangle+\langle\phi(t),y\rangle+\mu(t)\,.
\end{split}
\end{equation}
As before, $P,\Gamma:[0,T]\to\Sigma(H)$, $\Upsilon:[0,T]\to\mathcal L(H)$, $\psi,\phi:[0,T]\to H$ are $C^1$ in time, whereas $\mu:[0,T]\to \R$ is a $C^1$ function.

Observe that, for the quadratic term in $x$, we adopt the same notation $P(t)$. This is justified a posteriori, since we will show that the two processes satisfy the same equation, which admits a unique solution.

To shorten the computations, we simply write $y$ instead of $\int_H\xi\,dm(\xi)$. Computing the derivatives, we get
$$
D_xU(t,x,m)=P(t)x+\Upsilon(t)y+\psi(t), \quad D^2_{xx}U(t,x,m)=P(t)\,,
$$
$$
\dm U(t,x,m)(\xi)=\langle \Upsilon^*(t)x+\Gamma(t)y+\phi(t),\xi-y\rangle\,,$$
$$
D_m U(t,x,m,\xi)=D_\xi \dm U(t,x,m,\xi)=\Upsilon^*(t)x+\Gamma(t)y+\phi(t),
$$
$$
D_\xi D_m U(t,x,m,\xi)=D^2_{\xi\xi}\dm U(t,x,m,\xi)=0\,,
$$
$$
D_xD_mU(t,x,m,\xi)=\Upsilon^*(t),\qquad D^2_{mm}U(t,x,m,\xi,\xi')=\Gamma(t).
$$
We apply these formulas in the Master Equation to obtain
\begin{align*}
-&\frac{1}{2}\langle P'(t)x,x\rangle-\langle \Upsilon'(t)y,x\rangle -\frac{1}{2}\langle \Gamma'(t)y,y\rangle -\langle\psi'(t),x\rangle-\langle\phi'(t),y\rangle-\mu'(t)\\
&-\mbox{tr}\big((\mathcal A(t)+\mathcal A_0(t)) P(t)\big)-\langle (A-D_1(t))x+ {D_2(t)y}-c_3(t), P(t)x+\Upsilon(t)y+\psi(t)\rangle\\
&+\frac{1}{2}\langle \mathcal B(t)(P(t)x+\Upsilon(t)y+\psi(t)), P(t)x+\Upsilon(t)y+\psi(t)\rangle\\
&-\int_H\langle \Upsilon^*(t)x+\Gamma(t)y+\phi(t), (A-D_1(t))\xi+ {D_2(t)y}-c_3(t)\rangle\, dm(\xi)\\
&+\int_H\langle \mathcal B(t)(P(t)\xi+\Upsilon(t)y+\psi(t)), \Upsilon^*(t)x+\Gamma(t)y+\phi(t)\rangle dm(\xi)\\
&-2\int_H \mathrm{tr}(\mathcal A_0(t)\Upsilon^*(t))dm(\xi)-\iint_{H^2}\mathrm{tr}(\mathcal A_0(t)\Gamma(t))dm(\xi)dm(\xi')\\
&=\frac{1}{2}\langle Q(t)x,x\rangle +\langle S(t)y,x\rangle +\miezz\langle Z(t)y,y\rangle + \langle \eta(t),x\rangle +\langle \zeta(t),y\rangle + \lambda(t)\,,
\end{align*}
that is
\begin{align*}
-&\frac{1}{2}\langle P'(t)x,x\rangle-\langle \Upsilon'(t)y,x\rangle -\frac{1}{2}\langle \Gamma'(t)y,y\rangle -\langle\psi'(t),x\rangle-\langle\phi'(t),y\rangle-\mu'(t)\\
&-\mbox{tr}\big((\mathcal A(t)+\mathcal A_0(t))P(t)\big)-\langle (A-D_1(t))x+ {D_2(t)y}-c_3(t), P(t)x+\Upsilon(t)y+\psi(t)\rangle\\
&+\frac{1}{2}\langle \mathcal B(t)(P(t)x+\Upsilon(t)y+\psi(t)), P(t)x+\Upsilon(t)y+\psi(t)\rangle\\
&-\langle \Upsilon^*(t)x+\Gamma(t)y+\phi(t), (A-D_1(t))y+ {D_2(t)y}-c_3(t)\rangle\\
&+\langle \mathcal B(t)(P(t)y+\Upsilon(t)y+\psi(t)), \Upsilon^*(t)x+\Gamma(t)y+\phi(t)\rangle-\mathrm{tr}\big(\mathcal A_0(t)(2\Upsilon^*(t)+\Gamma(t))\big)\\
&=\frac{1}{2}\langle Q(t)x,x\rangle +\langle S(t)y,x\rangle +\miezz\langle Z(t)y,y\rangle + \langle \eta(t),x\rangle +\langle \zeta(t),y\rangle + \lambda(t)\,,
\end{align*}
with terminal condition
$$
U(T,x,m)=\frac{1}{2}\langle Q_Tx,x\rangle +\langle S_Ty,x\rangle+\miezz\langle Z_Ty,y\rangle+\langle \eta_T,x\rangle+\langle \zeta_T,y\rangle+\lambda_T\,.
$$
Comparing coefficients, we obtain (omitting the dependence on $t$)
\begin{align}
&\begin{cases}
P'+ P (A-D_{1})+(A^*-D_{1}^*) P - P \mathcal B  P +Q =0,\\
P (T)=Q_T\,;
\end{cases}\label{Riccati:P}\\
&\begin{cases}
\Upsilon'+ \Upsilon (A- D_1-\mathcal B  P )+\,(A^*- D_1^*- P \mathcal B ) \Upsilon-\Upsilon  \mathcal B  \Upsilon+\big(P +\Upsilon\big)D_2 +S =0\,,\\
\Upsilon(T)=S_T\,;
\end{cases}\label{Riccati:Upsilon}\\
&\begin{cases}
\Gamma'+\Gamma\big(A- D_1-\mathcal B (P+\Upsilon)\big)&\hspace{-0.3cm}+\,\big(A^*- D_1^*-(P+\Upsilon^*)\mathcal B \big)\Gamma-\Upsilon^*\mathcal B \Upsilon\\&\hspace{-0.3cm}+\, {(\Upsilon^*+\Gamma)D_{2}+D_{2}^*(\Gamma+\Upsilon)}+Z =0\,,\\
\Gamma(T)=Z_T\,;
\end{cases}\label{Equation:Gamma}\\
&\begin{cases}
\psi'+\big(A^*-D_{1}^*-( P + \Upsilon) \mathcal B\big) \psi -(P+\Upsilon) c_{3} +\eta=0\,,\\
\psi  (T)=\eta_T\,;
\end{cases}\label{Equation:psi}\\
&\begin{cases}
\phi'+\left(A^*-D_{1}^*-\left( P +\Upsilon^* \right) \mathcal B\right) \phi &\hspace{-0.3cm}-\,( \Upsilon^* + \Gamma )\mathcal B \psi  + {D_{2}^*(\psi+\phi)}\\&\hspace{-0.3cm}-\,(\Upsilon^*+\Gamma)c_{3}+\zeta=0\,,\\
\phi (T)=\zeta_T\,;
\end{cases}\label{Equation:phi}\\
&\begin{cases}
\mu'+\mbox{tr}((\mathcal A+\mathcal A_0)P)&\hspace{-0.3cm}+\,\mathrm{tr}\big(\mathcal A_0(2\Upsilon^*+\Gamma)\big)-\frac12\langle \mathcal B \psi  ,  \psi +2 \phi \rangle\\&\hspace{-0.3cm}-\,\langle c_{3},\psi+\phi\rangle+\lambda=0\,,\\
\mu (T)=\lambda_T\,,
\end{cases}\label{Equation:mu}
\end{align}
As already said, the equations \eqref{Riccati:P} and \eqref{eq:P-Riccati} are the same. Note also that equation \eqref{Equation:mu} is immediately solvable by integration:
\begin{equation}\label{explicitmu}
\begin{split}
\mu(t)&=\lambda_T+\int_t^T\mathrm{tr}\Big((\mathcal A(s)+\mathcal A_0(s))P(s)+\mathcal A_0(s)(2\Upsilon^*(s)+\Gamma(s))\Big)ds\\
&\quad+\int_t^T\Big(\lambda(s)-\frac12\langle\mathcal B(s)\psi(s),\psi(s)+2\phi(s)\rangle-\langle c_3(s),\psi(s)+\phi(s)\rangle\Big)ds
\end{split}
\end{equation}

Observe that the system \eqref{Riccati:P}-\eqref{Equation:phi} is uncoupled: the equation \eqref{Riccati:P} (and the same equation \eqref{eq:P-Riccati}) does not depend on the other coefficients and can be solved separately. Once it is solved, the equation \eqref{Riccati:Upsilon} for $\Upsilon$ depends only on $P$, and can be solved separately after \eqref{Riccati:P}. In this way we can obtain the solutions for all the other equations. This is not the case of \eqref{eq:r-backward}-\eqref{eq:Y-r-forward}, where the equations are strongly coupled.\\

As before, the equations \eqref{Riccati:P}-\eqref{Equation:phi} can be understood in a mild sense, and this allows us to give a definition of mild solution for the Master Equation.

\begin{definition}\label{def:LQMMaster}
Let Assumptions \ref{ipotesi} and \ref{ipotesi2} hold. We say that $U:[0,T]\times H\times\mathcal P_2(H)\to\R$ is a \emph{linear-quadratic-mean mild solution}, briefly an \emph{LQM mild
solution}, of the ME equation \eqref{Master}, if $U$ can be written in the form \eqref{U_explicit}, where
\begin{align*}
&(P,\Gamma)\in C_s([0,T];\Sigma(H)), &\Upsilon\in C_s([0,T];\mathcal L(H)),\\
&(\psi,\phi)\in C([0,T];H), &\mu\in C([0,T];\mathbb R)
\end{align*}
solve in a mild sense the system \eqref{Riccati:P}-\eqref{Equation:mu}, in the sense of Definitions \ref{defquadric} and \ref{deflinback}. 
\end{definition}

\begin{remark}\label{rem:casoregolare}
Observe that, if the $6-ple$ $(P,\Upsilon,\Gamma,\psi,\phi,\mu)$ is a classical solution of  \eqref{Riccati:P}-\eqref{Equation:mu}, then the function $U$ defined in \eqref{U_explicit} is smooth and is a classical solution of \eqref{Master}. The same happens with the $6-tuple$ $(P,Y,r,Z^r,s,Z^s)$ for the MFG system \eqref{mfg}.

This happens if $A:H\to H$ is a bounded linear operator, but, in general this is not true (see e.g. the works \cite{BaGo98, FaGo10,federico2024linearquadratic,GM,gozzi2024optimal,MRMR}).

\end{remark}

\section{Preliminary estimates}\label{sec:preliminari}
In this section we collect some estimates and preliminary results which will be useful throughout the rest of the paper.

First of all, thanks to \cite[Part II, Ch.\,1, Cor.\,2.1]{DPB}, we can take $M\ge 1$ and $\omega\in \R^+$ such that
\begin{equation}\label{est:semi}
\|e^{tA}\|_{\mathcal{L}(H)}  +  \|e^{tA^*}\|_{\mathcal{L}(H)}\le Me^{\omega t}\qquad\forall\,t\ge0.
\end{equation}
To simplify the notations, we call $M_T:=Me^{\omega T}$. Moreover, if $A$ generates a group, we call $M_T^-$ the constant $M_T$ associated to $-A$. 

Throughout the proofs, we consider $C$ as a positive constant who may change from line to line, unless specified differently.

The main results of this paper, i.e. the well-posedness of the MFG system \eqref{mfg} and of the Master Equation \eqref{Master}, strongly rely on various estimates and existence results of different Riccati equations and linear ODEs/SDEs which have a similar form. We therefore state in this section some results on general equations which will be useful and simplify the rest of the paper. We distinguish between the stochastic case and the deterministic one, starting with the latter.

\subsection{The deterministic case}\label{subsec:dc} As seen in Section \ref{sec:quadratic}, some equations appearing for the Mean Field Games and the Master Equation are deterministic. Hence, we aim to study here the following generic deterministic Riccati equation in $[t_0,T]$, with $t_0\in [0,T)$:
\begin{equation}\label{eq:genericRiccati}
\begin{cases}
\Lambda'+\Lambda A+A^*\Lambda -\Lambda\mathcal B\Lambda+B_1\Lambda+\Lambda B_2+J=0,\\
\Lambda(T)=J_T,
\end{cases}
\end{equation}
where $\Lambda\in C_s([t_0,T];\mathcal L(H))$.

The well-posedness result strongly needs the following Lemma, which gives a priori estimates for solutions of Riccati equations.

\begin{lemma}\label{LaStimaAPrioriPerEccellenza}
Let $t_0\in[0,T)$, and let $\Lambda\in C_s([t_0,T];\mathcal L(H))$ be a mild solution, in the sense of Definition \ref{defquadric}, of \eqref{eq:genericRiccati}. Assume that $J,\mathcal B, B_1, B_2\in L^\infty_s([t_0,T];\mathcal L(H))$, $J_T\in\mathcal L(H)$, $A$ satisfies Assumption \ref{ipotesi}. Then:
\begin{itemize}
\item [(i)] If $J=0$ and $J_T=0$, we have $\Lambda=0$; as a consequence, \eqref{eq:genericRiccati} has at most one solution in the space $C_s([t_0,T];\mathcal L(H))$.
\item [(ii)] If $\mathcal B=0$, then we have for a certain $C_\Lambda>0$ depending on $M_T$, $\|B_1\|_{L^\infty(\mathcal L(H))}$, $\|B_2\|_{L^\infty(\mathcal L(H))}$, $T$,
\begin{equation}\label{eq:stimalin}
\|\Lambda\|_{L^\infty(\mathcal L(H))}\le C_\Lambda\Big(\|J\|_{L^\infty(\mathcal L(H))}+\|J_T\|_{\mathcal L(H)}\Big).
\end{equation}
\item [(iii)] Assume that, for a certain $\lambda\in\R$ and for a.e. $t\in[t_0,T]$, it holds $M(t)\ge \lambda R$, where the operators $M(t)\in\mathcal L(H^2;H^2)$ and $R\in\mathcal L(H^2;H^2)$ are defined as
$$
M\begin{pmatrix}u\\v\end{pmatrix}=\begin{pmatrix}
J & B_1\\
-B_2 & \mathcal B
\end{pmatrix}\begin{pmatrix} u\\v\end{pmatrix},\qquad R\begin{pmatrix}u\\v\end{pmatrix}=\begin{pmatrix} 0 & I_H\\
I_H & 0\end{pmatrix}\begin{pmatrix} u \\ v\end{pmatrix}=\begin{pmatrix} v\\u\end{pmatrix}.
$$
Moreover, assume that at least one of these two assumptions is satisfied:
\begin{itemize}[label={}, align=left, leftmargin=*]
\item[(iii-a, \emph{symmetric case}):] $B_2=B_1^*$, $J,\mathcal B\in L^\infty(\Sigma^+(H))$, $J_T\in\Sigma^+(H)$;
\item[(iii-b, \emph{non-symmetric case}):] $A$ generates a strongly continuous group of operators $(e^{tA})_{t\in\R}$ and $J_T$ is coercive, i.e. there exists $\delta>0$ such that
$$
\langle J_Tx,x\rangle\ge\delta\|x\|_H^2\,.
$$
\end{itemize}
Then it holds, for a costants $C_\Lambda$ depending on all the norms of the coefficients and the parameters (also depending on $\delta$ and $M_T^-$ in the case $(iii-b)$)
\begin{equation}\label{eq:stimaquad}
\|\Lambda\|_{L^\infty(\mathcal L(H))}\le C_\Lambda.
\end{equation}
\end{itemize}
\end{lemma}

\begin{proof} Without loss of generality, assume $t_0=0$.

\emph{Step 1: (i)}

The first part of the proof is straightforward: if $J=0$, $J_T=0$, we have from Definition \ref{defquadric}
$$
\Lambda(t)=\int_t^T e^{(s-t)A^*}(-\Lambda(s)\mathcal B(s)\Lambda(s)+B_1(s)\Lambda(s)+\Lambda(s)B_2(s))e^{(s-t)A}\,ds.
$$
We call $M_\Lambda$ a bound from above for $\|\Lambda\|_{L^\infty(\mathcal L(H))}$. Then from the previous equation we have
$$
\|\Lambda(t)\|_{\mathcal L(H)}\le M_T^2\int_t^T \Big(M_\Lambda\|\mathcal B\|_{L^\infty(\mathcal L(H))}+\|B_1\|_{L^\infty(\mathcal L(H))}+\|B_2\|_{L^\infty(\mathcal L(H))}\Big)\|\Lambda(s)\|_{\mathcal L(H)}ds.
$$
Hence, from Gronwall's inequality, we have $\Lambda(t)=0$ for all $t\in[0,T]$.

From here we prove uniqueness of mild solutions for \eqref{eq:genericRiccati}. Let $\Lambda_1$ and $\Lambda_2$ be two bounded mild solutions of \eqref{eq:genericRiccati}. The the function $\Lambda:=\Lambda_1-\Lambda_2$ solves

$$
\begin{cases}
\Lambda'+\Lambda A+A^*\Lambda -\Lambda_1\mathcal B\Lambda -\Lambda\mathcal B\Lambda_2+B_1\Lambda+\Lambda B_2=0,\\
\Lambda(T)=0,
\end{cases}
$$
which is exactly \eqref{eq:genericRiccati} with $J=0$, $J_T=0$, $\mathcal B=0$, and $B_1$, $B_2$ replaced by $B_1-\Lambda_1\mathcal B$ and $B_2-\mathcal B\Lambda_2$. From the previous result, we have $\Lambda=0$, i.e. $\Lambda_1=\Lambda_2$. This proves uniqueness.\\

\emph{Step 2: (ii)}

If $\mathcal B=0$, the mild formulation for $\Lambda$ gives
$$
\Lambda(t)=e^{(T-t)A^*}J_T e^{(T-t)A}+\int_t^T e^{(s-t)A^*}\Big(J(s)+B_1(s)\Lambda(s)+\Lambda(s)B_2(s)\Big)e^{(s-t)A}\,ds.
$$
Arguing as before, we get
\begin{align*}
\nlh{\Lambda(t)}&\,\le M_T^2\big(\nlh{J_T}+T||J||_{L^\infty( \mathcal{L}(H))}\big)\\
&\quad+M_T^2\int_t^T\big(||B_1||_{L^\infty( \mathcal{L}(H))}+||B_2||_{L^\infty( \mathcal{L}(H))}\big)\nlh{\Lambda(s)}\,ds\,,
\end{align*}
which, using Gronwall's inequality, gives for a constant $C_\Lambda>0$
$$
\nlh{\Lambda(t)}\le C_\Lambda\big(\nlh{J_T}+||J||_{L^\infty( \mathcal{L}(H))}\big)\,.
$$
Taking the supremum over $t\in[t_0,T]$, we get \eqref{eq:stimalin}.\\

\emph{Step 3: (iii)}

We start with the symmetric case $(iii-a)$, where we require $B_2^*=B_1$, $J,\mathcal B\in L^\infty([0,T];\Sigma^+(H))$, $J_T\in\Sigma^+(H)$. In this case the ODE becomes
$$
\Lambda'+\Lambda (A+B_2)+(A+B_2)^*\Lambda-\Lambda\mathcal B\Lambda+J=0.
$$
Since $\Lambda\in C_s([0,T];\mathcal L(H))$, then thanks to \eqref{inclusionemanoncredosiaembedding} $\Lambda\in L^\infty_s([0,T];\mathcal L(H))$.

This implies $\Lambda^*\in L^\infty_s([0,T];\mathcal L(H))$. Observe also that $\Lambda^*$ solves the same Riccati equation of $\Lambda$. Thanks to Remark~\ref{rem:finedellafine}, we have $\Lambda^*\in C_s([0,T];\mathcal L(H))$. Applying step $(i)$ we get $\Lambda=\Lambda^*$, i.e. $\Lambda\in C_s([0,T];\Sigma(H))$.

Since $J$ and $J_T$ are non-negative and self-adjoint, from the integrated formulation given in \emph{Theorem 4.1} of \cite{CP74} we get $\Lambda\in\Sigma^+(H)$, which implies
$$
\norm{\Lambda}_{L^\infty(\mathcal L(H))}=\sup\limits_{t\in[t_0,T]}\sup\limits_{\norm{x}_H\le 1}\langle\Lambda(t)x,x\rangle.
$$

Observe that, in order to apply the previous theorem, we need that $A+B_2$ generates an evolution operator $\mathscr U(t,s)$, strongly continuous in $(t,s)$. This is guaranteed from \emph{Theorem A.3} of \cite{CP74}, since $A$ generates a strongly continuous semigroup for hypothesis and $B_2\in L_s^\infty([0,T];\mathcal L(H))$.

Now we estimate the scalar product from the mild formulation of \eqref{eq:genericRiccati}, given by Definition~\ref{defquadric}. We have
\begin{align*}
\langle\Lambda(t)x,x\rangle&=\langle J_T e^{(T-t)A}x,e^{(T-t)Ax}\rangle+\int_t^T \langle J(s)e^{(s-t)A}x,e^{(s-t)A}x\rangle ds\\
&\quad+\int_t^T\langle (B_1(s)\Lambda(s)+\Lambda(s)B_2(s))e^{(s-t)A}x,e^{(s-t)Ax}\rangle ds\\
&\quad-\int_t^T\langle \mathcal B(s)\Lambda(s)e^{(s-t)A}x,\Lambda(s)e^{(s-t)A}x\rangle ds\\
&\le C\Big(\norm{J_T}_{\mathcal L(H)}+\norm{J}_{L^\infty(\mathcal L(H))}+\int_t^T\norm{\Lambda(s)}_{\mathcal L(H)}ds\Big),
\end{align*}
where we have used the boundedness of $B_1$, $B_2$, $e^{tA}$ and the non-negativity of $\mathcal B$, which makes the last integral non-positive. Passing to the sup for $\norm x_H\le 1$, we get
$$
\norm{\Lambda(t)}_{\mathcal L(H)}\le C\Big(\norm{J_T}_{\mathcal L(H)}+\norm{J}_{L^\infty(\mathcal L(H))}\Big)+C\int_t^T\norm{\Lambda(s)}_{\mathcal L(H)}ds.
$$
Applying Gronwall inequality, we get \eqref{eq:stimaquad} and we conclude the first case.\\

For the second case, we argue in a different way. 

Consider $X\in C_s([0,T];\mathcal L(H))$ as the mild solution of the following problem
$$
\begin{cases}
	X'-(A+B_2-\mathcal B\Lambda)X=0,\\
	X(T)=I_H.
\end{cases}
$$
Set $C_\Lambda(t):=B_2(t)-\mathcal B(t)\Lambda(t)$. Since $B_2,\mathcal B,\Lambda\in L^\infty_s([0,T];\mathcal L(H))$, we have $C_\Lambda\in L^\infty_s([0,T];\mathcal L(H))$.

In order to write the solution in terms of an evolution family, we perform a time reversal. Let
$$
\widehat X(s):=X(T-s),\qquad s\in[0,T].
$$
Then $\widehat X$ solves
$$
\begin{cases}
	\widehat X'(s)
	=-\big(A+C_\Lambda(T-s)\big)\widehat X(s),\\
	\widehat X(0)=I_H.
\end{cases}
$$
Thanks to the hypotheses, $-A$ generates a strongly continuous semigroup. Moreover,
$s\mapsto-C_\Lambda(T-s)$ belongs to $L^\infty_s([0,T];\mathcal L(H))$. Hence the standard bounded non-autonomous perturbation result yields a strongly continuous evolution family $\mathscr U(s,r)$ associated with
$-A-C_\Lambda(T-s)$ (see also Theorem A.3 of \cite{CP74}). Therefore
$$
\widehat X(s)x=\mathscr U(s,0)x,
$$
and consequently
$$
X(t)x=\mathscr U(T-t,0)x,\qquad x\in H.
$$

We immediately have that $X$ is invertible, since the inverse $Z=X^{-1}$ is given by
$$
\begin{cases}
	Z'+Z(A+B_2-\mathcal B\Lambda)=0,\\
	Z(T)=I.
\end{cases}
$$
To prove that, we shortly call $\tilde A=A+B_2-\mathcal B\Lambda$. If $A$ is bounded, then we can write
$$
\begin{cases}
	(ZX)'=Z'X+ZX'=-Z\tilde AX+Z\tilde AX=0,\\
	Z(T)X(T)=I_H
\end{cases},\qquad
\begin{cases}
	(XZ)'=X'Z+XZ'=\tilde AXZ-XZ\tilde A,\\
	X(T)Z(T)=I_H.
\end{cases}
$$
In the first ODE, we immediately have $ZX=I_H$. In the second one, the constant operator $I_H$ solves the ODE, which has a unique solution. This means $XZ=I_H$. Hence, $Z=X^{-1}$ and we can also write
$$
X^{-1}(t)x=\mathscr U(T-t,0)^{-1}x.
$$

The general case follows approximating $A$ with the Hille-Yosida approximation $(A_n)_n$, and using the strong convergence of the semigroup of operators, see \cite[Proposition 3.6, page 136]{DPB}.

Define now $Y=\Lambda X$. To compute the equation satisfied by $Y$, we assume again $A\in\mathcal L(H)$, and the general case follows by Yosida approximation $(A_n)_n$. We have in this case
\begin{align*}
Y'&=\Lambda'X+\Lambda X'\\
&=-\Lambda AX-A^*\Lambda X+\Lambda\mathcal B\Lambda X-B_1\Lambda X-\Lambda B_2X-JX+\Lambda AX+\Lambda B_2X-\Lambda\mathcal B\Lambda X\\
&=-A^*Y-B_1Y-JX.
\end{align*}

Hence, the couple $(X,Y)$ is a mild solution of the following system
\begin{equation}\label{sistema}
\begin{cases}
X'-(A+B_2)X+\mathcal BY=0,\\
Y'+JX+(A^*+B_1)Y=0,\\
X(T)=I_H,\qquad Y(T)=J_T,
\end{cases}
\end{equation}
or, equivalently, the ODE in $H^2$
$$
\begin{cases}
\bo X'+\bo A\bo X +\bo B\bo X=\bo 0,\\
\bo X(T)=\begin{pmatrix} I_H\\ J_T\end{pmatrix},
\end{cases}
$$
where
$$
\bo X=\begin{pmatrix} X\\Y\end{pmatrix},\qquad \bo A=\begin{pmatrix}
-A & 0\\
0 & A^*
\end{pmatrix},\qquad \bo B=\begin{pmatrix} -B_2 & \mathcal B\\ J & B_1\end{pmatrix}.
$$
Since $A$ generates a strongly continuous group of operators, the operator $\bo A$ generates the group
\begin{equation}\label{semigruppodoppio}
e^{t\bo A}=\begin{pmatrix}
e^{-tA} & 0\\ 0 & e^{tA^*}
\end{pmatrix}.
\end{equation}
Moreover $\bo B\in L^\infty_s([0,T], \mathcal L(H^2))$. Hence, as before, $\bo A+\bo B$ generates an evolution operator $\bo{\mathcal U}(t,s)$, which implies as before the existence of the solution $\bo X$. Moreover
\begin{equation}\label{estbox}
\|\bo X\|_{L^\infty([0,T];\mathcal L(H, H^2))}\le C(\|I_H\|_{\mathcal L(H)}+\norm{J_T}_{\mathcal L(H)})\le C(1+\norm{J_T}_{\mathcal L(H)}).
\end{equation}
We want to obtain the a priori estimate \eqref{eq:stimaquad} for $\Lambda=YX^{-1}$. The estimate for $Y$ is already given in \eqref{estbox}.

For the estimate of $X^{-1}$, we cannot use the equation of $Z$, since it depends on $\Lambda$. Hence, we argue in the following way.

Consider the Yosida approximation $A_n$ of $A$ and the related solution $(X_n,Y_n)$ of \eqref{sistema}. Consider, for $x\in H$, the functions
$$
u_n,v_n:[0,T]\to H,\qquad u_n(t)=X_n(t)x,\qquad v_n(t)=Y_n(t)x.
$$
Of course, the couple $(u_n,v_n)$ solves
$$
\begin{cases}
u'_n-(A_n+B_2)u_n+\mathcal Bv_n=0,\\
v'_n+Ju_n+(A^*_n+B_1)v_n=0,\\
u_n(T)=x,\qquad v_n(T)=J_Tx.
\end{cases}
$$

Let the function $\phi_n:[t_0,T]\to\R$ be defined as
$$
\phi_n(t)=\langle u_n(t),v_n(t)\rangle.
$$
Computing the derivative, we get
\begin{align*}
\phi'_n&=\langle u'_n,v_n\rangle+\langle u_n,v'_n\rangle\\
&=\langle A_nu_n,v_n\rangle+\langle B_2u_n,v_n\rangle-\langle\mathcal Bv_n,v_n\rangle-\langle Ju_n,u_n\rangle-\langle A^*_nv_n,u_n\rangle-\langle B_1v_n,u_n\rangle\\
&=-\left\langle M\begin{pmatrix}u_n\\v_n\end{pmatrix},\begin{pmatrix}u_n\\v_n\end{pmatrix}\right\rangle\le-\lambda\left\langle R\begin{pmatrix}u_n\\v_n\end{pmatrix},\begin{pmatrix}u_n\\v_n\end{pmatrix}\right\rangle=-2\lambda\langle u_n,v_n\rangle=-2\lambda\phi_n.
\end{align*}
Since $\phi'_n+2\lambda\phi_n\le 0$, the function $e^{2\lambda t}\phi_n(t)$ is decreasing, which implies
$$
\langle X_n(t)x,Y_n(t)x\rangle=\phi(t)\ge e^{2\lambda (T-t)}\phi(T)\ge e^{-2|\lambda|T}\langle J_Tx,x\rangle\ge e^{-2|\lambda|T}\delta\|x\|_H^2.
$$
The left-hand side can be bounded from above:
$$
\langle X_n(t)x,Y_n(t)x\rangle\le C(1+\|J_T\|_{\mathcal L(H)})\|X_n(t)x\|_H\|x\|_H.
$$
For $y\in H$, choose $x=X_n(t)^{-1}y$. This implies, for a certain constant $C>0$,
$$
\|X_n(t)^{-1}y\|_H^2\le C\left(1+\|J_T\|_{\mathcal L(H)}\right)\|y\|_H\|X_n(t)^{-1}y\|_H.
$$
Since $A$ generates a strongly continuous group of operators, we have that $X_n(t)^{-1}y$ converges to $X(t)^{-1}y$, for all $y\in H$. Hence, we take the limit for $n\to+\infty$, and then we pass to the $sup$ for $\|y\|_H\le 1$ and for $t\in[t_0,T]$ to get
$$
\|X^{-1}\|_{L^\infty(\mathcal L(H))}\le C(1+\|J_T\|_{\mathcal L(H)}).
$$

Putting together the last estimate with \eqref{estbox}, we get \eqref{eq:stimaquad} and we conclude.
\end{proof}

The previous a priori estimate allows us to prove the following result.
\begin{proposition}\label{prop:wellposedness_general}
Let $T>0$, and consider the ODE \eqref{eq:genericRiccati}, where $A$ satisfies Assumption \ref{ipotesi} and $J,\mathcal B,B_1,B_2\in L^\infty_s([0,T];\mathcal L(H))$, $J_T\in\mathcal L(H)$. Then:
\begin{itemize}
\item [(i)] If $\mathcal B=0$, there exists a unique solution $\Lambda\in C_s([0,T];\mathcal L(H))$ of \eqref{eq:genericRiccati}, satisfying the bound \eqref{eq:stimalin};
\item [(ii)] If the hypotheses of Lemma \ref{LaStimaAPrioriPerEccellenza}, case $(iii)$, are satisfied, then there exists a unique solution $\Lambda\in C_s([0,T];\mathcal L(H))$ of \eqref{eq:genericRiccati}, satisfying the bound \eqref{eq:stimaquad}.
\end{itemize}
\end{proposition}

\begin{proof}
\emph{Step 1: short time existence}. We start proving that there exists $\tau>0$ such that \eqref{eq:genericRiccati} admits a unique mild solution $\Lambda\in C_s\big([T-\tau,T];\mathcal{L}(H)\big)$.

The idea is to use the Banach-Caccioppoli fixed point theorem. First of all, we consider, for $r>0$, $\tau>0$ which will be chosen later, the ball
\begin{equation}\label{BratRiccati}
\mathscr B_{r,\tau}=\Big\{\Lambda \in C_s([T-\tau,T];\mathcal L(H)): \,
\|\Lambda\|_{L^\infty([T-\tau;T];\mathcal L(H))}\leq r\Big\}\,.
\end{equation}
We endow $\mathscr B_{r,\tau}$ with the distance inherited from $L^{\infty}_s([T-\tau;T];\mathcal L(H))$.

We first observe that $\mathscr B_{r,\tau}$ is complete with respect to such metric.
Indeed, let $(\Lambda_n)_n\subset \mathscr B_{r,\tau}$ be a Cauchy sequence. Since $\mathcal L(H)$ is a Banach space, there exists a map $\Lambda:[T-\tau,T]\to\mathcal L(H)$ such that
\[
\sup_{t\in[T-\tau,T]}
\|\Lambda_n(t)-\Lambda(t)\|_{\mathcal L(H)}
\longrightarrow 0.
\]
We claim that $\Lambda\in C_s([T-\tau,T];\mathcal L(H))$. Indeed, for every $x\in H$,
\[
\sup_{t\in[T-\tau,T]}
\|\Lambda_n(t)x-\Lambda(t)x\|_H
\le
\sup_{t\in[T-\tau,T]}
\|\Lambda_n(t)-\Lambda(t)\|_{\mathcal L(H)}\,\|x\|_H
\longrightarrow 0.
\]
Hence $\Lambda_n(\cdot)x\to\Lambda(\cdot)x$ uniformly in
$C([T-\tau,T];H)$. Since each map
$t\mapsto\Lambda_n(t)x$ is continuous, it follows that
$t\mapsto\Lambda(t)x$ is continuous for every $x\in H$. Therefore $\Lambda\in C_s([T-\tau,T];\mathcal L(H))$.

Moreover, we have $\Lambda_n\in \mathscr B_{r,\tau}\subseteq B_r$ where $B_r$ is the closed ball of radius $r$, centered in the origin, of $L^\infty_s([T-\tau;T];\mathcal L(H))$. This immediately implies, for $n\to+\infty$,
\[
\sup_{t\in[T-\tau,T]}
\|\Lambda(t)\|_{\mathcal L(H)}
\le r.
\]
Thus $\Lambda\in \mathscr B_{r,\tau}$, and therefore
$\mathscr B_{r,\tau}$ is complete.

Consider the map $\mathscr F:\mathscr B_{r,\tau}\to C_s\big([T-\tau,T] ; \mathcal L(H)\big)$, defined as
\begin{equation}\label{defFgRiccati}
\begin{split}
\mathscr F(\Lambda)(t)
=\,& e^{(T-t)A^{*}}J_Te^{(T-t)A}\\
&+\int_{t}^{T} e^{(s-t)A^{*}}
\Big(
-\Lambda(s)\mathcal B(s)\Lambda(s)
+B_1(s)\Lambda(s)
+\Lambda(s)B_2(s)
+J(s)
\Big)
e^{(s-t)A}\,ds\,.
\end{split}
\end{equation}

According to Definition \ref{defquadric}, a mild solution to \eqref{eq:genericRiccati} is a fixed point of $\mathscr F$. To apply Banach-Caccioppoli theorem, we need to show that
\begin{itemize}
\item $\mathscr F(\mathscr B_{r,\tau})\subseteq \mathscr B_{r,\tau}$;
\item $\mathscr F$ is a contraction: namely, for some $0<\delta<1$ and for all $\Lambda_1$, $\Lambda_2\in \mathscr B_{r,\tau}$, it holds
\begin{equation}
\big\|\mathscr F(\Lambda_1)-\mathscr F(\Lambda_2)\big\|_{L^\infty([T-\tau,T];\mathcal L(H))}
\le \delta
\big\|\Lambda_1-\Lambda_2\big\|_{L^\infty([T-\tau,T];\mathcal L(H))}.
\label{eq:condcontrRiccati}
\end{equation}
\end{itemize}

From \eqref{defFgRiccati} we  note that $\mathscr F(\Lambda)\in C_s([T-\tau,T];\mathcal L(H))$ if $\Lambda\in \mathscr B_{r,\tau}$. Indeed, for fixed $x$, the function
$$
t\to \mathcal{F}(\Lambda)(t)x$$
is continuous by strong continuity of the semigroups and dominated convergence.
To prove that $\mathscr F(\mathscr B_{r,\tau})\subseteq \mathscr B_{r,\tau}$, we need to show that
$$
\|\mathscr F(\Lambda)(t)\|_{\mathcal L(H)}\leq r
\qquad\forall\,t\in[T-\tau,T]\,.
$$

Using the boundedness of the semigroup and \eqref{BratRiccati} we get from \eqref{defFgRiccati}
\begin{align*}
\|\mathscr F(\Lambda)(t)\|_{\mathcal L(H)}
\le&
\, M_T^2\|J_T\|_{\mathcal L(H)}
+\tau M_T^2\|J\|_{L^\infty(\mathcal L(H))}
\\
&
+\tau M_T^2
\Big(
\|\mathcal B\|_{L^\infty(\mathcal L(H))}r^2
+
(\|B_1\|_{L^\infty(\mathcal L(H))}
+\|B_2\|_{L^\infty(\mathcal L(H))})r
\Big)\\
\le &\ C_1\Big(1+\norlh{J_T}+\tau+\tau r^2\Big),
\end{align*}
where $C_1\in\R$ is a constant depending only on $M_T$, $||J||_{L^\infty(\mathcal L(H))}$, $||\mathcal B||_{L^\infty(\mathcal L(H))}$, $||B_1||_{L^\infty(\mathcal L(H))}$ and $||B_2||_{L^\infty(\mathcal L(H))}$ and $T$. Hence, we come up with the first necessary condition for $\tau$ and $r$:
\begin{equation}\label{eq:cond1trRiccati}
\begin{split}
C_1\Big(1+\norlh{J_T}+\tau+\tau r^2\Big)\leq r,
\end{split}
\end{equation}

As regards \eqref{eq:condcontrRiccati}, we have (up to changing $C_1$)
\begin{align*}
&
\|\mathscr F(\Lambda_1)(t)-\mathscr F(\Lambda_2)(t)\|_{\mathcal L(H)}
\\
\le&
\,M_T^2
\int_t^T
\Big\|
-\Lambda_1\mathcal B\Lambda_1
+\Lambda_2\mathcal B\Lambda_2
+B_1(\Lambda_1-\Lambda_2)
+(\Lambda_1-\Lambda_2)B_2
\Big\|_{\mathcal L(H)}
\,ds
\\
\le&
\,M_T^2
\int_t^T
\Big(
2r\|\mathcal B\|_{L^\infty(\mathcal L(H))}
+\|B_1\|_{L^\infty(\mathcal L(H))}+\|B_2\|_{L^\infty(\mathcal L(H))} \Big)\|\Lambda_1(s)-\Lambda_2(s)\|_{\mathcal L(H)}
\,ds
\\
\le&\,C_1\tau\big(1+r\big)d(\Lambda_1,\Lambda_2).
\end{align*}

Hence, the condition \eqref{eq:condcontrRiccati} is satisfied if $C_1\tau(1+r)<1$. Putting together the last condition with \eqref{eq:cond1trRiccati}, we find the following system for $\tau$ and $r$:
\begin{equation*}
\begin{cases}
C_1(1+\nlh{J_T}+\tau+\tau r^2)\le r,\\
C_1\tau(1+r)<1.
\end{cases}
\end{equation*}

We can easily see that if $r$ is sufficiently large and $\tau$ is sufficiently small, then the couple $(r,\tau)$ solves the previous system. In particular, we can take
$$
r=\mathcal R(||J_T||_{ \mathcal L(H)}),\qquad\tau=\mathcal T(||J_T||_{\mathcal{L}(H)}),
$$
with
$$
\mathcal R(s)=2C_1(1+s),\qquad \mathcal T(s)=\min\left\{\frac 1{2C_1(1+\mathcal R(s))},\frac{\mathcal R(s)}{2C_1(1+\mathcal R(s)^2)}\right\}.
$$

With this choice, we can apply the Banach-Caccioppoli fixed point theorem in the ball $\mathscr B_{r,\tau}$ and obtain a unique mild solution $\Lambda\in C_s\big([T-\tau,T];\mathcal L(H)\big)$.

\emph{Step 2: iteration}. Consider now the solution $\Lambda\in C_s([T-\tau,T];\mathcal L(H))$ found before. Thanks to Lemma \ref{LaStimaAPrioriPerEccellenza}, we have $\|\Lambda\|_{L^\infty([T-\tau,T];\mathcal{L}(H))}\le C_\Lambda$, where $C_\Lambda$ does not depend on $\tau$.

Then we can iterate the procedure, considering $T-\tau$ as the ending point for equation \eqref{eq:genericRiccati}. Hence, we can take a solution $\Lambda\in C_s([T-(\tau+\tau_1),T-\tau];\mathcal L(H))$. 

Observe that the final datum is $\Lambda(T-\tau)$, whose norm is bounded by $C_\Lambda$. This implies that we can take $\tau_1=\mathcal T(C_\Lambda)$. Then we stick the two solutions and we obtain a solution, that by some abuse of notation we call again $\Lambda$, on $[T-(\tau+\tau_1),T]$. This solution satisfies again the a priori estimate $\|\Lambda\|_{L^\infty([T-(\tau+\tau_1),T];\mathcal L(H))}\le C_\Lambda$ from Lemma \ref{LaStimaAPrioriPerEccellenza}.

From now on we can iterate the procedure with the same step. Actually, if we build a solution in $[T-(\tau+\tau_1+\tau_2),T-(\tau+\tau_1)]$, the final datum at $T-(\tau+\tau_1)$ is $\Lambda(T-(\tau+\tau_1))$, and again its norm is bounded by $C_\Lambda$. This implies that, again, we can take $\tau_2=\tau_1=\mathcal T(C_\Lambda)$.

Proceeding with the same method, we obtain for each $k\in\N$ a solution

$$
\Lambda\in C_s([T-(\tau+k\tau_1),T];\mathcal L(H))\,,\qquad\mbox{with }\|\Lambda\|_{L^\infty([T-(\tau+k\tau_1),T];\mathcal L(H))}\le C_\Lambda,.
$$
Then, for $k\ge\frac{T-\tau}{\tau_1}$, we reach $0$, and we obtain the desired solution in $[0,T]$, satisfying (depending on the case) the bound \eqref{eq:stimalin} or \eqref{eq:stimaquad}. This concludes the proof.
\end{proof}

\begin{remark}\label{rem:adjoint-Cs}
Let $\Lambda\in C_s([0,T];\mathcal L(H))$ be the unique mild
solution of \eqref{eq:genericRiccati} given by Proposition
\ref{prop:wellposedness_general}. Then $\Lambda^*\in C_s([0,T];\mathcal L(H)).$

Indeed, taking the adjoint in the mild formulation of \eqref{eq:genericRiccati}, we get that $\Lambda^*$ satisfies
$$
\begin{cases}
(\Lambda^*)'+A^*\Lambda^*+\Lambda^*A-\Lambda^*\mathcal B^*\Lambda^*+\Lambda^* B_1^*+B_2^*\Lambda^*+J^*=0,\\
\Lambda^*(T)=J_T^*.
\end{cases}
$$
Then, thanks to Remark~\ref{rem:finedellafine}, we have $\Lambda^*\in C_s([0,T];\mathcal L(H))$.
\end{remark}
Throughout the paper we often work with a regular approximating system, obtained with the classical Yosida approximants $(A_n)_n$ of $A$. We therefore need a general stability result for a system of the form \eqref{eq:genericRiccati}.

\begin{proposition}\label{prop:stabilita}
Let $T>0$, and consider, for all $n\in\N$, a sequence of possibly unbounded operators $(A_n)_n$ and $A$, satisfying Assumption \ref{ipotesi} and the estimate \eqref{est:semi} uniformly in $n$, and such that $e^{tA_n}\to e^{tA}$ in $C_s([0,T];\mathcal L(H))$. We require the same hypotheses for the adjoints $(A^*_n)_n$ and $A^*$. Consider also some bounded processes $(\mathcal B_n)_n$, $\mathcal B$, $(B_1^n)_n$, $B_1$, $(B_2^n)_n$, $B_2$, $(J_n)_n$, $J\in C_s([0,T];\mathcal L(H))$, $(J_T^n)_n$, $J_T\in\mathcal L(H)$, such that
$$
\mathcal B_n\to \mathcal B, \quad B_1^n\to B_1, \quad B_2^n\to B_2, \quad J_n\to J\quad\mbox{in }C_s([0,T];\mathcal L(H)),\qquad J_T^n\to J_T\quad\mbox{in }\mathcal L(H).
$$
We call $\Lambda_n$ the unique solution of \eqref{eq:genericRiccati} related to the coefficients $(A_n,\mathcal B_n,B^n_1,B_2^n,J_n,J_T^n)$, and we define $\Lambda$ similarly. Assume that, for a constant $C_0$ not depending on $n$,
$$
||\mathcal{B}_n||_{L^\infty(\mathcal L(H))}+||B_1^n||_{L^\infty(\mathcal L(H))}+||B_2^n||_{L^\infty(\mathcal L(H))}+||J_n||_{L^\infty(\mathcal L(H))}+\norlh{J_T^n}+||\Lambda_n||_{L^\infty(\mathcal L(H))}\le C_0.
$$
Then it holds
$$
\Lambda_n\to\Lambda\qquad\mbox{in }C_s([0,T];\mathcal L(H)).
$$
\end{proposition}

Before proving the result, we observe that the existence of the solutions $\Lambda_n$ and $\Lambda$, and a uniform bound on $\Lambda_n$, are required in the previous proposition. These results are true, for example, if the Assumptions of Lemma \ref{LaStimaAPrioriPerEccellenza} are satisfied, and they derive from the same Lemma and Proposition \ref{prop:wellposedness_general}.

However, we can prove that with our hypotheses the strong convergence implies a uniform convergence in the compact sets of $H$. This is stated in the following Lemma, which plays an important role in the proof of the Proposition.

\begin{lemma}\label{ascoli}
Let $(R_n)_n, R\in C_s([0,T];\mathcal L(H))$, and assume that
$$
\|R_n\|_{L^\infty(\mathcal L(H))}\le C_0,\qquad R_n\to R\quad\mbox{in }C_s([0,T];\mathcal L(H)).
$$
Then the strong convergence is uniform on the compact sets of $H$, i.e. for any $K\subset H$ compact we have
$$
\lim\limits_{n\to+\infty}\sup\limits_{x\in K}\|R_nx-Rx\|_{C(H)}=0.
$$
\end{lemma}
\begin{proof}
The proof is straightforward. Consider the sequence
$$
R_n: K\to C([0,T];H),\qquad x\to R_n(\cdot)x.
$$
Observe that $K$ is compact (hence bounded), and the sequence is equibounded and equicontinuous: if $C_K>0$ is such that $K\subset B_{C_K}(0)$, we have
$$
\|R_n(\cdot)x\|_{C(H)}\le \|R_n\|_{L^\infty(\mathcal L(H))}\|x\|_H\le C_0 C_K,\qquad \|R_n(\cdot)x-R_n(\cdot)y\|_{C(H)}\le C_0\|x-y\|_H.
$$
Moreover, for every fixed $x\in K$, the sequence $\{R_n(\cdot)x\}_{n\in\mathbb N}$ is relatively compact in $C([0,T];H)$, since $R_n(\cdot)x\to R(\cdot)x$ in $C([0,T];H)$ by the convergence $R_n\to R$ in $C_s([0,T];\mathcal L(H))$. Hence, by the Arzel\`a--Ascoli theorem, $\{R_n\}_{n\in\mathbb N}$ is relatively compact in $C(K;C([0,T];H))$, and we have $R_{n_k}\to \tilde R$ in $C(K;C([0,T];H))$ for some subsequence $(R_{n_k})_k$. Of course we must have $\tilde R=R$, and  since the argument can be applied to any subsequence, we get that every subsequence has a subsubsequence converging to R and this implies $R_n\to R$ in $C(K;C([0,T];H))$, which concludes the proof.
\end{proof}

Now we are ready to prove the previous proposition.

\begin{proof}[Proof of Proposition \ref{prop:stabilita}]
To prove the result, we subtract the mild formulation from Definition \ref{defquadric} for $\Lambda_n$ and $\Lambda$. Evaluating them at $x\in H$, we get (omitting for brevity the dependence on $s$)
\begin{align*}
&\Lambda_n(t)x-\Lambda(t)x\\
&=e^{(T-t)A_n^*}J_T^ne^{(T-t)A_n}x-e^{(T-t)A^*}J_Te^{(T-t)A}x+\int_t^T \Big(e^{(s-t)A_n^*}J_ne^{(s-t)A_n}x-e^{(s-t)A^*}Je^{(s-t)A}x\Big)\,ds\\
&\ +\int_t^T\Big[ e^{(s-t)A_n^*}\Big(B_1^n\Lambda_n +\Lambda_nB_2^n-\Lambda_n\mathcal B_n\Lambda_n\Big)e^{(s-t)A_n}-e^{(s-t)A^*}(B_1\Lambda+\Lambda B_2-\Lambda\mathcal B\Lambda)e^{(s-t)A}\Big]x\,ds.
\end{align*}
This implies
$$
\|(\Lambda_n(t)-\Lambda(t))x\|_H\le \alpha_1+\alpha_2+\int_t^T\big[\alpha_3(s)+\alpha_4(s)+\alpha_5(s)\big]ds,
$$
where
\begin{align*}
\alpha_1&=\big\|e^{(T-t)A_n^*}J_T^ne^{(T-t)A_n}x-e^{(T-t)A^*}J_Te^{(T-t)A}x\big\|_{H},\\
\alpha_2&=T\sup_{0\le t\le s\le T}\big\|e^{(s-t)A_n^*}J_ne^{(s-t)A_n}x-e^{(s-t)A^*}Je^{(s-t)A}x\big\|_{H},\\
\alpha_3(s)&=\big\|\big(e^{(s-t)A^*_n}B_1^n\Lambda_n e^{(s-t)A_n}-e^{(s-t)A^*}B_1\Lambda e^{(s-t)A}\big)x\big\|_H,\\
\alpha_4(s)&=\big\|\big(e^{(s-t)A^*_n}\Lambda_n B_2^n e^{(s-t)A_n}-e^{(s-t)A^*}\Lambda B_2 e^{(s-t)A}\big)x\big\|_H,\\
\alpha_5(s)&=\big\|\big(e^{(s-t)A^*_n}\Lambda_n\mathcal B_n\Lambda_n e^{(s-t)A_n}-e^{(s-t)A^*}\Lambda\mathcal B\Lambda e^{(s-t)A}\big)x\big\|_H.
\end{align*}
The terms $\alpha_1$ and $\alpha_2$ are estimated in the same way, as are the terms $\alpha_3$, $\alpha_4$ and $\alpha_5$. Hence, we show here just the bounds from above for $\alpha_1$ and $\alpha_5$. For $\alpha_1$ we get
\begin{align*}
\alpha_1&\le \big\|e^{(T-t)A^*_n}J^n_T \big(e^{(T-t)A_n}-e^{(T-t)A}\big)x\big\|_H+\big\|e^{(T-t)A^*_n}\big(J^n_T-J_T\big)e^{(T-t)A}x\big\|_H\\
&\quad+\big\|\big(e^{(T-t)A^*_n}-e^{(T-t)A^*}\big)J_Te^{(T-t)A}x\big\|_H\\
&\le M_TC_0\| (e^{(T-t)A_n}-e^{(T-t)A})x\|_H+M_T\|(J_T^n-J_T)e^{(T-t)A}x\big\|_H\\
&\quad+\big\|\big(e^{(T-t)A^*_n}-e^{(T-t)A^*}\big)J_Te^{(T-t)A}x\big\|_H\\
&\le \omega_n(t,x),
\end{align*}
where $\omega_n(t,x)\to 0$ for $n\to+\infty$, uniformly in $t\in[0,T]$, and where we used the strong convergence of $J^n_T$, $e^{(T-t)A_n}$ and $e^{(T-t)A}$. Observe that, thanks to Lemma \ref{ascoli}, the convergence of $\omega_n$ towards $0$ is uniform on the compact sets of $H$.

As regards $\alpha_5$, we have in a similar way
\begin{align*}
\alpha_5(s)
&\le
\big\|e^{(s-t)A_n^*}\Lambda_n\mathcal B_n\Lambda_n
\big(e^{(s-t)A_n}-e^{(s-t)A}\big)x\big\|_H+
\big\|e^{(s-t)A_n^*}\Lambda_n\mathcal B_n
(\Lambda_n-\Lambda)e^{(s-t)A}x\big\|_H\\
&\quad+
\big\|e^{(s-t)A_n^*}\Lambda_n(\mathcal B_n-\mathcal B)
\Lambda e^{(s-t)A}x\big\|_H+
\big\|e^{(s-t)A_n^*}(\Lambda_n-\Lambda)
\mathcal B\Lambda e^{(s-t)A}x\big\|_H\\
&\quad+
\big\|\big(e^{(s-t)A_n^*}-e^{(s-t)A^*}\big)\Lambda\mathcal B\Lambda e^{(s-t)A}x\big\|_H\\
&\le M_TC_0^3\|(e^{(s-t)A_n}-e^{(s-t)A})x\|_H+M_TC_0^2\|(\Lambda_n-\Lambda)e^{(s-t)A}x\|_H\\
&\quad+M_TC_0\|(\mathcal B_n-\mathcal B)\Lambda e^{(s-t)A}x\|_H+M_T\|(\Lambda_n-\Lambda)\mathcal B\Lambda e^{(s-t)A}x\|_H\\
&\quad+\big\|\big(e^{(s-t)A_n^*}-e^{(s-t)A^*}\big)\Lambda\mathcal B\Lambda e^{(s-t)A}x\big\|_H\\
&\le \eps_n(s,t,x)
+C\big\|(\Lambda_n-\Lambda)e^{(s-t)A}x\big\|_H
+C\big\|(\Lambda_n-\Lambda)\mathcal B\Lambda e^{(s-t)A}x\big\|_H,
\end{align*}

where $\eps_n(s,t,x)\to 0$ for $n\to+\infty$ and for all $0\le t\le s\le T$, uniformly on the compact sets of $H$. Moreover, $\eps_n$ is uniformly bounded in $n,s,t$ for fixed $x\in H$.

To avoid too heavy notation, we define
\[
\Delta_n:=\Lambda_n-\Lambda.
\]
Hence, arguing in the same way for the terms $\alpha_2$, $\alpha_3$, $\alpha_4$, we get
\begin{equation*}
\begin{split}
\|\Delta_n(t)x\|_H&\le \omega_n(t,x)
+C\int_t^T \eps_n(s,t,x)\,ds+C\int_t^T
\|\Delta_n(s)L_1(s,t)x\|_H\,ds,\\
&\quad+C\int_t^T
\Big(
\|\Delta_n(s)L_2(s,t)x\|_H
+\|\Delta_n(s)L_3(s,t)x\|_H
\Big)\,ds,
\end{split}
\end{equation*}
where
\[
L_1(s,t)=e^{(s-t)A},\qquad
L_2(s,t)=B_2(s)e^{(s-t)A},\qquad
L_3(s,t)=\mathcal B(s)\Lambda(s)e^{(s-t)A},
\]
and $\varepsilon_n$ collects all the terms where $\Delta_n$ does not appear. As said before, $\eps_n\to0$ for all $0\le t\le s\le T$, and it is uniformly bounded in $t$ and on the compact sets of $H$. Hence, up to changing the vanishing term $\omega_n$, we have
\begin{equation}\label{eq:passobase}
\begin{split}
\|\Delta_n(t)x\|_H&\le \omega_n(t,x)+C\int_t^T
\|\Delta_n(s)L_1(s,t)x\|_H\,ds,\\
&\quad+C\int_t^T
\Big(
\|\Delta_n(s)L_2(s,t)x\|_H
+\|\Delta_n(s)L_3(s,t)x\|_H
\Big)\,ds,\\
&\le \omega_n(t,x)+C_1\int_t^T \sup\limits_{i\in\{1,2,3\}}\|\Delta_n(s)L_i(s,t)x\|_H\,ds,
\end{split}
\end{equation}
for some $C_1>0$ not depending on $n$. We want to prove by induction that, for all $m\in\N$, there exists $\omega_{n,m}:[0,T]\times H\to\R^+$ such that
\begin{equation}\label{eq:induzione}
\begin{split}
\|\Delta_n(t)x\|_H&\le \omega_{n,m}(t,x)\\
&\quad+C_1^m\int_{t\le s_1\le \cdots\le s_m\le T}\sup\limits_{\bo i\in\{1,2,3\}^m}\|\Delta_n(s_m)L_{\bo i}(s_1,\cdots,s_m,t)x\|_H\,ds_1\cdots ds_m,
\end{split}
\end{equation}
where we require  $\omega_{n,m}(t,x)\to 0$ for $n\to+\infty$, uniformly in $t$ and on the compact sets of $H$, for every fixed $m$, and with $\omega_{n,m}$ uniformly bounded in $n,t$ and locally in $x$, for fixed $m\in\N$, and where we define recursively
$$
L_{(i_1,\dots,i_m)}(s_1,\dots,s_m,t)=L_{i_m}(s_m,s_{m-1})L_{(i_1,\dots,i_{m-1})}(s_1,\dots,s_{m-1},t).
$$
The case $m=1$ is exactly \eqref{eq:passobase}. Now, assume that \eqref{eq:induzione} is true for $m$. Then, using \eqref{eq:passobase} inside the integral in \eqref{eq:induzione} (here the role of $x$ is played by $L_{\bo i}(s_1,\dots,s_m,t)x$), we get
\begin{gather*}
\|\Delta_n(t)x\|_H\le \omega_{n,m}(t,x)+C_1^m\int_{t\le s_1\le \cdots\le s_m\le T}\sup\limits_{\bo i\in\{1,2,3\}^m}\omega_n\big(s_m,L_{\bo i}(s_1,\dots,s_m,t)x\big)\,ds_1\dots ds_m\\
+\,C_1^{m+1}\int_{t\le s_1\le \cdots\le s_m\le T}\sup\limits_{\bo i\in\{1,2,3\}^m}\int_{s_m}^T\sup\limits_{i\in\{1,2,3\}}\|\Delta_n(s)L_{i}(s,s_m)L_{\bo i}(s_1,\dots,s_m,t)x\|_H\, ds\, ds_1\dots ds_m.
\end{gather*}
The last term is easily bounded from above by
$$
C_1^{m+1}\int_{t\le s_1\le \cdots\le s_{m+1}\le T}\sup\limits_{\bo i\in\{1,2,3\}^{m+1}}\|\Delta_n(s_{m+1})L_{\bo i}(s_1,\cdots,s_m,s_{m+1},t)x\|_H\,ds_1\dots ds_{m+1}.
$$

As regards the first integral term, observe that the integrand converges to $0$ for any fixed $s_1,\dots, s_m$ and it is bounded thanks to the definition of $\omega_n$. Then, using the uniform convergence on the compact set described below, the integral converges to $0$ uniformly with respect to $t$. Moreover, since for any $K\subset H$ compact the set
\[
\bigcup_{\bo i\in\{1,2,3\}^m}\Big\{L_{\bo i}(s_1,\dots,s_m,t)x\mid t\le s_1\le \dots\le s_m\le T,\ x\in K\Big\}
\]
is a finite union of images of compact sets through a continuous function, it is compact. This implies that the integral converges uniformly on the compact sets, and we can absorb this term into $\omega_{n,m+1}$. Hence, we get
\begin{equation*}
\begin{split}
\|\Delta_n(t)x\|_H&\le \omega_{n,m+1}(t,x)\\
&\quad+C_1^{m+1}\int_{t\le s_1\le \cdots\le s_{m+1}\le T}\sup\limits_{\bo i\in\{1,2,3\}^{m+1}}\|\Delta_n(s_{m+1})L_{\bo i}(s_1,\cdots,s_{m+1},t)x\|_H\,ds_1\cdots ds_{m+1},
\end{split}
\end{equation*}
which proves \eqref{eq:induzione}. Now, consider $C_2$ such that $\|L_i(s,t)\|_{\mathcal L(H)}\le C_2$. Using the uniform bound of $\Delta_n$ and the definition of $L_{\bo i}$, we get
$$
\|\Delta_n(t)x\|_H\le \omega_{n,m}(t,x)+2C_0C_1^{m}C_2^m\|x\|_H\mathrm{Vol.}(\{0\le s_1\le \dots\le s_m\le T\}),
$$
i.e. 
$$
\|\Delta_n(t)x\|_H\le \omega_{n,m}(t,x)+\frac{2C_0(C_1C_2T)^m}{m!}\|x\|_H.
$$
For fixed $m\in\N$, we pass to the sup in $t\in[0,T]$ and to the limsup in $n\to+\infty$, using the uniform convergence of $\omega_{n,m}$ in $t$. This implies
$$
\limsup\limits_{n\to+\infty}\|\Lambda_n(\cdot)x-\Lambda(\cdot)x\|_{C([0,T];H)}\le \frac{2C_0(C_1C_2T)^m}{m!}\|x\|_H,
$$
for all $m\in\N$. Letting $m\to+\infty$, we conclude the proof. Since $x\in H$ is arbitrary, this is exactly the convergence in
\(C_s([0,T];\mathcal L(H))\).
\end{proof}

\begin{remark}\label{backlindet}
The computations done in the previous results can be easily extended to the case of an operator $\gamma\in C([t_0,T];H)$, which satisfies in a mild sense the following equation
\begin{equation}\label{eqlin}
\begin{cases}
\gamma'+A\gamma+B\gamma+\ell=0,\\
\gamma(T)=\ell_T,
\end{cases}
\end{equation}
where $A:D(A)\subseteq H$ generates a strongly continuous semigroup $e^{tA}$, $B\in L^\infty_s([0,T];\mathcal L(H))$, $\ell\in C([0,T];H)$, $\ell_T\in H$. The existence and uniqueness result is guaranteed from \cite[Proposition 3.4, page 136]{DPB}. Moreover, for a constant $C_\gamma$ depending on $||B||_{L^\infty(\mathcal L(H))}$, $T$ and $M_T$, we obtain as in Lemma \ref{LaStimaAPrioriPerEccellenza}
\begin{equation}\label{eq:estlin2}
\nch\gamma\le C_\gamma\Big(\nch\ell+\nh{\ell_T}\Big).
\end{equation}
Moreover, the stability result can be proved as in Proposition \ref{prop:stabilita}. If $B_n,B\in C_s([0,T];\mathcal L(H))$, $e^{tA_n}\to e^{tA}$, $B_n\to B$ in $C_s([0,T];\mathcal L(H))$, $\ell_n\to\ell$ in $C([0,T];H)$, $\ell^n_T\to \ell_T$ in $H$ and the norms are uniformly bounded, then $\gamma_n\to\gamma$ in $C([0,T];H)$. We omit the proofs of these results, since they are the same as the previous proofs in $C_s([0,T];\mathcal L(H))$ and since they are provided in a general stochastic case in the next subsection.
\end{remark}

\subsection{The stochastic case}\label{subsec:sc} We need also some basic results for mild solutions of linear SDEs of the form \eqref{SDEfor} (forward) and \eqref{SDEgen} (backward).

As already said, the forward case is simpler, since the diffusion term is not part of the unknown, and it has been extensively studied. We recollect the results in the following proposition, proving also the stability of the SDE in a general framework.

We rewrite here the SDE \eqref{SDEfor} for better readability:
\begin{equation*}
\begin{cases}
dz(t)=\Big[Az(t)+B(t)z(t)+j(t)\Big]dt+\sigma_t dW_t^0,\\
z(0)=z_0.
\end{cases}\tag{\ref{SDEfor}}
\end{equation*}

Henceforth, given a Banach space $V$ and $p\ge1$, we adopt the notations $\mathcal{Q}_{\mathscr P_0}^p$ and $\mathcal{S}_{\mathscr P_0}^p$ to denote the Banach spaces
\begin{equation*}
\begin{split}
&\mathcal{Q}_{\mathscr P_0}^p(V):= L^p_{\mathscr P_0}([0,T]\times\Omega;V)\cap L^p(\Omega;C([0,T];V)),\\
&\mathcal{S}_{\mathscr P_0}^p(V):= L^p_{\mathscr P_0}([0,T]\times\Omega;V)\cap C([0,T];L^p(\Omega;V)),
\end{split}
\end{equation*}
i.e. the spaces of progressively measurable processes $X(\cdot)$, in the first case a.s. continuous in time. They are endowed with the norm (omitting the dependence on $\mathscr P_0$ to avoid too heavy notation)
\begin{align*}
\|X\|_{\mathcal{Q}^{p}(V)}:=\|X\|_{L^p(\Omega;C(V))}=\Big(\E\Big[\sup\limits_{t\in [0,T]}\|X(t)\|_V^p\Big]\Big)^{\frac1p},\\
\|X\|_{\mathcal{S}^{p}(V)}:=\|X\|_{C(L^p(\Omega;V))}=\Big(\sup\limits_{t\in [0,T]}\E\Big[\|X(t)\|_V^p\Big]\Big)^{\frac1p}.
\end{align*}
Of course we have $\mathcal Q^p(V)\subseteq\mathcal S^p(V)$ and $\|X\|_{\mathcal S^p(V)}\le \|X\|_{\mathcal Q^p(V)}$ for any $X\in \mathcal Q^p(V)$.

\begin{proposition}\label{prop:exforward}
Let $T>0$, $A:D(A)\to H$ satisfying Assumption \ref{ipotesi}, $B\in L^\infty_s([0,T];\mathcal L(H))$, $j\in L^2_{\mathscr P_0}([0,T]\times\Omega;H)$, $\sigma\in C([0,T];\mathcal L_2(K_0;H))$, $z_0\in L^2_{\mathscr F_0^0}(H)$. Then the problem \eqref{SDEfor} admits a unique solution in $\mathcal Q^2(H)$. Moreover, consider $C_0>0$ such that
\begin{equation}\label{boundcoeff}
\|e^{tA}\|_{L^\infty(\mathcal L(H))}+||B||_{L^\infty(\mathcal L(H))}+\|j\|_{L^2(H)}+\|\sigma\|_{C(\mathcal L_2(K_0;H))}+\|z_0\|_{L^2(H)}\le C_0.
\end{equation}
Then it holds
\begin{equation}\label{boundforward}
\|z\|_{\mathcal Q^2(H)}\le C_1,
\end{equation}
for a certain $C_1$ depending on $C_0$ and $T$.
\end{proposition}

The well-posedness result has been proved in the deterministic case by \cite[Proposition 3.4, page 136]{DPB} and in the stochastic case by \cite[Proposition 1.127 and Proposition 1.147]{FGS}, \cite[Theorem 7.2]{DaPraZa}. We sketch here the proof, since the hypotheses required in \cite{FGS,DaPraZa} are not satisfied here. However, there are just minor modifications to handle in the idea of the proof.

\begin{proof}
Throughout the proof, we use the short notation $\|\cdot\|_{\mathcal L_2}$ instead of $\|\cdot\|_{\mathcal L_2(K_0;H)}$. There is no possibility of mistake since in these results the unique space of Hilbert-Schmidt operators is $\mathcal L_2(K_0;H)$.

Consider $\Phi:\mathcal Q^2(H)\to\mathcal Q^2(H)$ given by
$$
\Phi(z)(t)=e^{tA}z_0+\int_0^t e^{(t-s)A}(B(s)z(s)+j(s))ds+\int_0^t e^{(t-s)A}\sigma_s dW_s^0.
$$
Observe that $\Phi$ is well-posed: since $z$, $\sigma$ and $j$ are progressively measurable, so is $\Phi(z)$. Moreover, since $\sigma\in C([0,T];\mathcal L_2(K_0;H))$ is deterministic, it belongs to $\mathcal Q^p(\mathcal L_2(K_0;H))$ for every $p>2$, and by \cite[Proposition 1.112]{FGS} the stochastic integral admits a continuous modification. Finally, since $j$ is integrable, its convolution with the $C_0$ semigroup has continuous trajectories. This proves that $\Phi(z)$ has a.s. continuous trajectories. To conclude, we have
\begin{align*}
\EQ{\Phi(z)}&\le C\|z_0\|_{L^2(H)}^2+C\Big(||B||_{L^\infty(\mathcal L(H))}^2\|z\|^2_{\mathcal Q^2(H)}+\|j\|_{L^2(H)}^2\Big)\\
&\quad+C\E\left[\sup_{t\in[0,T]}\left\|\int_0^t e^{(t-s)A}\sigma_s\,dW_s^0\right\|_H^2\right].
\end{align*}
The last term can be bounded again using \cite[Proposition 1.112]{FGS}: for any $p>2$,
\[
\E\left[\sup_{t\in[0,T]}\left\|\int_0^t e^{(t-s)A}\sigma_s\,dW_s^0\right\|_H^p\right]\le C_p\E\left[\int_0^T\|\sigma_s\|_{\mathcal L_2}^p ds\right]\le  C\|\sigma\|_{C(\mathcal L_2)}^p,
\]
hence
\begin{equation}\label{BDGp=2}
\E\left[\sup_{t\in[0,T]}\left\|\int_0^t e^{(t-s)A}\sigma_s\,dW_s^0\right\|_H^2\right]\le \E\left[\sup_{t\in[0,T]}\left\|\int_0^t e^{(t-s)A}\sigma_s\,dW_s^0\right\|_H^p\right]^{\frac2p}\le C\|\sigma\|_{C(\mathcal L_2)}^2,
\end{equation}
which proves that $\Phi(z)\in\mathcal Q^2(H)$ if $z\in\mathcal Q^2(H)$.

A solution of \eqref{SDEfor} is a fixed point of $\Phi$. We use Banach-Caccioppoli to prove existence and uniqueness of a fixed point. Observe that we have
\begin{align*}
\EQ{(\Phi(z)-\Phi(w))}&=\E\left[\sup_{t\in[0,T]}\norm{\int_0^t e^{(t-s)A}B(s)(z(s)-w(s))\,ds}^2\right]\\
&\le T^2M_T^2||B||_{L^\infty(\mathcal L(H))}^2\EQ{(z-w)},
\end{align*}
which implies
$$
\|\Phi(z)-\Phi(w)\|_{\mathcal Q^2(H)}^2\le T^2C_0^4\|z-w\|_{\mathcal Q^2(H)}^2,
$$
where $C_0$ is given from \eqref{boundcoeff}. If we choose $T_0=(\sqrt 2C_0^2)^{-1}$, we get that $\Phi$ is a contraction on $[0,T_0]$, and so by Banach-Caccioppoli's Theorem there exists a unique fixed point. This gives the existence and uniqueness of solutions in $[0,T_0]$.

Then, we iterate the procedure. Since we can always choose the same time $T_0=(\sqrt 2C_0^2)^{-1}$, then we obtain another solution in $[T_0,2T_0]$. Sticking together the two solutions, we get a solution in $[0,2T_0]$. Then we continue the iteration obtaining for all $k\in\N$ a (unique) solution in $[0,kT_0]$. When $k>\frac{T}{T_0}$ we get the desired solution in $[0,T]$.

Now we prove \eqref{boundforward}. If $z$ solves \eqref{SDEfor} in the sense of Definition \ref{defmildX}, we get
\begin{align*}
\|z(\tau)\|_H^2&\le 4\,\|e^{\tau A}z_0\|_H^2+4T\int_0^\tau\Big(\|e^{(\tau-s)A}B(s)z(s)\|_H^2\,ds+4T\int_0^\tau\|e^{(\tau-s)A}j(s)\|_H^2\,ds\\
&\quad+4\left\|\int_0^\tau e^{(\tau-s)A}\sigma_s dW_s^0\right\|_H^2\\
&\le 4C_0^4\Bigg[2+T\int_0^\tau\sup\limits_{r\in[0,s]}\|z(r)\|_H^2\,ds+\left\|\int_0^\tau e^{(\tau-s)A}\sigma_sdW_s^0\right\|_H^2\Bigg].
\end{align*}
Passing to the sup for $\tau\in[0,t]$ and computing the average, one has for a certain $C_1$ depending on $T$ and $C_0$
\begin{align*}
\E\Big[\sup\limits_{\tau\in[0,t]}\|z(\tau)\|_H^2\Big]&\le C_1\Bigg(1+\int_0^t \E\Big[\sup\limits_{r\in[0,s]}\|z(r)\|_H^2\Big]ds+\E\Bigg[\sup\limits_{\tau\in[0,t]}\left\|\int_0^\tau e^{(\tau-s)A}\sigma_s dW_s^0\right\|_H^2\Bigg]\Bigg)\\
&\le C_1\Bigg(1+\int_0^t \E\Big[\sup\limits_{r\in[0,s]}\|z(r)\|_H^2\Big]ds+\|\sigma\|_{C(\mathcal L_2)}^2\Bigg)\\
&\le C_1\Bigg(1+\int_0^t \E\Big[\sup\limits_{r\in[0,s]}\|z(r)\|_H^2\Big]ds\Bigg),
\end{align*}
where the constant $C_1$ is changed in the last line and where we have used again \eqref{BDGp=2}. Applying Gronwall's inequality and evaluating at $t=T$, we get \eqref{boundforward} and we conclude the proof.
\end{proof}

\begin{remark}\label{rem:generalizzazione}
We stress the fact that we limit ourselves here to the case $\sigma$ deterministic, since in the MFG system and in the Master Equation we do not have a stochastic diffusion. However, the proposition can be easily modified in the following way:
\begin{itemize}
\item[-] If $\sigma\in \mathcal Q^2(\mathcal L_2(K_0;H))$, then $z\in \mathcal S^2(H)$, with
$$
\sup\limits_{t\in[0,T]}\E\big[\|z(t)\|_H^2\big]\le C_1;
$$
\item[-] If, for some $p>2$, $\sigma\in\mathcal Q^p(\mathcal L_2(K_0;H))$, then $z\in\mathcal Q^2(H)$, and \eqref{boundforward} holds; if moreover $z_0\in L^p_{\mathscr F_0^0}(H)$ and $j\in L^p_{\mathscr P_0}([0,T]\times\Omega;H)$, then $z\in\mathcal Q^p(H)$ and \eqref{boundforward} holds replacing $2$ by $p$.
\end{itemize}
We do not prove this generalization, since the proof can be performed in the same way as the one given before and of \cite[Theorem 7.2]{DaPraZa}.
\end{remark}

We also need a stability result, similar to the one provided in Proposition \ref{prop:stabilita} for the deterministic case.
\begin{proposition}\label{prop:stabilitaforward}
Let $T>0$, and take $A:D(A)\to H$ satisfying Assumption \ref{ipotesi}, and consider the Yosida approximants $\{A_n\}_n\subset \mathcal L(H)$. Consider also $(B_n)_n\in C_s([0,T];\mathcal L(H))$, $j_n\in L^2_{\mathscr P_0}([0,T]\times\Omega;H)$, $\sigma^n\in C([0,T];\mathcal L_2(K_0;H))$ and $z_0^n\in L^2_{\mathscr F_0^0}(H)$, such that
$$
(B_n,j_n,\sigma^n,z_0^n)\to(B,j,\sigma,z_0) \quad\mbox{in }C_s([0,T];\mathcal L(H))\times L^2_{\mathscr P_0}([0,T]\times\Omega;H)\times C(\mathcal L_2(K_0;H))\times L^2_{\mathscr F_0^0}(H).
$$
Calling $z_n$ the unique solution of \eqref{SDEfor} related to the coefficients $(A_n,B_n,j_n,\sigma^n,z_0^n)$, and defining $z$ similarly, we assume that \eqref{boundcoeff} holds for a constant $C_0$ not depending on $n$. Then it holds
$$
z_n\to z\qquad\mbox{in }\mathcal Q^2(H).
$$
\end{proposition}

\begin{proof}
As the proof is simpler than the one provided in Proposition \ref{prop:stabilita}, we just sketch it.

Subtracting the mild formulation of $z_n$ and $z$, we get
\begin{align*}
\|z_n(\tau)-z(\tau)\|_H^2&\le C\|e^{\tau A_n}z_0^n- e^{\tau A}z_0\|_H^2+C\int_0^\tau \|e^{(\tau-s)A_n}B_n(s)z_n(s)-e^{(\tau-s)A}B(s)z(s)\|_H^2\,ds\\
&\quad+C\int_0^\tau\| e^{(\tau-s)A_n}j_n(s)-e^{(\tau-s)A}j(s)\|_H^2\,ds\\
&\quad+C\left\|\int_0^\tau\big(e^{(\tau-s)A_n}\sigma_s^n-e^{(\tau -s)A}\sigma_s\big)dW_s^0\right\|_H^2.
\end{align*}
We pass to the $\sup$ for $\tau\in[0,t]$ and then we compute the average, obtaining
\begin{align*}
\E\Big[\sup\limits_{\tau\in[0,t]}\|z_n(\tau)-z(\tau)\|_H^2\Big]&\le C\Big(\omega_n+ \int_0^t\E\Big[\sup\limits_{r\in[0,s]}\|z_n(r)-z(r)\|_H^2\Big]\,ds\Big)\\
&\quad+C\,\E\left[\sup_{\tau\in[0,T]}\left\|\int_0^\tau\Big(e^{(\tau-s)A_n}\sigma^n_s-e^{(\tau-s)A}\sigma_s\Big)\,dW_s^0\right\|_H^2\right],
\end{align*}
where
\begin{align*}
\omega_n:=\;&
\big\|e^{\cdot A_n}z_0^n-e^{\cdot A}z_0\big\|_{\mathcal Q^2(H)}^2+\E\left[\int_0^T\sup_{\tau\in[s,T]}\big\|e^{(\tau-s)A_n}j_n(s)-e^{(\tau-s)A}j(s)\big\|_H^2\,ds\right]\\
&+\E\left[\int_0^T\sup_{\tau\in[s,T]}\big\|e^{(\tau-s)A_n}B_n(s)z(s)-e^{(\tau-s)A}B(s)z(s)\big\|_H^2\,ds\right].
\end{align*}
As in Proposition \ref{prop:stabilita}, using the strong convergence of the
Yosida semigroups, uniform on compact time intervals, the convergence of the
coefficients, and dominated convergence, we obtain $\omega_n\to0 \qquad\text{as }n\to+\infty$. For the stochastic integral, we have for any $p>2$
\begin{align*}
&\E\left[\sup_{\tau\in[0,T]}\left\|\int_0^\tau\Big(e^{(\tau-s)A_n}\sigma^n_s-e^{(\tau-s)A}\sigma_s\Big)\,dW_s^0\right\|_H^2\right]\\
&\ \le 2\E\left[\sup_{\tau\in[0,T]}\left\|\int_0^\tau
e^{(\tau-s)A_n}\bigl(\sigma^n_s-\sigma_s\bigr)\,dW_s^0\right\|_H^2\right]+
2\E\left[\sup_{\tau\in[0,T]}\left\|\int_0^\tau\bigl(e^{(\tau-s)A_n}-e^{(\tau-s)A}\bigr)\sigma_s\,dW_s^0\right\|_H^2\right]\\
&\ \le C\|\sigma^n-\sigma\|_{C(\mathcal L_2)}^2+C\,\E\left[\sup_{\tau\in[0,T]}\left\|\int_0^\tau\bigl(e^{(\tau-s)A_n}-e^{(\tau-s)A}\bigr)\sigma_s\,dW_s^0\right\|_H^p\right]^{\frac2p},
\end{align*}
where we have used \eqref{BDGp=2}, whose constant is uniform in $n$
thanks to the uniform semigroup bound for the Yosida approximants, to estimate the first integral. Now, the first term goes to $0$ thanks to our hypotheses, and the second term goes to $0$ thanks to \cite[Proposition 1.112]{FGS}. Thus, up to redefining the vanishing term $\omega_n$, we get
$$
\E\Big[\sup\limits_{\tau\in[0,t]}\|z_n(\tau)-z(\tau)\|_H^2\Big]\le C\Big(\omega_n+ \int_0^t\E\Big[\sup\limits_{r\in[0,s]}\|z_n(r)-z(r)\|_H^2\Big]\,ds\Big).
$$

Applying Gronwall's inequality, we get
$$
\|z_n-z\|_{\mathcal Q^2(H)}^2=\EQ{(z_n-z)}\le C\omega_n e^{CT}.
$$
Letting $n\to+\infty$, we obtain $z_n\to z$ in $\mathcal Q^2(H)$, which concludes the proof.
\end{proof}

Now we turn to the study of the backward SDE \eqref{SDEgen}, henceforth BSDE, rewritten here for better readability:
\begin{equation*}
\begin{cases}
dz(t)=-\Big[Az(t)+B(t)z(t)+j(t)\Big]dt-\sigma_t^z\,dW_t^0,\\
z(T)=z_T,
\end{cases}\tag{\ref{SDEgen}}
\end{equation*}

As already said, since \eqref{SDEgen} is backward, a prescribed diffusion coefficient would not yield in general an adapted solution. Hence, \(\sigma^z\) is part of the unknown and must be determined together with \(z\) in order to ensure adaptedness.

\begin{proposition}\label{prop:exbackward}
Let $T>0$, $A:D(A)\to H$ satisfying Assumption \ref{ipotesi}, $B\in L^\infty_s([0,T];\mathcal L(H))$, $j\in L^2_{\mathscr P_0}([0,T]\times\Omega;H)$, $z_T\in L^2_{\mathscr F_T^0}(H)$. Then the problem \eqref{SDEgen} admits a unique solution $(z,\sigma^z)\in L^2_{\mathscr P_0}([0,T]\times\Omega;H\times\mathcal L_2(K_0;H))$. Moreover, consider $C_0>0$ such that
\begin{equation}\label{boundcoeffback}
||e^{tA}||_{L^\infty(\mathcal L(H))}+||B||_{L^\infty(\mathcal L(H))}+\|j\|_{L^2(H)}+\|z_T\|_{L^2(H)}\le C_0.
\end{equation}
Then it holds
\begin{equation}\label{boundbackward}
\|z\|_{\mathcal Q^2(H)}+\|\sigma^z\|_{L^2(\mathcal L_2)}\le C_1,
\end{equation}
for a certain $C_1$ depending on $C_0$.
\end{proposition}

Observe that, differently from the forward case, we do not know here that $z\in\mathcal Q^2(H)$. However, the norm $\|\cdot\|_{\mathcal Q^2(H)}$ is well defined even if $z\not\in\mathcal Q^2(H)$, for example if $z\in L^2(\Omega;L^\infty([0,T];H))$, and an estimate like \eqref{boundbackward} is still possible.

\begin{proof}
Thanks to \cite[Theorem 3.1]{HuPeng1991}, \eqref{SDEgen} admits a unique solution $(z,\sigma^z)\in L^2_{\mathscr P_0}([0,T]\times\Omega;H\times\mathcal L_2(K_0;H))$. We just have to show that \eqref{boundbackward} holds.

From the mild formulation of $(z,\sigma^z)$ we have
$$
z(t)=e^{A(T-t)}z_T +\int_t^T e^{A(s-t)}(B(s)z(s)+j(s))ds+\int_t^T e^{A(s-t)}\sigma^z_sdW_s^0.
$$
Computing the conditional expectation with respect to $\mathscr F_t^0$, we get
\begin{equation}\label{mildcondizionato}
z(t)=e^{A(T-t)}\E[z_T\mid \mathscr F_t^0]+\int_t^T e^{A(s-t)}\Big(B(s)\E[z(s)\mid\mathscr F_t^0]+\E[j(s)\mid\mathscr F_t^0]\Big)ds.
\end{equation}
From here we get, using \eqref{boundcoeffback} and Jensen's inequality,
\begin{align}\label{diventopazzo}
\|z(t)\|_H^2&\le C_0^2\cex{\|z_T\|_H^2}+C_0^2T\int_t^T\Big(C_0^2\cex{\|z(s)\|_H^2}+
\cex{\|j(s)\|_H^2}\Big)\,ds.
\end{align}
Taking the expectation we get
\begin{align*}
\E\big[\|z(t)\|_H^2\big]\le C_0^2\E\big[\|z_T\|_H^2\big]+C_0^2T\int_t^T\left(C_0^2\E\big[\|z(s)\|_H^2\big]+\E\big[\|j(s)\|_H^2\big]\right)\,ds.
\end{align*}
Applying Gronwall's inequality, we obtain for a certain $C_1$ depending on $C_0$

\begin{equation}\label{primastimaz}
\|z\|_{\mathcal S^2(H)}^2=\sup_{t\in[0,T]}\E\big[\|z(t)\|_H^2\big]\le C_1.
\end{equation}
Now, we define
\[
Y:=\|z_T\|_H+\int_0^T\big(C_0\|z(s)\|_H+\|j(s)\|_H\big)\,ds.
\]
From \eqref{mildcondizionato} we deduce
\[
\|z(t)\|_H\le C_0\cex{Y}.
\]
Therefore, using Doob's inequality on the martingale $(\E[Y\mid\mathscr F_t^0])_t$,
\begin{align*}
\|z\|_{\mathcal Q^2(H)}^2=\E\left[\sup_{t\in[0,T]}\|z(t)\|_H^2\right]\le C_0^2\E\left[\sup_{t\in[0,T]}\cex{Y}^2\right]\le 4C_0^2\E[Y^2],
\end{align*}

Finally, using \eqref{primastimaz},
\begin{align*}
\E[Y^2]\le C\left(\E\big[\|z_T\|_H^2\big]+\E\left[\int_0^T\|z(s)\|_H^2\,ds\right]+\E\left[\int_0^T\|j(s)\|_H^2\,ds\right]\right)\le C_1,
\end{align*}
where the constant $C_1$ has been changed in the last inequality. Hence
\begin{equation}\label{boundzQ}
\|z\|_{\mathcal Q^2(H)}\le C_1.
\end{equation}

To conclude, we apply \cite[Lemma 2.1]{HuPeng1991} to obtain
\begin{align*}
\|\sigma^z\|_{L^2(\mathcal L_2)}^2&\le C\left(\E\big[\|z_T\|_H^2\big]+\E\Big[\int_0^T\|B(s)z(s)+j(s)\|_H^2\,ds\Big]\right)\\
&\le C\left(C_0^2+C_0^2\int_0^T\E\big[\|z(s)\|_H^2\big]\,ds\right)\le C_1,
\end{align*}
where in the last inequality we have used \eqref{primastimaz}. This, together with \eqref{boundzQ}, proves \eqref{boundbackward} and concludes the proof.
\end{proof}

\begin{remark}
It is possible to prove that, if $(z,\sigma^z)$ solves \eqref{SDEgen} and \eqref{boundcoeffback} holds, then $z\in\mathcal S^2(H)$. Since we know that \eqref{boundbackward} holds, it is enough to show that
$$
\lim\limits_{r\to t}\E[\|z(t)-z(r)\|_H^2]=0.
$$
This can be done writing $z(t)-z(r)$ from \eqref{mildcondizionato} and using the same ideas as Proposition \ref{prop:exbackward}. We do not provide here the proof, since it is similar to what has been done before and it is not relevant for the well-posedness of the MFG system.
\end{remark}

We conclude this section with a stability result for the BSDE.
\begin{proposition}\label{stabilitaback}
Let $T>0$, take $A:D(A)\to H$ satisfying Assumption \ref{ipotesi}, and consider the Yosida approximants $\{A_n\}_n\subset \mathcal L(H)$. Consider also $(B_n)_n\in C_s([0,T];\mathcal L(H))$, $j_n\in\mathcal Q^2(H)$ and $z_T^n\in L^2_{\mathscr F_T^0}(H)$, such that
$$
(B_n,j_n,z_T^n)\to(B,j,z_T) \quad\mbox{in }C_s([0,T];\mathcal L(H))\times L^2_{\mathscr P_0}([0,T]\times\Omega;H)\times L^2_{\mathscr F_T^0}(H).
$$
Calling $(z_n,\sigma^z_n)$ the unique solution of \eqref{SDEgen} related to the coefficients $(A_n,B_n,j_n,z_T^n)$, and defining $(z,\sigma^z)$ similarly, we assume that \eqref{boundcoeffback} holds for a constant $C_0$ not depending on $n$. Then it holds
$$
(z_n,\sigma_n^z)\to (z,\sigma^z)\qquad\mbox{in }L^2_{\mathscr P_0}([0,T]\times\Omega;H\times\mathcal L_2(K_0;H)).
$$
Moreover, we have
\begin{equation}\label{convznback}
\lim\limits_{n\to+\infty}\norm{z_n-z}_{\mathcal Q^2(H)}=0.
\end{equation}
\end{proposition}

\begin{proof}
Using the formulation \eqref{mildcondizionato} for $z_n$ and $z$, we obtain
\begin{align*}
z_n(t)-z(t)
={}&e^{A_n(T-t)}\E[z_T^n-z_T\mid\mathscr F_t^0]+(e^{(T-t)A_n}-e^{(T-t)A})\E[z_T\mid\mathscr F_t^0]\\
&+\int_t^T e^{(s-t)A_n}B_n(s)\E[z_n(s)-z(s)\mid\mathscr F_t^0]\,ds\\
&+\int_t^T e^{(s-t)A_n}(B_n(s)-B(s))\E[z(s)\mid\mathscr F_t^0]\,ds\\
&+\int_t^T e^{(s-t)A_n}\E[j_n(s)-j(s)\mid\mathscr F_t^0]\,ds\\
&+\int_t^T (e^{(s-t)A_n}-e^{(s-t)A})\E[B(s)z(s)+j(s)\mid\mathscr F_t^0]\,ds.
\end{align*}
Therefore, for a constant $C$ independent of $n$,
\begin{equation}\label{convznS}
\|z_n(t)-z(t)\|_H\le C\E\left[\omega_n+\int_t^T\|z_n(s)-z(s)\|_H\,ds\mid\mathscr F_t^0\right],
\end{equation}
where
\begin{align*}
\omega_n:={}&\|z_T^n-z_T\|_H+\sup_{r\in[0,T]}\|(e^{rA_n}-e^{rA})z_T\|_H+\int_0^T\sup_{r\in[0,T]}\|(e^{rA_n}-e^{rA})(B(s)z(s)+j(s))\|_H\,ds\\
&+\int_0^T\|j_n(s)-j(s)\|_H\,ds+\int_0^T\|(B_n(s)-B(s))z(s)\|_H\,ds.
\end{align*}
Using Jensen's inequality for conditional expectations, we obtain
\[
\E[\|z_n(t)-z(t)\|_H^2]\le C\E[\omega_n^2]+C\int_t^T\E[\|z_n(s)-z(s)\|_H^2]\,ds.
\]
Therefore, by Gronwall's inequality,
\begin{equation*}
\sup_{t\in[0,T]}\E[\|z_n(t)-z(t)\|_H^2]\le C\E[\omega_n^2].
\end{equation*}

Thanks to our hypotheses, we have $\E[\omega_n^2]\to0$. Then
\begin{equation}\label{ultima}
\sup_{t\in[0,T]}\E[\|z_n(t)-z(t)\|_H^2]\to0.
\end{equation}

We now set
\[
\widetilde\omega_n:=\omega_n+\int_0^T\|z_n(s)-z(s)\|_H\,ds.
\]
From \eqref{convznS}, we have
\[
\|z_n(t)-z(t)\|_H\le C\E[\widetilde\omega_n\mid\mathscr F_t^0].
\]
Hence, by Doob's inequality,
\begin{align*}
\|z_n-z\|_{\mathcal Q^2(H)}^2\le C\E\left[\sup_{t\in[0,T]}\E[\widetilde\omega_n\mid\mathscr F_t^0]^2\right]\le C\E[\widetilde\omega_n^2]\to0,
\end{align*}
where the last convergence is due to \eqref{ultima}. This proves \eqref{convznback}, and in particular we have
\[
z_n\to z\qquad\mbox{in }L^2_{\mathscr P_0}([0,T]\times\Omega;H).
\]

To conclude, we have to prove the convergence of $\sigma_n^z$. First of all, we call
\[
f_n(s):=B_n(s)z_n(s)+j_n(s),\qquad f(s):=B(s)z(s)+j(s).
\]
Then, from the previous estimates, $f_n\to f$ in $L^2_{\mathscr P_0}([0,T]\times\Omega;H)$. By the martingale representation theorem, there exist
$\eta_n,\eta\in L^2_{\mathscr P_0}([0,T]\times\Omega;\mathcal L_2(K_0;H))$
and adapted $\kappa_n,\kappa$ such that
\[
z_T^n=\E[z_T^n]+\int_0^T\eta_n(r)\,dW_r^0,\qquad z_T=\E[z_T]+\int_0^T\eta(r)\,dW_r^0,
\]
and, for almost every $s\in[0,T]$,
\[
f_n(s)=\E[f_n(s)]+\int_0^s\kappa_n(s,r)\,dW_r^0,\qquad f(s)=\E[f(s)]+\int_0^s\kappa(s,r)\,dW_r^0.
\]
Subtracting the two formulations we get
\begin{equation}\label{conveta}
\E\left[\int_0^T\|\eta_n(r)-\eta(r)\|_{\mathcal L_2}^2\,dr\right]\le \E[\|z_T^n-z_T\|_H^2]\to0,
\end{equation}
\begin{equation}\label{convkappa}
\E\left[\int_0^T\int_0^s\|\kappa_n(s,r)-\kappa(s,r)\|_{\mathcal L_2}^2\,dr\,ds\right]\to0.
\end{equation}

The construction of the solution in \cite[Lemma 2.1]{HuPeng1991} yields
\begin{gather*}
\sigma_n^z(t)=-e^{(T-t)A_n}\eta_n(t)-\int_t^T e^{(s-t)A_n}\kappa_n(s,t)\,ds,\\
\sigma^z(t)=-e^{(T-t)A}\eta(t)-\int_t^T e^{(s-t)A}\kappa(s,t)\,ds.	
\end{gather*}
Therefore,
\begin{align*}
\sigma_n^z(t)-\sigma^z(t)
={}&-e^{(T-t)A_n}(\eta_n(t)-\eta(t))-(e^{(T-t)A_n}-e^{(T-t)A})\eta(t)\\
&-\int_t^T e^{(s-t)A_n}(\kappa_n(s,t)-\kappa(s,t))\,ds\\
&-\int_t^T (e^{(s-t)A_n}-e^{(s-t)A})\kappa(s,t)\,ds.
\end{align*}
Using \eqref{conveta}, \eqref{convkappa}, the uniform boundedness of the semigroups and the dominated convergence theorem, we obtain
\[
\E\left[\int_0^T\|\sigma_n^z(t)-\sigma^z(t)\|_{\mathcal L_2(K_0;H)}^2\,dt\right]\to0,
\]
which concludes the proof.
\end{proof}

\section{The MFG system: well-posedness}\label{sec:mfg}
We are now ready to prove existence and uniqueness of solutions for the MFG system \eqref{mfg} in the LQ case, in the sense of Definition \ref{def:LQsolMFG}.

The main idea is to prove that the system \eqref{eq:P-Riccati}-\eqref{eq:s} is equivalent to a simpler decoupled system, whose existence and uniqueness of solutions is a natural consequence of the results obtained in the previous section.

We start with the existence result.

\begin{theorem}\label{thm:exMFG}
Let Assumptions \ref{ipotesi}, \ref{ipotesi2}, \ref{ipotesi3} hold true. Then there exists a LQM mild solution
\begin{equation}\label{reg1}
(P,Y,r,Z^r,s,Z^s)\in C_s([0,T];\mathcal L(H))\times\mathcal Q^2(H)\times L^2_{\mathscr P_0}([0,T]\times\Omega;(H\times\mathcal L_2(K_0;H))^2)
\end{equation}
of the MFG system \eqref{mfg}, in the sense of Definition \ref{def:LQsolMFG}.
\end{theorem}

\begin{proof}
We need to prove that there exists a mild solution of the system \eqref{eq:P-Riccati}-\eqref{eq:s}. We consider the $6$-tuple
\begin{equation}\label{reg2}
(P,K,\ell,Y,s,Z^s)\in C_s([0,T];\mathcal L(H)^2\times H)\times \mathcal Q^2(H)\times L^2_{\mathscr P_0}([0,T]\times\Omega;H\times\mathcal L_2(K_0;H)),
\end{equation}
which is a mild solution of \eqref{eq:P-Riccati} and
\begin{align}
&\begin{cases}
K'+K(A-D_1+D_2)+(A^*-D_1^*)K-K\mathcal B K+(Q+S)=0,\\
K(T)=Q_T + S_T,
\end{cases}\label{eq:K-Riccati}\\
&\begin{cases}
\ell'+(A^*-D_1^*-K\mathcal B)\ell-Kc_3+\eta=0,\\
\ell(T)=\eta_T,
\end{cases}\label{eq:elle}\\
&\begin{cases}
dY_t=\Bigl[(A-D_1+D_2-\mathcal BK)Y_t-\mathcal B\ell-c_3\Bigr]\dd t
+\sigma_0\dd W_t^0,\\
Y_0=\bar m_0,
\end{cases}\label{eq:Y2maèsemprelei}\\
&\begin{cases}
ds(t)=-\Big[\mathrm{tr}(\mathcal AP+\mathcal A_0(2K^*-P))+\langle (K-P)Y+\ell,D_2Y-c_3\rangle+\frac12\langle ZY,Y\rangle\Big]dt\\
\qquad \ -\Big[\langle\zeta,Y\rangle-\frac12\langle\mathcal B(K-P)Y+\ell,(K-P)Y+\ell\rangle+\lambda\Big]dt+\langle Z_t^s,dW_t^0\rangle,\\
s(T)=\frac12\langle Z_TY_T,Y_T\rangle+\langle \zeta_T,Y_T\rangle+\lambda_T.
\end{cases}\label{eq:s2maèsemprelei}
\end{align}
We want to prove that, if we define
\begin{equation}\label{relazioni}
Z^r=(K-P)\sigma_0,\qquad r=(K-P)Y+\ell,
\end{equation}
then the $6$-tuple $(P,Y,r,Z^r,s,Z^s)$ satisfies the regularity of \eqref{reg1} and is a LQM mild solution of \eqref{mfg}.

First of all, observe that the equations \eqref{eq:P-Riccati} and \eqref{eq:K-Riccati} are decoupled from the others, and their existence of solution is guaranteed from Proposition \ref{prop:wellposedness_general}. In particular, \eqref{eq:P-Riccati} is a symmetric equation satisfying hypotheses of Lemma \ref{LaStimaAPrioriPerEccellenza}. For the Riccati equation \eqref{eq:K-Riccati}, we distinguish three cases.

If hypothesis $(2-a)$ of Assumption \ref{ipotesi3} is satisfied, then we are in the symmetric case \emph{(iii-a)} of Lemma \ref{LaStimaAPrioriPerEccellenza}, with operators $B_1=\frac12 d_2 I_H-D_1^*$ and $B_2=\frac12 d_2 I_H-D_1^*-D_1$. This allows us to apply Proposition \ref{prop:wellposedness_general} and obtain the existence and uniqueness of $K$.

If, instead, hypothesis $(2-b)$ of Assumption \ref{ipotesi3} is satisfied, then we are in the non-symmetric case \emph{(iii-a)} of the same Lemma, and we apply again Proposition \ref{prop:wellposedness_general}. Finally, if hypotheses $(2-c)$ is satisfied, we are in the case \emph{(i)} of the same Lemma, and we immediately conclude $K=0$.

Moreover, the coefficients appearing in \eqref{eq:K-Riccati}
and their adjoints are strongly continuous. Indeed, by
Assumptions \ref{ipotesi2} and \ref{ipotesi3} the operators $
\mathcal B, D_1, D_2,Q,  S$ are continuous in operator norm, and hence both they and their
adjoints belong to $C_s([0,T];\mathcal L(H))$.
Therefore, by Remark~\ref{rem:adjoint-Cs},
\[
K^*\in C_s([0,T];\mathcal L(H)).
\]
Since $Y_0=\bar m_0$ and $\sigma_0$ are deterministic, Remark \ref{rem:generalizzazione} applied with $p=4$ yields $Y\in Q^4(H)$. Moreover, since $P$, $K$ and $\ell$ are deterministic and bounded, from $r=(K-P)Y+\ell$ we also have $r\in Q^4(H)$. Hence, the driver and the terminal condition in \eqref{eq:s2maèsemprelei} are square-integrable, and Proposition \ref{prop:exbackward} yields the existence and uniqueness of the solution $(s,Z^s)$.
Hence, we can solve first \eqref{eq:elle} thanks to Remark \ref{backlindet} and then \eqref{eq:Y2maèsemprelei} thanks to Proposition \ref{prop:exforward}, \eqref{eq:s2maèsemprelei} thanks to Proposition \ref{prop:exbackward}. The same results imply the desired regularity \eqref{reg2}.

Substituting the definition of $r$ into \eqref{eq:Y2maèsemprelei} and \eqref{eq:s2maèsemprelei}, we get \eqref{eq:Y-r-forward} and \eqref{eq:s}. Moreover from \eqref{relazioni} we get $(r,Z^r)\in L^2_{\mathscr P_0}([0,T]\times\Omega;H\times\mathcal L_2(K_0;H))$. We just have to prove that $(r,Z^r)$ solves \eqref{eq:r-backward}.

To do so, we consider the Yosida approximants $(A_n)_n$ of $A$ and the solutions $P_n,K_n,\ell_n, Y^n$ of
\begin{align}
&\begin{cases}
P'_n+P_n(A_n-D_1)+(A_n^*-D_1^*)P_n-P_n\mathcal BP_n+Q=0,\\
P_n(T)=Q_T,
\end{cases}\label{eq:Pn}\\
&\begin{cases}
K_n'+K_n(A_n-D_1+D_2)+(A_n^*-D_1^*)K_n-K_n\mathcal B K_n+(Q+S)=0,\\
K_n(T)=Q_T + S_T,
\end{cases}\label{eq:Kn}\\
&\begin{cases}
\ell_n'+(A_n^*-D_1^*-K_n\mathcal B)\ell_n-K_nc_3+\eta=0,\\
\ell_n(T)=\eta_T,
\end{cases}\label{eq:ellen}\\
&\begin{cases}
dY^n_t=\Bigl[(A_n-D_1+D_2-\mathcal BK_n)Y^n_t-\mathcal B\ell_n-c_3\Bigr]\dd t
+\sigma_0\dd W_t^0,\\
Y^n_0=\bar m_0,
\end{cases}\label{eq:Yn}
\end{align}
The existence and uniqueness of solutions for \eqref{eq:Pn} and \eqref{eq:Kn} is guaranteed again from Proposition \ref{prop:wellposedness_general}. Moreover, thanks to Proposition \ref{prop:stabilita}, we have $P_n\to P$ and $K_n\to K$ in $C_s([0,T];\mathcal L(H))$.

As regards \eqref{eq:ellen}, the existence and uniqueness of solutions is guaranteed from Remark \ref{backlindet}, and for the same remark we have $\ell_n\to\ell$ in $C([0,T];H)$. Finally, \eqref{eq:Yn} admits a unique solution for Proposition \ref{prop:exforward}, and thanks to Proposition \ref{prop:stabilitaforward} we have $Y^n\to Y$ in $\mathcal Q^2(H)$.

Then we define $r_n=(K_n-P_n)Y^n+\ell_n$. Since $A_n\in\mathcal L(H)$, the solutions of \eqref{eq:Pn}-\eqref{eq:Yn} are classical, and we can compute directly the equation satisfied by $r_n$. We have

\begin{align*}
dr_n={}&(K'_n-P'_n)Y^n\,dt+(K_n-P_n)dY^n+\ell_n' \,dt\\
={}&\Big[-K_n(A_n-D_1+D_2)Y^n-(A_n^*-D_1^*)K_nY^n
+K_n\mathcal BK_nY^n-(Q+S)Y^n\Big]dt\\
&+\Big[P_n(A_n-D_1)Y^n+(A_n^*-D_1^*)P_nY^n
-P_n\mathcal BP_nY^n+QY^n\Big]dt\\
&+\Big[(K_n-P_n)(A_n-D_1+D_2-\mathcal BK_n)Y^n
-(K_n-P_n)\mathcal B\ell_n-(K_n-P_n)c_3\Big]dt\\
&+\Big[K_nc_3-(A_n^*-D_1^*-K_n\mathcal B)\ell_n-\eta\Big]dt
+(K_n-P_n)\sigma_0\,dW_t^0\\
={}&-\Big[(A^*_n-D_1^*)\big[(K_n-P_n)Y^n+\ell_n\big]+SY^n-P_n\mathcal B\big[(K_n-P_n)Y^n+\ell_n\big]\Big]dt\\
&+\Big[P_nc_3-P_nD_2Y^n-\eta\Big]dt+(K_n-P_n)\sigma_0dW_t^0\\
={}&-\Big[(A^*_n-D_1^*-P_n\mathcal B)r_n+(P_nD_2+S)Y^n-P_nc_3+\eta\Big]dt+(K_n-P_n)\sigma_0 dW_t^0.
\end{align*}

For the terminal condition, we have
$$
r_n(T)=(Q_T+S_T-Q_T)Y^n_T+\eta_T=S_TY^n_T+\eta_T.
$$

Since $r_n=(K_n-P_n)Y_n+\ell_n$, we have $r_n\to (K-P)Y+\ell=r$ in $\mathcal Q^2(H).$ To justify the convergence of the products, observe that
\[
(K_n-P_n)Y^n-(K-P)Y
=
(K_n-P_n)(Y^n-Y)
+\big[(K_n-P_n)-(K-P)\big]Y.
\]
The first term converges to zero in $Q^2(H)$ thanks to the uniform
boundedness of $(K_n-P_n)_n$ in
$C_s([0,T];\mathcal L(H))$ and to $Y^n\to Y$ in $Q^2(H)$.
For the second term, for almost every $\omega\in\Omega$, the set
\[
K_\omega:=\{Y_t(\omega):t\in[0,T]\}
\]
is compact in $H$, since $Y(\omega)$ has continuous trajectories.
Hence, by Lemma~\ref{ascoli},
\[
\sup_{t\in[0,T]}
\left\|
\big[(K_n(t)-P_n(t))-(K(t)-P(t))\big]Y_t(\omega)
\right\|_H
\longrightarrow0.
\]
Using the uniform boundedness of $(K_n-P_n)_n$ and dominated
convergence, we conclude that
\[
(K_n-P_n)Y^n\longrightarrow (K-P)Y
\qquad\text{in }Q^2(H).
\] Moreover, we have $(K_n-P_n)\sigma_0\to Z^r$ in $L^2_{\mathscr P_0}([0,T]\times\Omega;\mathcal L_2(K_0;H))$ and $r_n(T)\to S_TY_T+\eta_T$ in $L^2_{\mathscr F_T^0}(H)$.
Indeed, setting
\[
R_n:=(K_n-P_n)-(K-P),
\]
we have $R_n\to0$ in $C_s([0,T];\mathcal L(H))$ and
\[
\sup_n\|R_n\|_{L^\infty(\mathcal L(H))}<+\infty.
\]
For every $t\in[0,T]$ and every orthonormal basis $(e_j)_j$ of $K_0$,
\[
\|R_n(t)\sigma_0(t)\|_{\mathcal L_2(K_0;H)}^2
=
\sum_{j=1}^{\infty}
\|R_n(t)\sigma_0(t)e_j\|_H^2.
\]
By the strong convergence of $R_n(t)$, each term in the sum
converges to zero, while
\[
\|R_n(t)\sigma_0(t)e_j\|_H^2
\le
C\|\sigma_0(t)e_j\|_H^2.
\]
Hence, by dominated convergence in the series,
\[
\|R_n(t)\sigma_0(t)\|_{\mathcal L_2(K_0;H)}
\longrightarrow0
\qquad\text{for every }t\in[0,T].
\]
Since
\[
\|R_n(t)\sigma_0(t)\|_{\mathcal L_2(K_0;H)}
\le
C\|\sigma_0(t)\|_{\mathcal L_2(K_0;H)},
\]
another application of the dominated convergence theorem in time
gives the desired convergence. This implies, using Proposition \ref{stabilitaback}, that we can take the limit in the equation satisfied by $r_n$ and obtain that $(r,Z^r)$ satisfies \eqref{eq:r-backward} in the sense of Definition \ref{deflinback}.
\end{proof}

\begin{remark}\label{rem:digiotto}
Thanks to \eqref{eq:stimalin}, \eqref{eq:stimaquad}, \eqref{eq:estlin2}, \eqref{boundforward} and \eqref{boundbackward} we also know that
\begin{equation}\label{stimamfgmodified}
||P||_{L^\infty(\mathcal L(H))}+||K||_{L^\infty(\mathcal L(H))}+\nch{\ell}+\norm{Y}_{\mathcal Q^2(H)}+\norm{s}_{\mathcal Q^2(H)}+\norm{Z^s}_{L^2(\mathcal L_2)}\le C_0,
\end{equation}
for a constant $C_0$ depending on the norms of the coefficients $e^{tA}$, $D_1$, $D_2$, $\mathcal B$, $Q$, $S$, $\eta$, $c_3$, $\sigma_0$, $Q_T$, $S_T$, $\eta_T$, $\bar m_0$, and on $\delta$ and $\beta$ in Assumption \ref{ipotesi2}, if case \emph{(2b)} is satisfied. This implies, for a possibly different constant $C_0$,
\begin{equation}\label{stimamfg}
||P||_{L^\infty(\mathcal L(H))}+\norm{r}_{\mathcal Q^2(H)}+\norm{Z^r}_{L^2(\mathcal L_2)}+\norm{Y}_{\mathcal Q^2(H)}+\norm{s}_{\mathcal Q^2(H)}+\norm{Z^s}_{L^2(\mathcal L_2)}\le C_0.
\end{equation}

Moreover, during the proof we have shown that $P_n\to P$ in $C_s([0,T];\mathcal L(H))$, $Y_n\to Y$ and $r_n\to r$ in $\mathcal Q^2(H)$, and $Z^{r,n}\to Z^r$ in $L^2_{\mathscr P_0}([0,T]\times\Omega;\mathcal L_2(K_0;H))$. Repeating the same arguments with $p=4$, we obtain $Y_n\to Y$ and $r_n\to r$ in $\mathcal Q^4(H)$; hence, Proposition \ref{stabilitaback} applied to the equation for $s^n$ yields $s^n\to s$ in $\mathcal Q^2(H)$ and $Z^{s,n}\to Z^s$ in $L^2_{\mathscr P_0}([0,T]\times\Omega;\mathcal L_2(K_0;H))$.

In particular, we easily have
\begin{equation}\label{madonnabenedettadellincoronata}
\begin{split}
&(P_n,Y_n,r_n,Z^{r,n},s^n,Z^{s,n})\to (P,Y,r,Z^{r},s,Z^{s})\qquad\mbox{in }C_s(\mathcal L(H))\times\mathcal Q^2(H)^2\times L^2_{\mathscr P_0}(\mathcal L_2\times H\times\mathcal L_2),\\
&u_n(\cdot,x)\to u(\cdot,x)\qquad\mbox{in }\mathcal Q^1(\R)\quad\forall\, x\in H,\\
&D_xu_n(\cdot,x)\to D_x u(\cdot,x)\qquad\mbox{in }\mathcal Q^2(H)\quad\forall\, x\in H.
\end{split}	
\end{equation}
\end{remark}

We conclude with the uniqueness result.

\begin{theorem}\label{thm:uniqMFG}
Let Assumptions \ref{ipotesi}, \ref{ipotesi2}, \ref{ipotesi3} hold. Then the LQM mild solution of the MFG system \eqref{mfg} is unique in the class given by \eqref{reg1}.
\end{theorem}

\begin{proof}
Consider the solution $(P,Y,r,Z^r,s,Z^s)$ obtained in the previous theorem, and related to the solution $(P,K,\ell,Y,s,Z^s)$ of \eqref{eq:P-Riccati} and \eqref{eq:K-Riccati}-\eqref{eq:s2maèsemprelei}.

We take another solution $(\tilde P,\tilde Y,\tilde r,\tilde Z^r,\tilde s,\tilde Z^s)$ of \eqref{eq:P-Riccati}-\eqref{eq:s}, satisfying \eqref{reg1}, and we want to prove that
$$
(\tilde P,\tilde Y,\tilde r,\tilde Z^r,\tilde s,\tilde Z^s)=(P,Y,r,Z^r,s,Z^s).
$$
Since \eqref{eq:P-Riccati} is decoupled from the others, we immediately have for Proposition \ref{prop:wellposedness_general} $\tilde P=P$. Then, we define
$$
\tilde\ell\in L^2_{\mathscr P_0}([0,T]\times\Omega;H),\qquad \tilde\ell:=(P-K)\tilde Y+\tilde r.
$$
To compute the equation satisfied by $\tilde\ell$, we argue as always with an approximating argument. Hence, we consider the Yosida approximants $(A_n)_n$, the functions $P_n,K_n$ as in \eqref{eq:Pn},\eqref{eq:Kn}, and we define $\tilde Y^n$, $\tilde r_n, \tilde Z^{r,n}$ as the solutions of
\begin{align}
&\begin{cases}\label{eq:rntil}
d\tilde r_n(t)=-\Bigl[(A_n^*-D_1^*-P_n\BB) \tilde r_n(t)+( P_nD_2+S)\tilde Y_t-P_nc_3+\eta\Bigr]\dd t
+\tilde Z^{r,n}_t\dd W_t^0,\\
\tilde r_n(T)= S_T\tilde Y_T+\eta_T,
\end{cases}\\
&\begin{cases}\label{eq:Yntil}
d\tilde Y^n_t=\Bigl[(A_n-D_1+D_2)\tilde Y^n_t-\BB P_n\tilde Y^n_t-\BB  \tilde r_n(t)-c_3\Bigr]\dd t
+\sigma_0\dd W_t^0,\\
\tilde Y^n_0=\bar m_0.
\end{cases}
\end{align}
Observe that the equation for $\tilde r_n$ does not depend on $\tilde Y_n$. Hence, the solutions $(\tilde r_n,\tilde Z^{r,n})$ and $\tilde Y^n$ are uniquely well defined, thanks to Proposition \ref{prop:exbackward} and \ref{prop:exforward}.

Since $A_n\in\mathcal L(H)$ and $P_n,K_n$ are continuously differentiable,
the processes $\tilde Y^n$ and $\tilde r_n$ are strong Itô solutions. Then we define $\tilde\ell_n:=(P_n-K_n)\tilde Y_n+\tilde r_n$ and we can use the Itô product formula to compute the equation satisfied by it. We have

\begin{align*}
d\tilde\ell_n={}&(P'_n-K'_n)\tilde Y^n\,dt+(P_n-K_n)\,d\tilde Y^n+d\tilde r_n\\
={}&\Big[-P_n(A_n-D_1)-(A_n^*-D_1^*)P_n+P_n\mathcal BP_n-Q\Big]\tilde Y^n\,dt\\
&+\Big[K_n(A_n-D_1+D_2)+(A^*_n-D_1^*)K_n-K_n\mathcal B K_n+Q+S\Big]\tilde Y^n\,dt\\
&+(P_n-K_n)\Big[(A_n-D_1+D_2)\tilde Y^n-\mathcal BP_n\tilde Y^n-\mathcal B\tilde r_n-c_3\Big]dt\\
&-\Big[(A_n^*-D_1^*-P_n\mathcal B)\tilde r_n
+(P_nD_2+S)\tilde Y-P_nc_3+\eta\Big]dt+\Big[(P_n-K_n)\sigma_0+\tilde Z^{r,n}\Big]dW_t^0\\
={}&\Big[-(A_n^*-D_1^*)(P_n-K_n)\tilde Y^n
+K_n\mathcal B(P_n-K_n)\tilde Y^n+(P_nD_2+S)\tilde Y^n\Big]dt\\
&+\Big[-(A_n^*-D_1^*)\tilde r_n+K_n\mathcal B\tilde r_n
+K_nc_3-\eta-(P_nD_2+S)\tilde Y\Big]dt+\Big[(P_n-K_n)\sigma_0+\tilde Z^{r,n}\Big]dW_t^0.
\end{align*}

Recalling the definition of $\tilde\ell_n$, we get the following BSDE:
\begin{equation}\label{eq:elln}
\begin{cases}
d\tilde\ell_n=-\Big[(A_n^*-D_1^*-K_n\mathcal B)\tilde\ell_n-K_nc_3+\eta\Big]dt+(P_nD_2+S)(\tilde Y^n-\tilde Y)\,dt\\[3mm] \hspace{1.09cm}+\Big[(P_n-K_n)\sigma_0+\tilde Z^{r,n}\Big]dW_t^0,\\[3mm]
\tilde\ell_n(T)=S_T(\tilde Y_T-\tilde Y^n_T)+\eta_T.
\end{cases}
\end{equation}

Thanks to Proposition \ref{prop:stabilita}, we have $P_n\to P$, $K_n\to K$ in $C_s([0,T];\mathcal L(H))$. Then, Proposition \ref{stabilitaback} applied to \eqref{eq:rntil} implies that $(\tilde r^n,\tilde Z^{r,n})\to(\tilde r,\tilde Z^r)$ in $L^2_{\mathscr P_0}([0,T]\times\Omega;H\times\mathcal L_2(K_0;H))$; hence Proposition \ref{prop:stabilitaforward} applied to \eqref{eq:Yntil} gives $\tilde Y^n\to\tilde Y$ in $\mathcal Q^2(H)$. This implies $\tilde\ell_n\to\tilde\ell$ in $L^2_{\mathscr P_0}([0,T]\times\Omega;H)$.

Moreover, observe that $P_n-K_n\to P-K$ in $C_s([0,T];\mathcal L(H))$ implies that $(P_n-K_n)\sigma_0\to (P-K)\sigma_0$ in $\mathcal Q^2(\mathcal L_2(K_0;H))$. Hence, applying again Proposition \ref{stabilitaback} to \eqref{eq:elln} we get the equation satisfied by $\tilde\ell$:
\begin{equation}\label{eq:elltil}
\begin{cases}
d\tilde\ell=-\Big[(A^*-D_1^*-K\mathcal B)\tilde\ell-Kc_3+\eta\Big]dt+\Big[(P-K)\sigma_0+\tilde Z^{r}\Big]dW_t^0,\\
\tilde\ell(T)=\eta_T.
\end{cases}
\end{equation}
By Proposition \ref{prop:exbackward}, equation \eqref{eq:elltil}
admits a unique adapted mild solution pair in $L^2_{\mathscr P_0}([0,T]\times\Omega;H\times \mathcal L_2(K_0;H))$. Since $(\ell,0)$ is also a mild solution of \eqref{eq:elltil}, uniqueness gives
$$
\tilde\ell=\ell,
\qquad
\tilde Z^r=(K-P)\sigma_0=Z^r,
$$
where the last equality comes from \eqref{relazioni}. From $\tilde\ell=\ell$ we get
\begin{equation}\label{trucchettodabaraccone}
P\tilde Y+\tilde r=PY+r+K(\tilde Y-Y).
\end{equation}
Substituting it into the equation of $\tilde Y$, we get
$$
\begin{cases}
d\tilde Y=\Big[(A-D_1+D_2)\tilde Y-\mathcal B(PY+r)-\mathcal BK(\tilde Y-Y)-c_3\Big]dt+\sigma_0\, dW_t^0,\\
\tilde Y_0=\bar m_0.
\end{cases}
$$
As before, the solution is unique for Proposition \ref{prop:exforward}. Since $Y$ solves the previous equation, we must have $\tilde Y=Y$. Now from \eqref{trucchettodabaraccone} we immediately have $\tilde r=r$. Finally, the backward SDEs solved by $(\tilde s,\tilde Z^s)$ and $(s,Z^s)$ are now the same, which implies $\tilde s=s$, $\tilde Z^s=Z^s$ and concludes the proof.
\end{proof}

\section{The Master Equation: well-posedness}\label{sec:ME}
The well-posedness of the Master Equation is a  consequence of all the results obtained in Section \ref{sec:preliminari}. We report here the straightforward proof to obtain existence and uniqueness of solutions for \eqref{Riccati:P}-\eqref{Equation:mu}.

\begin{theorem}\label{thm:wpMaster}
Let Assumptions \ref{ipotesi}, \ref{ipotesi2}, \ref{ipotesi3} hold true. Then there exists a unique LQM mild solution $U$ of the Master Equation \eqref{Master}, in the sense of Definition \ref{def:LQMMaster}. Moreover we have
\begin{equation}\label{stimaME}
||P||_{L^\infty(\mathcal L(H))}+||\Upsilon||_{L^\infty(\mathcal L(H))}+||\Gamma||_{L^\infty(\mathcal L(H))}+\nch{\psi}+\nch{\phi}+\norminf{\mu}\le C_0,
\end{equation}
where $C_0$ is a constant only depending on $M_T$ and all the coefficients of the systems \eqref{Riccati:P}-\eqref{Equation:mu}. In the non-symmetric case, $C_0$ may also depend on $\delta,\beta$ and $M_T^-$.
\end{theorem}

\begin{proof}
The Riccati equation for $P$ \eqref{Riccati:P} is the same as \eqref{eq:P-Riccati}, so the existence and uniqueness of solutions has already been proved in Theorem~\ref{thm:exMFG}.  

The uniqueness of solutions for the ODE \eqref{Riccati:Upsilon} is done in Lemma~\ref{LaStimaAPrioriPerEccellenza}, case $(i)$. For the existence, the ODE does not fall into the hypotheses of Proposition \ref{prop:wellposedness_general}. However, we claim that the solution is given by
\begin{equation}\label{formaUps}
\Upsilon(t)=K(t)-P(t),
\end{equation}
where $K$ is defined in \eqref{eq:K-Riccati}. Actually, we have $\Upsilon(T)=K(T)-P(T)=S_T$ and
\begin{align*}
\Upsilon'=K'-P'={}&-K(A-D_1+D_2)-(A^*-D_1^*)K+K\mathcal B K-Q-S\\
&+P(A-D_1)+(A^*-D_1)^*P-P\mathcal BP+Q\\
={}&-(K-P)(A-D_1-\mathcal BP)-(A^*-D_1^*-P\mathcal B)(K-P)\\
&+(K-P)\mathcal B(K-P)-KD_2-S\\
={}&-\Upsilon (A-D_1-\mathcal BP)-(A^*-D_1^*-P\mathcal B)\Upsilon+\Upsilon\mathcal B\Upsilon-(P+\Upsilon)D_2-S,
\end{align*}
which implies \eqref{Riccati:Upsilon}.

Now we consider the ODE \eqref{Equation:Gamma} for $\Gamma$.
It can be written as
\[
\Gamma'
+\Gamma A+A^*\Gamma
+B_1\Gamma+\Gamma B_2+J=0,
\]
where
\[
B_1
=
-D_1^*-(P+\Upsilon^*)\mathcal B+D_2^*, \quad
B_2
=
-D_1-\mathcal B(P+\Upsilon)+D_2,
\]
\[
J
=
-\Upsilon^*\mathcal B\Upsilon
+\Upsilon^*D_2
+D_2^*\Upsilon
+Z.
\]
By the $L^\infty$ bounds of $P$, $\Upsilon$ (and consequently $\Upsilon^*$), of the original coefficients and their adjoints, we obtain
\[
B_1,B_2,J\in L^\infty_s([0,T];\mathcal L(H)).
\]
Therefore Proposition~\ref{prop:wellposedness_general}, case~(i),
applies and yields a unique mild solution
\[
\Gamma\in C_s([0,T];\mathcal L(H)),
\]
with estimate \eqref{eq:stimalin} holding. This, together with the bounds \eqref{stimamfgmodified} and \eqref{formaUps}, proves \eqref{stimaME} for the first three terms $(P,\Upsilon,\Gamma)$.

We also observe that $\Gamma$ is self-adjoint. Indeed, with the
notation introduced above, we have $B_2=B_1^*.$
In fact, since $P=P^*$ and $\mathcal B=\mathcal B^*$,
\[
B_2^*
=
\bigl(-D_1-\mathcal B(P+\Upsilon)+D_2\bigr)^*
=
-D_1^*-(P+\Upsilon^*)\mathcal B+D_2^*
=
B_1.
\]
Moreover,
\[
J^*
=
\bigl(
-\Upsilon^*B\Upsilon
+\Upsilon^*D_2
+D_2^*\Upsilon
+Z
\bigr)^*
=
J,
\]
since $B=B^*$ and $Z=Z^*$. Finally, the terminal datum satisfies $Z_T^*=Z_T.$
Since $B_1,B_2,J\in L^\infty_s([0,T];\mathcal L(H))$ and $B_1^*=B_2, B_2^*=B_1, J^*=J,$
Remark~\ref{rem:adjoint-Cs} yields
\[
\Gamma^*\in C_s([0,T];\mathcal L(H)).
\]
Taking adjoints in the mild formulation of the equation for $\Gamma$,
we see that $\Gamma^*$ solves
\[
(\Gamma^*)'
+\Gamma^*A+A^*\Gamma^*
+B_1\Gamma^*
+\Gamma^*B_2
+J=0,
\qquad
\Gamma^*(T)=Z_T.
\]
Thus $\Gamma^*$ solves the same equation as $\Gamma$. By the
uniqueness result of Proposition~\ref{prop:wellposedness_general},
case~(i), we conclude that $\Gamma^*=\Gamma.$

For the other terms the argument is analogous. Once obtained the unique $(P,\Upsilon,\Gamma)$, existence and uniqueness for \eqref{Equation:psi} is given by Remark~\ref{backlindet}, and in the same way we work with \eqref{Equation:phi}. Moreover, \eqref{eq:estlin2} gives the desired estimate \eqref{stimaME} for $(\psi,\phi)$. Finally, \eqref{Equation:mu} is uniquely solved with the formula \eqref{explicitmu}, which also allows to complete the estimate \eqref{stimaME} and to conclude the proof.
\end{proof}

\section{Verification}\label{sec:8}

Since the solutions constructed for the MFG system \eqref{mfg} and the Master Equation \eqref{Master} are obtained only in a mild sense, we need to verify that they actually represent the objects they are meant to describe.

More precisely, we prove that the mild solution of the Master Equation generates the value function of the MFG system along the corresponding flow of measures and that the mild solution of the MFG system coincides with the value function of the representative agent.

We start with the first statement.

\begin{theorem}\label{thm:v1}
Let Assumptions \ref{ipotesi}, \ref{ipotesi2}, \ref{ipotesi3} hold true. Then, if $(u,m)$ and $U$ are respectively the LQM mild solutions of the MFG system \eqref{mfg} and the Master Equation \eqref{Master}, it holds
\begin{equation}\label{traiettorie}
u(t,x)=U(t,x,m(t))\qquad\forall\,(t,x)\in[0,T]\times H.
\end{equation}
In particular, the following identities between the coefficients hold:
\begin{equation}\label{identità}
\begin{split}
&K=P+\Upsilon,\qquad r=\Upsilon Y+\psi,\hspace{3.15cm} Z^r=\Upsilon\sigma_0,\\
& \ell=\psi,\hspace{1.75cm}s=\frac12\langle\Gamma Y,Y\rangle+\langle\phi,Y\rangle+\mu,\qquad Z^s=\sigma_0^*\big(\Gamma Y+\phi\big).
\end{split}
\end{equation}
\end{theorem}

Observe that, due to the flow property of the MFG system, \eqref{traiettorie} is equivalent to \eqref{eq:ME_mfg} used before to formally derive the formulation of the Master Equation.

\begin{proof}
We start observing that the identities in \eqref{identità} imply \eqref{traiettorie}. If \eqref{identità} holds we have
\begin{align*}
U(t,x,m(t))=\tilde U(t,x,Y_t)={}&\frac12\langle P(t)x,x\rangle+\langle \Upsilon(t)Y_t,x\rangle+\frac12\langle\Gamma(t)Y_t,Y_t\rangle\\
&+\langle\psi(t),x\rangle+\langle\phi(t),Y_t\rangle+\mu(t)\\
={}&\frac12\langle P(t)x,x\rangle+\langle r(t),x\rangle+s(t)=u(t,x).
\end{align*}
Hence, we just need to prove \eqref{identità}. Observe that $K=P+\Upsilon$ has already been proved in Theorem \ref{thm:wpMaster}, equation \eqref{formaUps}. From here we observe that the ODEs \eqref{eq:elle} and \eqref{Equation:psi} are the same, and for uniqueness we obtain $\ell=\psi$.

Using this information in \eqref{relazioni}, we get $Z^r=\Upsilon\sigma_0$ and $r=\Upsilon Y+\psi$. To conclude, we just have to prove that, if we define the couple
$$
(\tilde s,\tilde Z^s)=\Big(\frac12\langle\Gamma(t)Y_t,Y_t\rangle+\langle\phi(t),Y_t\rangle+\mu(t), \sigma_0^*(t)\big(\Gamma(t)Y_t+\phi(t)\big)\Big),
$$
then $(\tilde s, \tilde Z^s)\in L^2_{\mathscr P_0}([0,T]\times\Omega;\R\times\mathcal L_2(K_0;\R))$ and solves the backward SDE \eqref{eq:s} in the sense of Definition \ref{deflinback}.

Observe that, since $Y_0=\bar m_0$ and $\sigma_0$ are both deterministic, they belong respectively to $L^4_{\mathscr F_0^0}(H)$ and $\mathcal Q^4(\mathcal L_2(K_0;H))$. Hence, thanks to Remark~\ref{rem:generalizzazione}, we have $Y\in\mathcal Q^4(H)$. Since $\Gamma$ is deterministic, this ensures
$$
\frac12\langle\Gamma Y,Y\rangle\in\mathcal Q^2(\R).
$$
The other terms are obviously square integrable, and they are all adapted to $\mathscr P_0$. This proves that $(\tilde s,\tilde Z^s)\in L^2_{\mathscr P_0}$.

To prove that $(\tilde s, \tilde Z^s)$ solves \eqref{eq:s}, observe first of all that $\tilde s(T)=s(T)$. Then, we argue as usually with an approximating argument: hence, we consider the Yosida approximants $(A_n)_n$ of $A$ and the solutions $Y_n$, $\Gamma_n$,$\phi_n$ and $\mu_n$ of
\begin{align*}
&\begin{cases}
dY_t^n=\Big[(A_n-D_1+D_2)Y_t^n-\mathcal B(PY_t^n+r(t))-c_3\Big]dt+\sigma_0dW_t^0,\\
Y_0^n=\bar m_0\,;
\end{cases}\\
&\begin{cases}
\Gamma_n'+\Gamma_n\big(A_n- D_1-\mathcal B (P+\Upsilon)\big)&\hspace{-0.3cm}+\,\big(A^*_n- D_1^*-(P+\Upsilon^*)\mathcal B \big)\Gamma_n-\Upsilon^*\mathcal B \Upsilon\\&\hspace{-0.3cm}+\, {(\Upsilon^*+\Gamma_n)D_{2}+D_{2}^*(\Gamma_n+\Upsilon)}+Z =0\,,\\
\Gamma_n(T)=Z_T\,;
\end{cases}\\
&\begin{cases}
\phi_n'+\left(A_n^*-D_{1}^*-\left( P +\Upsilon^* \right) \mathcal B\right) \phi_n &\hspace{-0.3cm}-\,( \Upsilon^* + \Gamma_n )\mathcal B \psi  + {D_{2}^*(\psi+\phi_n)}\\&\hspace{-0.3cm}-\,(\Upsilon^*+\Gamma_n)c_{3}+\zeta=0\,,\\
\phi_n (T)=\zeta_T\,;
\end{cases}\\
&\begin{cases}
\mu_n'+\mbox{tr}((\mathcal A+\mathcal A_0)P)&\hspace{-0.3cm}+\,\mathrm{tr}\big(\mathcal A_0(2\Upsilon^*+\Gamma_n)\big)-\frac12\langle \mathcal B \psi  ,  \psi +2 \phi_n \rangle\\&\hspace{-0.3cm}-\,\langle c_{3},\psi+\phi_n\rangle+\lambda=0\,,\\
\mu _n(T)=\lambda_T\,.
\end{cases}
\end{align*}
Thanks to Proposition \ref{prop:stabilitaforward} and Remark \ref{rem:generalizzazione}, we have $Y^n$ uniformly bounded in $\mathcal Q^p(H)$ for all $p\ge1$ and $Y^n\to Y$ in $\mathcal Q^2(H)$. This ensures $Y^n\to Y$ in $\mathcal Q^4(H)$. Using also Proposition \ref{prop:stabilita} and Remark \ref{backlindet}, we get $\Gamma_n\to\Gamma$ in $C_s([0,T];\Sigma(H))$ (as we did in the previous result, observe that $\Gamma_n$ is symmetric), $\phi_n\to\phi$ in $C([0,T];H)$, $\mu_n\to\mu$ in $C([0,T];\R)$.
Then we define
\begin{align*}
&\tilde s_n^1=\frac12\langle\Gamma_n Y^n,Y^n\rangle,\qquad \tilde s_n^2=\langle\phi_n,Y^n\rangle,\\
&\tilde s_n=\tilde s_n^1+\tilde s_n^2+\mu_n=\frac12\langle\Gamma_n Y^n,Y^n\rangle+\langle\phi_n,Y^n\rangle+\mu_n\,.
\end{align*}
By the previous convergences, arguing as in the proof of Theorem \ref{thm:exMFG}, we get $\tilde s_n\to\tilde s$ in $\mathcal Q^2(\R)$.

To compute the equation satisfied by $\tilde s_n$, we start by $\tilde s_n^1$. We have
\begin{align*}
d\tilde s_n^1={}&\frac12\langle\Gamma'_nY^n,Y^n\rangle dt+\langle\Gamma_nY^n,dY^n\rangle+\frac12\mathrm{tr}(\Gamma_n\sigma_0\sigma_0^*)\\
={}&\frac12\big\langle\Big[-\Gamma_n\big(A_n-D_1-\mathcal B(P+\Upsilon)\big)-\big(A^*_n-D_1^*-(P+\Upsilon^*)\mathcal B\big)\Gamma_n+\Upsilon^*\mathcal B\Upsilon\Big]Y^n,Y^n\big\rangle dt\\
&+\langle\Gamma_n Y^n,\sigma_0 dW_t^0\rangle+\frac12\big\langle\Big[-(\Upsilon^*+\Gamma_n)D_2-D_2^*(\Gamma_n+\Upsilon)-Z\Big]Y^n,Y^n\big\rangle dt\\
&+\big\langle\Gamma_n Y^n,(A_n-D_1+D_2)Y^n-\mathcal B(PY^n+r)-c_3\big\rangle dt+\mathrm{tr}(\mathcal A_0\Gamma_n)dt\\
={}&\Big[\langle\Gamma_n Y^n,\mathcal B(\Upsilon Y^n-r)-c_3\rangle-\langle\Upsilon Y^n,D_2 Y^n\rangle-\frac12\langle ZY^n,Y^n\rangle\Big]dt\\
&+\Big[\frac12\langle\mathcal B\Upsilon Y^n,\Upsilon Y^n\rangle+\mathrm{tr}(\mathcal A_0\Gamma_n)\Big]dt+\langle \sigma_0^*\Gamma_n Y^n,dW_t^0\rangle.
\end{align*}
As regards $\tilde s_n^2$, we have
\begin{align*}
d\tilde s_n^2={}&\langle\phi'_n,Y^n\rangle dt+\langle \phi_n,dY^n\rangle\\
={}&\big\langle -(A^*_n-D_1^*-(P+\Upsilon^*)\mathcal B)\phi_n+(\Upsilon^*+\Gamma_n)(\mathcal B\psi+c_3)-D_2^*(\psi+\phi_n)-\zeta,Y^n\big\rangle dt\\
&+\big\langle\phi_n,(A_n-D_1+D_2)Y^n-\mathcal B(PY^n+r)-c_3\big\rangle dt+\langle\phi_n,\sigma_0 dW_t^0\rangle\\
={}&\big(\big\langle\phi_n,\mathcal B(\Upsilon Y^n-r)-c_3\big\rangle+\big\langle(\Upsilon^*+\Gamma_n)(\mathcal B\psi+c_3)-D_2^*\psi-\zeta, Y^n\big\rangle)dt+\langle\sigma_0^*\phi_n, dW_t^0\rangle.
\end{align*}
Hence, the equation satisfied by $\tilde s_n$ takes the following form:
\begin{align*}
d\tilde s_n={}&d\tilde s_n^1+d\tilde s_n^2+\mu'_n dt\\
={}&\Big[\langle\Gamma_n Y^n,\mathcal B(\Upsilon Y^n-r)-c_3\rangle-\langle\Upsilon Y^n,D_2 Y^n\rangle-\frac12\langle ZY^n,Y^n\rangle\Big]dt\\
&+\Big[\frac12\langle\mathcal B\Upsilon Y^n,\Upsilon Y^n\rangle+\mathrm{tr}(\mathcal A_0\Gamma_n)\Big]dt+\langle \sigma_0^*\Gamma_n Y^n,dW_t^0\rangle\\
&+\big(\big\langle\phi_n,\mathcal B(\Upsilon Y^n-r)-c_3\big\rangle+\big\langle(\Upsilon^*+\Gamma_n)(\mathcal B\psi+c_3)-D_2^*\psi-\zeta, Y^n\big\rangle)dt+\langle\sigma_0^*\phi_n, dW_t^0\rangle\\
&+\Big[-\mbox{tr}((\mathcal A+\mathcal A_0)P)-\,\mathrm{tr}\big(\mathcal A_0(2\Upsilon^*+\Gamma_n)\big)+\frac12\langle \mathcal B \psi  ,  \psi +2 \phi_n \rangle+\,\langle c_{3},\psi+\phi_n\rangle-\lambda\Big]dt\\
={}&\Big[\langle\Gamma_n Y^n+\phi_n,\mathcal B(\Upsilon Y^n-r+\psi)\rangle-\langle \Upsilon Y^n+\psi,D_2Y^n-c_3\rangle-\frac12\langle Z Y^n,Y^n\rangle\Big]dt\\
&+\Big[\frac12\langle \mathcal B(\Upsilon Y^n+\psi),\Upsilon Y^n+\psi\rangle-\langle\zeta, Y^n\rangle-\mathrm{tr}\big((\mathcal A+\mathcal A_0)P\big)-\mathrm{tr}(2\mathcal A_0\Upsilon^*)-\lambda\Big]dt\\
&+\langle\sigma_0^*(\Gamma_n Y^n+\phi_n),dW_t^0\rangle.
\end{align*}
Observe that the last trace term can be rewritten in the following way:
$$
\mathrm{tr}_H(2\mathcal A_0\Upsilon^*)=\mathrm{tr}_H(\sigma_0\sigma_0^*\Upsilon^*)=\mathrm{tr}_H(\sigma_0(Z_t^r)^*)=\mathrm{tr}_{K_0}((Z_t^r)^*\sigma_0),
$$
where we have used the identity of $Z^r$ in \eqref{identità}, proved before.

Now, we use Proposition \ref{stabilitaback} to pass to the limit for $n\to+\infty$. Observe that we have
$$
\Upsilon Y^n+\psi\to \Upsilon Y+\psi=r \qquad\mbox{in }\mathcal Q^4(H).
$$
Moreover, since
$$
\tilde s_n(T)=\frac12\langle Z_TY_T^n,Y_T^n\rangle
+\langle \zeta_T,Y_T^n\rangle+\lambda_T,
$$
the convergence $Y^n\to Y$ in $\mathcal Q^4(H)$ yields
$$
\tilde s_n(T)\longrightarrow\frac12\langle Z_TY_T,Y_T\rangle+\langle\zeta_T,Y_T\rangle+\lambda_T=s(T)
\qquad\mbox{in }L^2(\Omega).
$$
Since we also know that $\sigma_0^*(\Gamma_n Y^n+\phi_n)\to \sigma_0^*(\Gamma Y+\phi)$ in $\mathcal Q^2(\mathcal L_2(K_0;\R))$, then the pair $\big(\tilde s,\sigma_0^*(\Gamma Y+\phi)\big)$, is a mild solution of
\begin{align*}
d\tilde s={}&-\Big[\langle r,D_2 Y-c_3\rangle+\frac12\langle ZY,Y\rangle-\frac12\langle\mathcal Br,r\rangle+\langle\zeta,Y\rangle\Big]dt
\\&-\Big[\mathrm{tr}\big((\mathcal A+\mathcal A_0)P\big)+\mathrm{tr}\big((Z_t^r)^*\sigma_0\big)+\lambda\Big]dt+\langle\sigma_0^*(\Gamma Y+\phi),dW_t^0\rangle,
\end{align*}
which is the same backward SDE as \eqref{eq:s}. For uniqueness (Proposition \ref{prop:exbackward}), we must have
$$
s=\tilde s=\frac12\langle\Gamma Y,Y\rangle+\langle \phi,Y\rangle+\mu,\qquad Z^s=\tilde Z^s=\sigma_0^*(\Gamma Y+\phi),
$$
which ends the proof of the identities in \eqref{identità} and concludes the proof of the Theorem.
\end{proof}

The last result is a verification theorem for the MFG system: even if the solution $(u,m)$ of \eqref{mfg} is just mild, the value function $u$ still represents the value function of the representative agent.

\begin{theorem}\label{thm:v2}
Suppose Assumptions \ref{ipotesi}, \ref{ipotesi2}, \ref{ipotesi3} hold true. Then, if $(u,m)$ solves \eqref{mfg} in the sense of Definition \ref{def:LQsolMFG}, $u$ satisfies \eqref{def:uvaluefunction}, with $J$ defined in \eqref{eq:J}. Moreover, the infimum in \eqref{def:uvaluefunction} is achieved, i.e. there exists an optimal trajectory $X^{opt}$ and an optimal control $\alpha^{opt}$, given in feedback form by
\begin{equation}\label{eq:opt}
\alpha^{opt}(\tau)= -R_0^{-1}\Big(B^*PX^{opt}_\tau+B^*r+L_0^*X^{opt}_\tau+M_0^*Y_\tau+\theta_0\Big),
\end{equation}
such that
$$
u(t,x)=J(t,x,\alpha^{opt}).
$$
Moreover, the optimal trajectory, for a player starting in position $x$ at time $t$, is given by
\begin{equation}\label{ottima}
\begin{cases}
dX^{opt}_\tau=\big[AX_\tau^{opt}-(D_1+\mathcal BP)X_\tau^{opt}+D_2Y_\tau-\mathcal Br-c_3]\,d\tau+\sigma(\tau) dW_\tau+\sigma_0(\tau) dW_\tau^0\,,\\
X^{opt}_t=x.
\end{cases}
\end{equation}
\end{theorem}

\begin{proof}
We fix $t\in[0,T)$, $x\in H$ and $\alpha\in L^2_{\mathscr P}([t,T]\times\Omega;V)$, and we call $X^\alpha$ the solution of \eqref{SDE} with initial condition $X(t)=x$ and with control $\alpha$. A more precise notation would be $X^{t,x,\alpha}$, but we omit the dependence on $t$ and $x$ for better readability.

Let $(u,m)$ be the LQM mild solution of \eqref{mfg}. We want to compute the SDE satisfied by $u(\tau,X_\tau^\alpha)$, for $\tau\in [t,T]$. Since $u$ and $X^\alpha$ are just mild solutions of \eqref{mfg} and \eqref{SDE}, we argue as always by approximation.

We take $(A_n)_n\subset\mathcal L(H)$ as the Yosida approximants of $A$, and consider the related solution $X^{\alpha,n}$ of \eqref{SDE} and the related solution $(u_n,m_n)$ of \eqref{mfg}, with $A$ replaced by $A_n$. Since all the solutions are classical (see also Remark \ref{rem:casoregolare}), we can directly compute the equation satisfied by $u_n(\tau,X_\tau^{\alpha,n})$.

Observe that we have
$$
u_n(\tau,X_\tau^{\alpha,n})=\frac12\langle P_n(\tau)X_\tau^{\alpha,n},X_\tau^{\alpha,n}\rangle+\langle r_n(\tau),X_\tau^{\alpha,n}\rangle+s_n(\tau).
$$

Computing the SDE, and omitting for simplicity the dependence on $\tau$, we have
\begin{align*}
du_n(\tau,X_\tau^{\alpha,n})={}&\frac12\langle P'_n\xan,\xan\rangle d\tau+\langle P_n\xan, d\xan\rangle+\mathrm{tr}\big((\mathcal A+\mathcal A_0)P_n\big)d\tau\\
&+\langle dr_n,\xan\rangle+\langle r_n, d\xan\rangle+\mathrm{tr}\big(Z^{r,n}\sigma_0^*\big)d\tau+ds_n\\
={}&\frac12\langle\big[-P_n(A_n-D_1)-(A^*_n-D_1^*)P_n+P_n\mathcal B P_n-Q\big]\xan,\xan\rangle d\tau\\
&+\langle P_n\xan,A_n\xan+B\alpha+B_0 Y^n\rangle d\tau+\langle P_n\xan,\sigma dW_\tau\rangle+\langle P_n\xan,\sigma_0 dW^0_\tau\rangle\\
&+\Big[\mathrm{tr}\big((\mathcal A+\mathcal A_0)P_n+Z^{r,n}\sigma_0^*\big)-\langle (A^*_n-D_1^*-P_n\mathcal B)r_n+(P_nD_2+S)Y^n,\xan\rangle\Big]d\tau\\
&+\langle P_nc_3-\eta,\xan\rangle d\tau+\langle \xan, Z^{r,n}dW_\tau^0\rangle+\langle r_n,A_n\xan+B\alpha+B_0Y^n\rangle d\tau\\
&+\langle r_n,\sigma dW_\tau\rangle+\langle r_n,\sigma_0 dW_\tau^0\rangle-\mathrm{tr}\big((\mathcal A+\mathcal A_0)P_n+Z^{r,n}\sigma_0^*\big)d\tau+\frac12\langle\mathcal Br_n,r_n\rangle d\tau\\
&-\Big[\langle r_n,D_2 Y^n-c_3\rangle+\frac12\langle Z Y^n,Y^n\rangle+\langle \zeta,Y^n\rangle+\lambda\Big]d\tau+\langle Z^{s,n},dW_\tau^0\rangle.
\end{align*}

Rearranging, we get

\begin{align*}
du_n(\tau,X_\tau^{\alpha,n})={}&\Big[-\frac12\langle Q\xan,\xan\rangle-\langle SY^n,\xan\rangle-\frac12\langle ZY^n,Y^n\rangle-\langle\eta,\xan\rangle-\langle\zeta,Y^n\rangle-\lambda\Big]d\tau\\
&\!\!\!\!+\langle P_n\xan+r_n,-A_n\xan+D_1\xan-D_2 Y^n+c_3\rangle d\tau\\
&\!\!\!\!+\Big[\frac12\langle\mathcal B\big(P_n\xan+r_n\big),P_n\xan+r_n\rangle+\langle P_n\xan+r_n, A_n\xan+B\alpha+B_0 Y^n\rangle\Big]d\tau\\
&\!\!\!\!+\langle\sigma^*\big(P_n\xan+r_n\big),dW_\tau\rangle+\langle \sigma_0^*\big(P_n\xan+r_n\big)+(Z^{r,n})^*\xan+Z^{s,n},dW_\tau^0\rangle.
\end{align*}

Recall that, from the definition of $u_n$ and from \eqref{eq:v}, we have

$$
D_xu_n(t,x)=P_nx+r_n,\qquad v_n(t,x)=(Z^{r,n})^*x+Z^{s,n}.
$$

Using also the definition of the Hamiltonian and denoting by $\mathcal H_n$ the Hamiltonian obtained by replacing $A$ with $A_n$ in \eqref{eq:LewisHamilton}, we can compactly write the SDE for $u_n$ as
\begin{equation}\label{itowentzell}
\begin{split}
du_n(\tau,X_\tau^{\alpha,n})={}&\mathcal H_n(\tau,\xan_\tau,D_x u_n(\tau,\xan_\tau),Y_\tau^n)d\tau+\langle D_x u_n(\tau,\xan_\tau),d\xan_\tau\rangle\\
&+\langle v_n(\tau,\xan_\tau),dW_\tau^0\rangle.
\end{split}
\end{equation}

Integrating for $\tau\in[t,T]$, we get
\begin{align*}
u_n(T,\xan_T)={}&u_n(t,x)+\int_t^T\mathcal H_n(\tau,\xan_\tau,D_x u_n(\tau,\xan_\tau),Y_\tau^n) d\tau\\
&+\int_t^T\langle D_x u_n(\tau,\xan_\tau),b_n(\tau,\xan_\tau,Y^n_\tau,\alpha_\tau)\rangle d\tau+\int_t^T\langle\sigma^*(\tau) D_xu_n(\tau,\xan_\tau), dW_\tau\rangle\\
&+\int_t^T\langle\sigma_0^*(\tau)D_x u_n(\tau,\xan_\tau)+v_n(\tau,\xan_\tau), dW_\tau^0\rangle,
\end{align*}
where we have adopted the notation $b_n$ for the drift coefficient, i.e.
$$
b_n(\tau,x,y,\alpha):=A_nx+B(\tau)\alpha+B_0(\tau)y.
$$
Due to the definition of $\mathcal H_n$ given in \eqref{def:hamilton}, we have
\begin{equation}\label{daichecisiamoquasi}
\mathcal H_n(\tau,x,p,y)+\langle p,b_n(\tau,x,y,\alpha)\rangle\ge -L(\tau,\alpha,x,y),
\end{equation}
where $L$, defined in \eqref{eqn:Fqua}, does not depend on $n$. This implies, recalling that $u_n(T,\xan_T)=G(\xan_T,Y^n_T)$, again not depending on $n$ (see \eqref{eqn:Gqua}),
\begin{align*}
u_n(t,x)\le{}& \int_t^T L(\tau,\alpha_\tau,\xan_\tau,Y^n_\tau)d\tau + G(\xan_T,Y^n_T)\\
&-\int_t^T\langle\sigma^*(\tau)D_xu_n(\tau,\xan_\tau),dW_\tau\rangle-\int_t^T\langle\sigma_0^*(\tau)D_x u_n(\tau,\xan_\tau)+v_n(\tau,\xan_\tau),dW_\tau^0\rangle.
\end{align*}
Taking the conditional expectation with respect to $\mathscr F_t^0$, we get
\begin{align*}
u_n(t,x)\le \E\left[\int_t^T L(\tau,\alpha_\tau,\xan_\tau,Y^n_\tau)d\tau + G(\xan_T,Y^n_T)\mid \mathscr F_t^0\right].
\end{align*}
Now we pass to the limit for $n\to+\infty$. Observe that, thanks to Proposition \ref{prop:stabilitaforward}, we have
\begin{equation}\label{diciannovesoldiperunalira}
\xan\to X^\alpha,\qquad Y^n\to Y\qquad\mbox{in }\mathcal Q^2(H).
\end{equation}
Since $G$ and $L$ have a quadratic structure with bounded coefficients and $\alpha\in L^2$, we can say that, for a.e. $\tau\in[t,T]$,
$$
L(\tau,\alpha_\tau,X_\tau^{\alpha,n},Y_\tau^n)\to L(\tau,\alpha_\tau,X_\tau^\alpha,Y_\tau),\qquad G(\xan_T,Y_T^n)\to G(X^\alpha_T,Y_T)\qquad\mbox{in } L^1(\Omega;\R).
$$
Using \eqref{diciannovesoldiperunalira} and the Cauchy--Schwarz inequality we obtain
\[
\int_t^T L(\tau,\alpha_\tau,X_\tau^{\alpha,n},Y_\tau^n)\,d\tau
\to
\int_t^T L(\tau,\alpha_\tau,X_\tau^\alpha,Y_\tau)\,d\tau
\qquad\mbox{in }L^1(\Omega;\R).
\]
Indeed, the terms depending only on $\alpha$ cancel when taking the
difference, while the mixed terms involving $\alpha$ are controlled by its
$L^2_{\mathscr P}$ norm.

By the continuity of the conditional expectation in $L^1$,
\begin{align*}
&\E\Big[\int_t^T L(\tau,\alpha_\tau,\xan_\tau,Y^n_\tau)d\tau \mid \mathscr F_t^0\Big]\to \E\Big[\int_t^T L(\tau,\alpha_\tau,X^\alpha_\tau,Y_\tau)d\tau\mid \mathscr F_t^0\Big],\\
&\E\Big[G(\xan_T,Y^n_T)\mid\mathscr F_t^0\Big]\to \E\Big[G(X^\alpha_T,Y_T)\mid\mathscr F_t^0\Big]\qquad\mbox{in }L^1(\Omega;\R),
\end{align*}
where the convergence also holds a.s., up to a non-relabeled subsequence.

For the left-hand side, thanks to the quadratic structure of $u_n$ and the results contained in Proposition \ref{prop:stabilita}, Remark \ref{rem:generalizzazione} and Propositions \ref{prop:stabilitaforward}, \ref{stabilitaback}, we immediately have $u_n(\cdot,x)\to u(\cdot,x)$ in $\mathcal Q^1(\R)$, and, for the fixed $x\in H$, the convergence also holds a.s., up to a non-relabeled subsequence. This implies, passing to the limit,
\begin{equation}\label{eq:quasifinito}
u(t,x)\le \E\left[\int_t^T L(\tau,\alpha_\tau,X^\alpha_\tau,Y_\tau)d\tau + G(X^\alpha_T,Y_T)\mid \mathscr F_t^0\right].
\end{equation}
Since $\alpha$ is arbitrary, we get
$$
u(t,x)\le\inf\limits_{\alpha\in L^2_{\mathscr P}}\E\left[\int_t^T L(\tau,\alpha_\tau,X^\alpha_\tau,Y_\tau)d\tau + G(X^\alpha_T,Y_T)\mid \mathscr F_t^0\right]=\inf\limits_{\alpha\in L^2_{\mathscr P}}J(t,x,\alpha)\,.
$$
Finally, observe that equality holds in \eqref{eq:quasifinito} if it holds in \eqref{daichecisiamoquasi}, which happens, thanks to \eqref{eq:alfaottimo}, if
$$
\alpha^{opt}(\tau,x,p,y)=-R_0^{-1}(\tau)\big(B^*(\tau)p+L_0^*(\tau)x+M_0^*(\tau)y+\theta_0(\tau)\big).
$$
This gives the optimal control $\alpha^{opt,n}$ in feedback form:
\begin{equation}\label{infinito}
\alpha^{opt,n}_\tau=-R_0^{-1}(\tau)\Big(B^*(\tau)D_xu_n(\tau,X^{\alpha^{opt,n},n}_\tau)+L_0^*(\tau)X^{\alpha^{opt,n},n}_\tau+M_0^*(\tau)Y^n_\tau+\theta_0(\tau)\Big).
\end{equation}
In order to pass to the limit, we already know the convergence of $Y^n$ thanks to \eqref{diciannovesoldiperunalira}. We cannot use the same estimate for $X^{\alpha^{opt,n},n}$, since $\alpha^{opt,n}$ now depends on $n$.

For brevity, we write $X^{opt,n}$ instead of $X^{\alpha^{opt,n},n}$. The SDE satisfied by $X^{opt,n}$ is given by
\begin{equation*}
\begin{cases}
dX^{opt,n}_\tau=\big[A_nX_\tau^{opt,n}+B\alpha^{opt,n}_\tau+B_0Y^n_\tau\big]\,d\tau+\sigma(\tau) dW_\tau+\sigma_0(\tau) dW_\tau^0\,,\\
X^{opt,n}_t=x,
\end{cases}
\end{equation*}
which means, using \eqref{infinito} and \eqref{EQUAZIONI},
\begin{equation*}
\begin{cases}
dX^{opt,n}_\tau=\big[A_nX_\tau^{opt,n}-(D_1+\mathcal BP_n)X_\tau^{opt,n}+D_2Y_\tau^n-\mathcal Br_n-c_3]\,d\tau+\sigma(\tau) dW_\tau+\sigma_0(\tau) dW_\tau^0\,,\\
X^{opt,n}_t=x.
\end{cases}
\end{equation*}

Then, thanks to \eqref{madonnabenedettadellincoronata} and Proposition \ref{prop:stabilitaforward}, we have $X^{opt,n}\to X^{opt}$ in $\mathcal Q^2(H)$, where $X^{opt}$ solves \eqref{ottima}. Moreover, thanks to \eqref{boundforward}, $X^{opt,n}$ is bounded in $\mathcal Q^2(H)$ uniformly in $n$.

\noindent Moreover, observe that we also have
\begin{equation}\label{ciao}
\begin{split}
&\norm{D_xu_n(\cdot, X^{opt,n}_\cdot)-D_xu(\cdot,X^{opt}_\cdot)}_{\mathcal Q^2(H)}=\norm{P_nX^{opt,n}-PX^{opt}+r_n-r}_{\mathcal Q^2(H)}\\
&\qquad\qquad\le C\norm{P_nX^{opt,n}-PX^{opt}}_{\mathcal Q^2(H)}+C\norm{r_n-r}_{\mathcal Q^2(H)}.
\end{split}
\end{equation}
We write
$$
P_nX^{opt, n}-PX^{opt}=P_n(X^{opt,n}-X^{opt})+(P_n-P)X^{opt}.$$
The first term goes to zero since $P_n$ is uniformly bounded and $X^{opt,n}$ converges to $X^{opt}$ in $\mathcal Q^2(H)$. For the second, for almost every $\omega$, write 
$$
K_\omega= \{X_\tau^{opt}(\omega)\mid \tau\in [t,T]\}
$$
which is compact in $H$ since the trajectory is continuous. By Lemma \ref{ascoli}, we get
$$
\sup_{\tau\in[t,T]}\|(P_n(\tau)-P(\tau))X_\tau^{opt}(\omega)\|_H\to0.
$$
By uniform boundedness, the left-hand side is dominated by a random variable in $L^2(\Omega)$, and we can apply the dominated convergence theorem, concluding 
\begin{equation}\label{ciao2}
\norm{P_nX^{opt,n}-PX^{opt}}_{\mathcal Q^2(H)} \to 0.
\end{equation}
Then by \eqref{ciao}, \eqref{ciao2} and \eqref{madonnabenedettadellincoronata} we get
$$
\norm{D_xu_n(\cdot, X^{opt,n}_\cdot)-D_xu(\cdot,X^{opt}_\cdot)}_{\mathcal Q^2(H)}\to 0.
$$
This implies $\alpha_\cdot^{opt,n}\to\alpha_\cdot^{opt}$ in $\mathcal Q^2(V)$, where $\alpha^{opt}$ is given by \eqref{eq:opt}.

For every $n$, the choice $\alpha^{opt,n}$ yields equality in
\eqref{daichecisiamoquasi}. Hence,
\[
u_n(t,x)=
\E\left[
\int_t^T L(\tau,\alpha_\tau^{opt,n},X_\tau^{opt,n},Y_\tau^n)\,d\tau
+G(X_T^{opt,n},Y_T^n)
\mid\mathscr F_t^0
\right].
\]
Using the previous convergences and arguing as above, we can pass to the
limit and obtain
\[
u(t,x)=J(t,x,\alpha^{opt}),
\]
which concludes the proof.
\end{proof}

\begin{remark}
We want to stress the fact that, at least in the finite dimensional cases, equation \eqref{itowentzell} can be obtained in an easier way using the so-called \emph{Itô-Wentzell formula}: if $u(\cdot,x)$ is a process solving the SDE
$$
du(t,x)=F(t,x)dt+\langle G(t,x),dW_t\rangle,
$$
and $X_\cdot$ is a process satisfying
$$
dX_t=b_t dt+\Sigma_t dW_t,
$$
then it holds
$$
d(u(t,X_t))=(du)(t,X_t)+\langle D_xu(t,X_t),dX_t\rangle+\frac12\mathrm{tr}(\Sigma_t\Sigma_t^*D^2_{xx}u(t,X_t))dt+\mathrm{tr}(\mathrm{Jac}G(t,X_t)\Sigma_t)dt.
$$
We decided to adopt a direct proof of \eqref{itowentzell}, using the quadratic structure of $u_n$, since we are not aware of any generalization of the previous formula for SDEs in infinite dimensional Hilbert spaces.
\end{remark}

\section{Acknowledgments}
The authors are supported by the 2026 GNAMPA project “Problemi di controllo e modelli
mean-field con vincoli di stato e dimensione infinita: teoria e applicazioni”, and have been part of the 2022 PRIN project ``Impact of the Human Activities on the Environment and Economic Decision Making in a Heterogeneous Setting: Mathematical Models and Policy Implications".

The first author has been supported by the 2025 GNAMPA project ``Modelli economici ambientali via MFGs sui grafi e grafoni: studio analitico e assimilazione dei dati", and she is part of the project ``Transition-Field Mean Field Approaches for Environmental Transition".

\printbibliography

\end{document}